\documentclass[12pt,a4paper]{article}

\usepackage{amssymb}

\usepackage{amsmath}

\newtheorem{theorem}{Theorem}[section]
\newtheorem{definition}[theorem]{Definition}
\newtheorem{lemma}[theorem]{Lemma}
\newtheorem{corollary}[theorem]{Corollary}
\newtheorem{proposition}[theorem]{Proposition}

\newtheorem{theorem(Composition-Diamond Lemma)}[theorem]{Theorem(Composition-Diamond Lemma)}

\begin{document}

\title{Composition-Diamond lemma for associative $n$-conformal algebras\footnote{Supported
by the NNSF of China (No.10771077) and the
NSF of Guangdong Province (No.06025062).}}

\author{
L. A. Bokut\footnote{Supported by RFBR 01-09-00157, LSS--344.2008.1
and SB RAS Integration grant No. 2009.97 (Russia).}\\
{\small \ School of Mathematical Sciences, South China Normal
University}\\
{\small Guangzhou 510631, P. R. China}\\
{\small Sobolev Institute of Mathematics, Russian Academy of
Sciences}\\
{\small Siberian Branch, Novosibirsk 630090, Russia}\\
{\small  bokut@math.nsc.ru}\\
\\
 Yuqun
Chen\\
{\small \ School of Mathematical Sciences, South China Normal
University}\\
{\small Guangzhou 510631, P. R. China}\\
{\small yqchen@scnu.edu.cn}\\
\\
Guangliang Zhang\\
 {\small \ Nanhai College, South China Normal
University}\\
{\small Guangzhou 510631, P. R. China}\\
{\small  zgl541@tom.com}}

\date{}
\maketitle

\maketitle \noindent\textbf{Abstract:} In this paper, we study the
concept of associative $n$-conformal algebra over a field of
characteristic 0 and establish Composition-Diamond lemma for a free
associative $n$-conformal algebra. As an application, we construct
Gr\"{o}bner-Shirshov bases for Lie $n$-conformal algebras presented
by generators and defining relations.


\noindent {\bf AMS} Mathematics Subject Classification(2000): 17B69,
16S15, 13P10

\section{Introduction}

 Gr\"{o}bner  and  Gr\"{o}bner--Shirshov bases theories were invented independently by A.I. Shirshov
\cite{Sh62c} for Lie algebras and H. Hironaka \cite{Hi64} and B.
Buchberger \cite{Bu65,Bu70} for associative-commutative algebras.

Shirshov's paper \cite{Sh62c} based on his papers \cite{Sh62b}
(Gr\"{o}bner--Shirshov bases theory for (anti) commutative algebras,
the reduction algorithm for (anti-) commutative algebras) and
\cite{Sh58} (Lyndon--Shirshov words (these words were defined some
earlier \cite{Lyndon54}, but incidentally that was unknown for 25
years in Russia and these words were called Shirshov's regular
words, see, for example, \cite{Bah87,Bah03,BMPZ92,Bo72,Ufn95,Mik96},
see also \cite{Cohn65}), Lyndon--Shirshov basis of a free Lie
algebra (see also \cite{CFL58}). The latter Shirshov's papers
\cite{Sh58,Sh62b} were based on his Thesis \cite{Sh53a}, A.G. Kurosh
(adv), published in three papers \cite{Sh53b} (on free Lie algebras:
$K_d$-Lemma (Lazard--Shirshov elimination process), the subalgebra
theorem (Shirshov--Witt theorem)), \cite{Sh54} (on free (anti-)
commutative algebras: linear bases, the subalgebra theorems),
\cite{Sh62a} (on free Lie algebras: a series of bases with
well-ordering of basic Lie words such that $[w]=[[u][v]]>[v]$, see
also \cite{Vie78}; the series is called now Hall sets \cite{Reu93}
or Hall--Shirshov bases). Shirshov's Thesis, in turn, was in line
with a Kurosh's paper \cite{Ku60} (on free non-associative algebras:
the subalgebra theorem). Also Shirshov's paper \cite{Sh62b} was in a
sense of a continuation of a paper by A.I. Zhukov \cite{Zhu50}, a
student of Kurosh (on free non-commutative algebras: decidability of
the word problem for non-associative algebras). The difference with
the Zhukov's approach was that Zhukov did not use any linear
ordering of non-associative words, but just the partial deg-ordering
to compere two words by the degree (length).

It would be not a big exaggeration to say that Shirshov's paper
\cite{Sh62c} was between line of the Kurosh's program of study free
algebras of different classes of non-associative algebras.

Shirshov's paper \cite{Sh62c} contained implicitly the
Gr\"{o}bner--Shirshov bases theory for associative algebras too
because he constantly used that any Lie polynomial is at the same
time a non-commutative polynomial. For example, the maximal term of
a Lie polynomial is defined as its maximal word as a non-commutative
polynomial, definition of a Lie composition (Lie S-polynomial) of
two Lie polynomials begins with the definition of their composition
as non-commutative polynomials and follows by putting some special
Lie brackets on it, and so on. The main Composition (-Diamond) Lemma
for associative polynomials is actually proved in the paper and we
need only to ``forget" about Lie brackets in the proof of this lemma
for Lie polynomials (\cite{Sh62c} Lemma 3). Explicitly Composition
(-Diamond) Lemma was formulated much later in papers L.A. Bokut
\cite{Bo76} and G. Bergman \cite{Be78}.

We formulate Shirshov's Composition-Diamond Lemma for associative
algebras following his paper \cite{Sh62c} with the only change of
``Lie polynomials" to ``non-commutative polynomials". Let $k\langle
X\rangle$ be a free associative algebra over a field $k$ on a set
$X$ such that the free monoid $X^*$ is well-ordered. For a
polynomial $f$, by $\overline f$, Shirshov \cite{Sh62c} denotes the
maximal word of $f$. Let $f,g$ are two monic polynomials (may be
equal), $w\in X^*$ such that $w=acb,\ \overline f=ac,\ \overline
g=cd$, where $a,b,c$ are words and $c$ is not empty. Then
$(f,g)_w=fb-ag$ is called an (associative) composition of $f,g$ (it
is the original Shirshov's notation, now we use $(f,g)_w)$, clearly
$\overline {fb-ag} <w$. For Lie polynomials $f,g$, Shirshov puts
some special brackets $[fb]-[ag]$ such that $\overline {[fb]-[ag]}
<w$ too. Let $S$ be a reduced set in $k\langle X\rangle$ and $S^*$
is a reduced set that is obtained from $S$ by (transfinite)
induction applying the following elementary operation:  to join to
$S$ a composition of two elements of $S$ and to apply the reduction
algorithm to the result set (until one will get a reduced set with
only trivial compositions, after the reduction). By nowadays
terminology, $S^*$ is a Gr\"{o}bner--Shirshov basis of the ideal
generated by S, and that process of adding compositions is the
Shirshov's algorithm. He calls $S$ a $stable$ set if, at each step,
the degree of the composition $(f,g)_w$, after the reduction, is
bigger than degrees of $f,g$ (or $(f,g)_w$ is zero after the
reduction). Of course, if S is a finite (or recursive)  stable set,
then $S^*$ is a recursive set, and from the next lemma the word
problem is solvable in the algebra with defining relations $S$. Now
suppose that $S$ is a Gr\"{o}bner--Shirshov basis in the sense that
$S$ is a reduced set and any composition of elements of $S$ is
trivial after the reduction ($S$ is closed under compositions). The
$S^*$ is fulfilled this condition (simply because any composition,
after the reduction, is in $S^*$). Hence S is a stable set in the
sense of Shirshov. Then Lemma 3 in \cite{Sh62c} has the following
form.

\ \

\noindent\textbf{Shirshov's Composition-Diamond Lemma for
associative algebras}. Let $S\subset k\langle X\rangle$ is a
Gr\"{o}bner--Shirshov basis of the ideal $Id(S)$. Let $f\in Id(S)$.
Then $\overline f=a\overline s b$ for some $s\in S, a,b\in X^*$
(hence $S$-irreducible words $Irr(S)$, that does not contain
subwords which are maximal words of polynomials from $S$, is a
$k$-basis of the algebra $k\langle X|S\rangle$ with defining
relations $S$).

It is easy to see that converse is also true.

\ \

Last years there were quite a few results on Gr\"{o}bner--Shirshov
bases for semi-simple Lie (super) algebras, irreducible modules,
Kac--Moody algebras, Coxeter groups, braid groups, quantum groups,
conformal algebras, free inverse semigroups, Kurosh's
$\Omega$-algebras, Loday's dialgebras, Leibniz algebras, Rota-Baxter
algebras, Vinberg-Milnor's right-symmetric algebras, and so on, see,
for example, papers
\cite{Bo07,Bo09,BChainikovShum07,BFKS08,BoKLM99,BokutShiao01,KL00a,KL00b,KLLO02,KLLP07},
surveys \cite{BokutChen08,BChibrikov05,BoFKK00,BK00, BK05a}.
Actually, conformal algebras, dialgebras, Rota-Baxter algebras are
examples of $\Omega$-algebras. For non-associative
$\Omega$-algebras, Composition-Diamond Lemma was proved in
\cite{DH08}. The case of associative $\Omega$-algebras (associative
algebras with any set $\Omega$ of multi-linear operations) was
treated in \cite{BCQ08} with an application to free ($\lambda$-)
Rota-Baxter algebras (the latter is associative algebras with linear
operation $P(x)$ and the identity $P(x)P(y)=P(P(x)y) + P(xP(y))
+\lambda P(xy)$, where $\lambda$ is a fix element of a ground field,
see, for example, \cite{EGuo08}). Composition-Diamond Lemma for
dialgebras \cite{BCL08} has an application to the PBW theorem for
universal enveloping dialgebras of  Leibniz algebras  (see
\cite{Aymon03}).

This paper is a continuation of the paper \cite{BFK04} on
Gr\"{o}bner--Shirshov bases for conformal algebras. V. Kac
\cite{Kac96} defined conformal algebras under the influence of
Vertex algebras (Belavin-Polyakov-Zamolodchikov \cite{BPZZ84},
Borcherds \cite{Bor86}) and the Operator Product Expansion in
mathematical physics. Structure theory of associative and Lie
conformal algebras were studied in
\cite{BAK01,Ko91,Ko00a,Ko00b,Ko02, Ko03,Ko07a,Ko07b,Ko08a,Ko08b}. M.
Roitman \cite{Ro99,Ro00,Ro05} instigated free associative and Lie
conformal algebras. In particular, he found a linear basis of a free
associative conformal algebra and proved that an analog of PBW
theorem is not valid for Lie conformal algebras.

 Composition-Diamond Lemma for associative conformal algebras \cite{BFK04}
   has some difference from that lemma for
associative and Lie algebras (as well as from Buchberger Theorem for
commutative algebras). First of all, there are compositions of left
and right multiplications, just opposite to the classical cases.
Secondly, the condition on a set of relations $S$, that $S$ is
closed under compositions, is not equivalent to the condition that
$Irr(S)$ is a basis of the algebra with defining relations $S$.
Later on, it became clear that the same situation with both
conditions is valid for dialgebras and right-symmetric algebras, as
well as for $n$-conformal algebras.

The class of $n$-conformal algebras is a particular example of the
class of (H-) pseudo-algebras in the sense of Bakalov--D'Andrea--Kac
\cite{BAK01} and Beilinson--Drinfeld \cite{BD04} for the case the
Hopf algebra H is a polynomial algebra on $n$ variables.

In this paper, we state Composition-Diamond Lemma for associative
$n$-conformal algebras. A proof follows from the proof of
Composition-Diamond Lemma for associative conformal algebras in
\cite{BFK04} almost word by word.

\begin{definition}
 Let $k$ be a field with characteristic $0$ and $C$ a vector
 space over $k$. Let $Z_+$ be the non-negative integer number, $Z$ the integer ring and $n$
 a positive integer number.
 We associate to each $(m_1,\cdots,m_n) \in Z_+^n$
 a bilinear product on $C$, denoted by $\langle
 \overrightarrow{m}\rangle$ where
$\overrightarrow{m}=(m_1,\cdots,m_n)\in Z_+^n$. Let
 $D_i:C\rightarrow C$ be linear mappings such that  $D_iD_j=D_jD_i, \ 1\leqslant i,j \leqslant n$.
 Then $C$ is an $n$-conformal algebra with derivations $D=\{D_1,\dots, D_n\}$ if the following axioms are
 satisfied:
\begin{enumerate}
\item[(i)]\ If $a,b \in C$, then there is an $\overrightarrow{N}(a,b)\in
Z_+^n$ such that
 $a \langle \overrightarrow{m}\rangle b=0$ if $\overrightarrow{m}\not\prec
 \overrightarrow{N}(a,b)$, where for any $\overrightarrow{m}=(m_1,\cdots,m_n),
\overrightarrow{l}=(l_1,\cdots,l_n)\in Z_+^n$,
\begin{eqnarray*}
&&\overrightarrow{m}=(m_1,\cdots,m_n)\prec
\overrightarrow{l}=(l_1,\cdots,l_n)\Leftrightarrow m_i\leq l_i, \
i=1,\dots,n \\
&& \mbox{and there exists }\ \ i_0, \ 1\leq i_0\leq n \ \ \mbox{
such that }\ \ \ m_{i_0}< l_{i_0}.
\end{eqnarray*}
 $\overrightarrow{N}:\ C\times C\rightarrow Z_+^n$ is called the locality
 function.
\item[(ii)]\ For any $a,b \in C, \ \overrightarrow{m} \in Z_+^n$,
      $D_i(a \langle \overrightarrow{m}\rangle b)=D_ia \langle \overrightarrow{m}\rangle b+
      a \langle \overrightarrow{m}\rangle D_ib,\ i=1,\dots,n.$
\item[(iii)]\ For any $a,b \in C, \ \overrightarrow{m} \in Z_+^n$,
    $D_ia \langle \overrightarrow{m}\rangle b=-m_ia \langle \overrightarrow{m}-\overrightarrow{e_i}\rangle
    b,\  i=1,\dots,n$, where
    $\overrightarrow{e_i}=(\underbrace{0,\cdots,0}_{i-1},1,0,\cdots,0)$.
\end{enumerate}
\end{definition}

For $\overrightarrow{m},\overrightarrow{s}\in Z_+^n$, put
$(-1)^{\overrightarrow{s}}=(-1)^{s_1+\cdots+s_n}$ and
$\binom{\overrightarrow{m}}{\overrightarrow{s}}=\binom{m_1}{s_1}\cdots\binom{m_n}{s_n}$.

\ \

An $n$-conformal algebra $C$ is associative if in addition the
associativity condition holds.

{\bf Associativity Condition:} For  any $a,b,c \in C$ and
$\overrightarrow{m},\overrightarrow{m'} \in Z_+^n$,
$$
   (a\langle \overrightarrow{m}\rangle b)\langle \overrightarrow{m}'\rangle
   c=\sum_{\overrightarrow{s}\in
   Z_+^n}(-1)^{\overrightarrow{s}}\binom{\overrightarrow{m}}{\overrightarrow{s}}a\langle
   \overrightarrow{m}-\overrightarrow{s}\rangle (b\langle
   \overrightarrow{m}'+\overrightarrow{s} \rangle c),
   $$
   or equivalently a $right$ analogy
   $$
   a\langle \overrightarrow{m}\rangle (b\langle \overrightarrow{m}'\rangle
   c)=\sum_{\overrightarrow{s}\in
   Z_+^n}(-1)^{\overrightarrow{s}}\binom{\overrightarrow{m}}{\overrightarrow{s}}(a\langle
   \overrightarrow{m}-\overrightarrow{s}\rangle b)\langle
   \overrightarrow{m}'+\overrightarrow{s}\rangle c.
   $$

Let $C(B,\overrightarrow{N},D_1,\cdots,D_n)$ be an associative
$n$-conformal algebra generated by $B$ with the locality function
$\overrightarrow{N}:\ B\times B\rightarrow Z_+^n$. Then
$C(B,\overrightarrow{N},D_1,\cdots,D_n)$ is called the free
associative $n$-conformal algebra with the data
$(B,\overrightarrow{N},D_1,\cdots,D_n)$ if for any associative
$n$-conformal algebra $C'$ with derivations $D=\{D_1,\dots, D_n\}$
and any map $\varepsilon: B\rightarrow C'$ such that
$\varepsilon(a)\langle \overrightarrow{m}\rangle \varepsilon(b)=0$
if $\overrightarrow{m}\not\prec\overrightarrow{N}(a,b)$, there
exists a unique homomorphism
$f:C(B,\overrightarrow{N},D_1,\cdots,D_n)\rightarrow C'$ such that
$f(b)=\varepsilon(b)$ for all $b\in B$.

\section{Free associative $n$-conformal algebra}

\begin{definition}
Let $B$ be a non-empty set and $D=\{D_1,\dots, D_n\}$. Denote by
$D^\omega(B)=\{D_1^{i_1}\cdots D_n^{i_n}b\ | \ b\in
B,i_1,\cdots,i_n\in Z_+\}$. Then, we define the lengths of words
(non-associative) on $D^\omega(B)$ inductively:
\begin{enumerate}
\item[(i)]\ any $D_1^{i_1}\cdots D_n^{i_n}b,\ b\in B,\ i_1,\cdots,i_n\in Z_+ $
is a word of length 1,
\item[(ii)]\ if $(u)$ and $(v)$ are two words of lengths $k$ and $l$
respectively, then $ ((u) \langle \overrightarrow{m}\rangle (v))$ is
a word of length $k+l$.
\end{enumerate}
For any word $u$ on $D^\omega(B)$, we denote its length by $|u|$.
\end{definition}

A word $u$ is $right$ $normed$ if it has the form
\begin{equation}\label{e1.1}
[u]=x_1\langle \overrightarrow{m}^{(1)}\rangle(x_2\langle
\overrightarrow{m}^{(2)}\rangle\cdots(x_k\langle
\overrightarrow{m}^{(t)}\rangle x_{t+1})\cdots)
\end{equation}
where $x_i\in D^\omega(B), \overrightarrow{m}^{(j)}\in
Z_+^n,1\leqslant i \leqslant t+1,1\leqslant j \leqslant t$.

\ \

We fix a locality function $\overrightarrow{N}:\ B\times
B\rightarrow Z_+^n$.

\begin{definition}
A right normed word $[u]$ of the form (\ref{e1.1}) is called a
$normal$ $word$ if
$$
x_j\in B, \ \overrightarrow{m}^{(j)} \prec
\overrightarrow{N}(a_j,a_{j+1}) \ \mbox{ for } \ j=1,\dots, t,
 \ \mbox{ and }\   x_{t+1}=D_1^{i_1}\cdots D_n^{i_n}b\in D^\omega(B).
$$
If this is the case, then we denote the index of $[u]$ by
$ind(u)=(i_1,\cdots,i_n)$.
\end{definition}

Here and after, $[u]$ would mean the right normed bracketing
\begin{eqnarray*}
&&[x_1\langle \overrightarrow{m}^{(1)}\rangle x_2\langle
\overrightarrow{m}^{(2)}\rangle\cdots x_t\langle
\overrightarrow{m}^{(t)}\rangle x_{t+1}]\\
&=&x_1\langle \overrightarrow{m}^{(1)}\rangle (x_2\langle
\overrightarrow{m}^{(2)}\rangle\cdots (x_t\langle
\overrightarrow{m}^{(t)}\rangle x_{t+1})\cdots).
\end{eqnarray*}

\ \

In the following Lemmas \ref{l2.1} and \ref{l2.2}, we assume that
$C_1$ is an associative $n$-conformal algebra generated by all
non-associative words on $D^\omega(B)$.

We denote by $(u)$ a nonassociative word on $u$, i.e., $(u)$ is some
bracketing of $u$.

\begin{lemma}\label{l2.1}
Any word $(u)$ on $D^\omega(B)$ is a linear combination of right
normed words of the same length $|u|$.
\end{lemma}
\noindent{\bf Proof.} We use induction on $|u|$. The result holds
trivially for $|u|=1$. Let $|u|>1$. We assume that
\begin{equation}\label{e2.1}
(u)=[v]\langle \overrightarrow{m} \rangle[w]
\end{equation}
where $[v]$ and $[w]$ are right normed words. If $|v|=1$, then
(\ref{e2.1}) is right normed, and we are done. Hence, let $|v|>1$.
Then
$$
[v]=x\langle \overrightarrow{m}'\rangle[v_1],\ for \ some \  x\in
D^\omega(B),
$$
and so
$$
([v]\langle \overrightarrow{m}\rangle[w])=(x\langle
\overrightarrow{m}'\rangle[v_1])\langle \overrightarrow{m}\rangle[w]
=\sum_{\overrightarrow{s}\in
   Z_+^n}(-1)^{\overrightarrow{s}}\binom{\overrightarrow{m}}{\overrightarrow{s}}x\langle
   \overrightarrow{m}'-\overrightarrow{s}\rangle ([v_1]\langle
   \overrightarrow{m}+\overrightarrow{s}\rangle[w]).
$$
Now, the result follows the induction. \ \  $\square$

\begin{lemma}\label{l2.2}
Any word $(u)$ on $D^\omega(B)$ is a linear combination of normal
words of the same length $|u|$.
\end{lemma}
\noindent{\bf Proof.} Due to Lemma \ref {l2.1}, we may assume that
$(u)$ is right normed, i.e., $(u)=[u]$. We proceed by induction on
$|u|$. If $|u|=2$, then
$$
[u]=D_1^{i_1}\cdots D_n^{i_n}a_1\langle \overrightarrow{m}\rangle
D_1^{j_1}\cdots D_n^{j_n}a_2=\alpha_1\cdots\alpha_n a_1\langle
\overrightarrow{m}-\overrightarrow{i}\rangle D_1^{j_1}\cdots
D_n^{j_n}a_2
$$
where
$\alpha_s=(-1)^{i_s}[m_s]^{i_s},[m_s]^{i_s}=m_s(m_s-1)\cdots(m_s-i_s+1),1\leqslant
s \leqslant n$ and the right side is zero if $m_t-i_t<0$ for some
$t$. Thus, in this case, it suffices to deal with a word
\begin{equation}\label{e2.2}
[u]=a_1\langle \overrightarrow{m}\rangle D_1^{j_1}\cdots
D_n^{j_n}a_2, \ a_1,a_2\in B, \ j_1,\cdots,j_n\in Z_+, \
\overrightarrow{m}\in Z_+^n
\end{equation}
If $\overrightarrow{m}\prec
\overrightarrow{N}=\overrightarrow{N}(a_1,a_2)$, then (\ref{e2.2})
is normal, and we are done. We assume that
$\overrightarrow{N}(a_1,a_2)=(l_1,\cdots,l_n)$ and $m_t\geqslant
l_t$ for some $t$. If $j_t=0$, then
\begin{eqnarray*}
&&[u]=a_1\langle \overrightarrow{m}\rangle D_1^{j_1}\cdots
D_{t-1}^{j_{t-1}}D_{t+1}^{j_{t+1}}\cdots D_n^{j_n}a_2\\
&=&\beta
\sum_{p_1,\cdots,p_{t-1},p_{t+1},\cdots,p_n\in
   Z_+}\binom{m_1}{p_1}\cdots\binom{m_n}{p_n}D_1^{j_1-p_1}\cdots D_n^{j_n-p_n}(a_1\langle
   \overrightarrow{m}-p_1\overrightarrow{e_1}\cdots-p_n\overrightarrow{e_n}\rangle a_2 )\\
&=&0
\end{eqnarray*}
where $\beta=[j_1]^{p_1}\cdots[j_{t-1}]^{p_{t-1}}[j_{t+1}]^{p_{t+1}}\cdots[j_n]^{p_n}$.\\
If $j_t>0$, then
\begin{eqnarray*}
0&=&D_t^{j_t}(a_1\langle \overrightarrow{m}\rangle D_1^{j_1}\cdots
D_{t-1}^{j_{t-1}}D_{t+1}^{j_{t+1}}\cdots D_n^{j_n}a_2)\\
&=&[u]+\sum_{p>0}(\gamma_p a_1\langle
   \overrightarrow{m}-p\overrightarrow{e_t}\rangle D_t^{j_t-p}D_1^{j_1}\cdots
D_{t-1}^{j_{t-1}}D_{t+1}^{j_{t+1}}\cdots D_n^{j_n}a_2)
\end{eqnarray*}
where $\gamma_p=(-1)^p[m_t]^p\binom{j_t}{p}$. By induction on $j_t$
we conclude the case for $|u|=2$.

 Now, let $|u|\geqslant3$ and
suppose $[u]$ has the form (\ref{e1.1}). Using axioms of associative
$n$-conformal algebra and induction, we may assume that
$$
 [u]=a_1\langle
\overrightarrow{m}\rangle[u_1]
$$
where $[u_1]=a_2\langle \overrightarrow{m}'\rangle[u_2]$ is a normal
word of length $|u|-1$.

 Now, if
$\overrightarrow{m}<\overrightarrow{N}=\overrightarrow{N}(a_1,a_2)$,
then $[u]$ is already normal. Otherwise
$$
[u]=a_1\langle \overrightarrow{m}\rangle(a_2\langle
\overrightarrow{m}'\rangle[u_2]) =-\sum_{\overrightarrow{s}\in
   Z_+^n \backslash 0}(-1)^{\overrightarrow{s}}\binom{\overrightarrow{m}}{\overrightarrow{s}}a_1\langle
   \overrightarrow{m}-\overrightarrow{s}\rangle (a_2\langle
   \overrightarrow{m}'+\overrightarrow{s}\rangle [u_2]).
$$
By induction on $\overrightarrow{m}$, we can conclude that each
summand in right hand side above is normal. \ \  $\square$

\ \

\noindent{\bf Remark}. The proof of the Lemma \ref{l2.2} provides an
algorithm for presenting any given word $(u)$ as a linear
combination of normal words. Here is the algorithm:
\begin{enumerate}
\item [(i)]\
If $|u|=1$, then $(u)$ is normal.
\item [(ii)]\ If $|u|=2$, then $u=a\langle \overrightarrow{m}\rangle D_1^{i_1}\cdots
D_n^{i_n}b$ for some $a,b\in B,\overrightarrow{m}\in Z_+^n$ and
$i_1,\cdots,i_n\in Z_+$. If $i_1=\cdots=i_n=0$, $u$ is normal if
$\overrightarrow{m}\prec  \overrightarrow{N}(a,b)$, or zero if
$\overrightarrow{m}\not\prec  \overrightarrow{N}(a,b)$. Otherwise
assume that $i_t>0, 1\leqslant t \leqslant n$. Then
\begin{eqnarray*}
u&=&a\langle \overrightarrow{m}\rangle D_1^{i_1}\cdots D_n^{i_n}b\\
&=&D_t(a\langle \overrightarrow{m}\rangle D_1^{i_1}\cdots
D_{t}^{i_{t}-1}\cdots D_n^{i_n}b)+ m_t a\langle
\overrightarrow{m}-\overrightarrow{e_t}\rangle D_1^{i_1}\cdots
D_{t}^{i_{t}-1}\cdots D_n^{i_n}b.
\end{eqnarray*}
Now use induction
to rewrite both terms as a linear combination of normal words.
\item [(iii)]\ Suppose $|u|\geqslant3$. Then by using the arguments from the proof of
Lemma \ref{l2.1}, one can present $(u)$ in the form
$$
(u)=a_1\langle \overrightarrow{m}\rangle(a_2\langle
\overrightarrow{m}'\rangle(v))
$$
where $a_1,a_2\in B$ and $\overrightarrow{m},\overrightarrow{m}'\in
Z_+^n$. If $\overrightarrow{m}\prec \overrightarrow{N}(a_1,a_2)$,
apply induction to rewrite $a_2\langle
\overrightarrow{m}'\rangle(v)$ as a linear combination of normal
words of length $|u|-1$, and then multiply each term from the left
side by $a_1\langle \overrightarrow{m}\rangle$ to obtain the desired
linear combinations. If $\overrightarrow{m}\not\prec
\overrightarrow{N}(a_1,a_2)$, then
$$
(a_1\langle \overrightarrow{m}\rangle a_2)\langle
\overrightarrow{m}'\rangle(v)=0
$$
and so
$$
a_1\langle \overrightarrow{m}\rangle(a_2\langle
\overrightarrow{m}'\rangle(v)) =-\sum_{\overrightarrow{s}\in
   Z_+^n\backslash 0}(-1)^{\overrightarrow{s}}\binom{\overrightarrow{m}}{\overrightarrow{s}}a_1\langle
   \overrightarrow{m}-\overrightarrow{s}\rangle (a_2\langle
   \overrightarrow{m}'+\overrightarrow{s}\rangle (v)).
$$
Now, for $\overrightarrow{s}\in
   Z_+^n\backslash 0$, $\overrightarrow{m}-\overrightarrow{s}\prec  \overrightarrow{m}$ and
   the induction hypothesis applies to each term $a_1\langle
   \overrightarrow{m}-\overrightarrow{s}\rangle (a_2\langle
   \overrightarrow{m}'+\overrightarrow{s}\rangle (v))$ to obtain linear
   combination of words of the form $a_1\langle
   \overrightarrow{l}\rangle (a_2\langle
   \overrightarrow{l}'\rangle (v))$ with $\overrightarrow{l}\prec  \overrightarrow{N}(a_1,a_2)$.
   Another treatment of these words finishes the algorithm.
\end{enumerate}

Let $C$ be a $k$-linear space spanned by normal words (as a basis).
For normal words $[u]$ and $[v]$ and any $\overrightarrow{m}\in
Z_+^n$, the multiplication $[u]\langle \overrightarrow{m}\rangle[v]$
is defined using the algorithm. Also, define $D_i[u]$ by Leibniz
rule and the rule
$$
D_ia\langle \overrightarrow{m}\rangle[u]=-m_ia\langle
\overrightarrow{m}-\overrightarrow{e_i}\rangle[u],\ for \ a\in B,\
1\leqslant i \leqslant n.
$$

\begin{theorem}
$C=C(B,\overrightarrow{N},D_1,\cdots,D_n)$ is a free associative
$n$-conformal algebra with the data
$(B,\overrightarrow{N},D_1,\cdots,D_n)$.
\end{theorem}
\noindent{\bf Proof.} It is sufficient to check the axioms for
associative $n$-conformal algebras.

(i) (Locality): $[u]\langle \overrightarrow{m}\rangle[v]=0$ if
$\overrightarrow{m}\not\prec \overrightarrow{N}(u,v)$. To see this,
we use induction on $|u|$.
\par Assume first that $|u|=1$ and $u=a\in B$. Now, we want to use
a second induction on $|v|$. If $[v]=b\in B$, then
$\overrightarrow{N}(u,v)=\overrightarrow{N}(a,b)$ by the definition,
and so $[u]\langle \overrightarrow{m}\rangle[v]=a\langle
\overrightarrow{m}\rangle b=0$ if $\overrightarrow{m}\not\prec
\overrightarrow{N}(u,v)$. If $[v]=D_1^{i_1}\cdots D_n^{i_n}b$ where
$i_t>0$ for some $t,\ 1\leqslant t \leqslant n$, then by Remark, we
have
$$
a\langle \overrightarrow{m}\rangle D_1^{i_1}\cdots D_n^{i_n}b
=D_t(a\langle \overrightarrow{m}\rangle D_1^{i_1}\cdots
D_{t}^{i_{t}-1}\cdots D_n^{i_n}b)+ m_t a\langle
\overrightarrow{m}-\overrightarrow{e_t}\rangle D_1^{i_1}\cdots
D_{t}^{i_{t}-1}\cdots D_n^{i_n}b.
$$
Now, put $\overrightarrow{N}(a,D_1^{i_1}\cdots
D_n^{i_n}b)=\overrightarrow{N}(a,D_1^{i_1}\cdots
D_{t}^{i_{t}-1}\cdots D_n^{i_n}b)+\overrightarrow{e_t}$. Then
$$
a\langle \overrightarrow{m}\rangle D_1^{i_1}\cdots D_n^{i_n}b=0,\
for \ \overrightarrow{m}\not\prec
\overrightarrow{N}(a,D_1^{i_1}\cdots D_n^{i_n}b).
$$

Let $|v|>1$, and $v=b\langle \overrightarrow{m}'\rangle[v_1]$. Take
$$
\overrightarrow{N}(a,v)=\overrightarrow{N}(a,b)+\overrightarrow{N}(b,v_1)-\overrightarrow{m}'.
$$
Then for $\overrightarrow{m}\not\prec  \overrightarrow{N}(a,v)$,
clearly $\overrightarrow{m}\not\prec  \overrightarrow{N}(a,b)$.
Therefore, we have (by the multiplication algorithm)
$$
a\langle \overrightarrow{m}\rangle b\langle
\overrightarrow{m}'\rangle[v_1] =\sum_{\overrightarrow{s}^1\in
   Z_+^n\backslash 0}\alpha_{\overrightarrow{s}^1}a\langle
   \overrightarrow{m}-\overrightarrow{s}^1\rangle (b\langle
   \overrightarrow{m}'+\overrightarrow{s}^1\rangle [v_1]).
$$
where $\overrightarrow{m}'+\overrightarrow{s}^1 \prec
\overrightarrow{N}(b,v_1)$. This means that
$$
\overrightarrow{s}^1\prec\overrightarrow{N}(b,v_1)-\overrightarrow{m}'
\ \ \ and  \ \ \ \overrightarrow{m}-\overrightarrow{s}^1 \not\prec
\overrightarrow{N}(a,b).
$$
Hence by the multiplication algorithm again we have
$$
\sum_{\overrightarrow{s}^1,\overrightarrow{s}^2\in Z_+^n\backslash
0}\alpha_{\overrightarrow{s}^1}\alpha_{\overrightarrow{s}^2}a
      \langle \overrightarrow{m}-\overrightarrow{s}^1-\overrightarrow{s}^2\rangle
      (b\langle \overrightarrow{m}'+\overrightarrow{s}^1+\overrightarrow{s}^2\rangle [v_1])
$$
where
$$
\overrightarrow{m}'+\overrightarrow{s}^1+\overrightarrow{s}^2 \prec
\overrightarrow{N}(b,v_1) \ and \
\overrightarrow{m}-\overrightarrow{s}^1-\overrightarrow{s}^2
\not\prec  \overrightarrow{N}(a,b).
$$
Repeating the above argument $k=\sum_{i=1}^nt_i-m'_i$ times where
$\overrightarrow{N}(b,v_1)=(t_1,,\cdots,t_n)$, we arrive at
$$
a\langle \overrightarrow{m}\rangle b\langle
\overrightarrow{m}'\rangle[v_1] =\sum_{\overrightarrow{s}^1,\cdots,
\overrightarrow{s}^k\in Z_+^n\backslash
0}\alpha_{\overrightarrow{s}^1}\cdots\alpha_{\overrightarrow{s}^k}a
      \langle \overrightarrow{m}-\Sigma_{i=1}^k \overrightarrow{s}^i\rangle
      (b\langle \overrightarrow{m}'+\Sigma_{i=1}^k \overrightarrow{s}^i\rangle [v_1]).
$$
In this case, $\overrightarrow{m}'+\Sigma_{i=1}^k
\overrightarrow{s}^i \not\prec  N(b,v_1)$ and so the last expression
is zero. By induction, we have $[u]\langle
\overrightarrow{m}\rangle[v]=0$ if $\overrightarrow{m} \not\prec
\overrightarrow{N}(u,v)$.

Next, we continue our argument on the assumption that $|u|>1$. Thus,
$[u]=a\langle \overrightarrow{m}'\rangle[u_1]$ for $a\in B$. Then
$$
(a\langle \overrightarrow{m}'\rangle[u_1])\langle
\overrightarrow{m}\rangle[v]=a\langle
\overrightarrow{m}'\rangle([u_1]\langle
\overrightarrow{m}\rangle[v])+\sum_{\overrightarrow{s}\in
   Z_+^n\backslash 0}\alpha_{\overrightarrow{s}} a\langle
   \overrightarrow{m}'-\overrightarrow{s}\rangle ([u_1]\langle
   \overrightarrow{m}+\overrightarrow{s}\rangle [v]).
$$
Set $\overrightarrow{N}(u,v)=\overrightarrow{N}(u_1,v)$. Then
$[u]\langle \overrightarrow{m}\rangle[v]=0$ if
$\overrightarrow{m}\not\prec  \overrightarrow{N}(u,v)$ and the
locality holds.

(ii) The identity
\begin{equation}\label{e2.4}
D_i([u]\langle \overrightarrow{m}\rangle[v])=D_i[u]\langle
\overrightarrow{m}\rangle[v]+[u]\langle \overrightarrow{m}\rangle
D_i[v],\ \ 1\leqslant i \leqslant n.
\end{equation}

Assume first $|u|=|v|=1$. Then $[u]=a\in B$ and $[v]=D_1^{i_1}\cdots
D_n^{i_n}b,\ i_1,\cdots,i_n\in Z_+,\ b\in B$. The case when
$\overrightarrow{m}\prec \overrightarrow{N}(a,b)$ holds is true from
the definition. So suppose that $\overrightarrow{m} \not\prec
\overrightarrow{N}$. If $i_1=\cdots=i_n=0$, then the left hand side
of (\ref{e2.4}) is zero while the right hand side of (\ref{e2.4}) is
$$
-m_ia\langle \overrightarrow{m}-\overrightarrow{e_i}\rangle
b+a\langle \overrightarrow{m}\rangle D_ib
$$
and
$$
a\langle \overrightarrow{m}\rangle D_ib=D_i(a\langle
\overrightarrow{m}\rangle b)+m_ia\langle
\overrightarrow{m}-\overrightarrow{e_i}\rangle b \ \ (by \ Remark)
$$
$$
=m_ia\langle \overrightarrow{m}-\overrightarrow{e_i}\rangle b \ \ \
\ \ \ \ \ \ \ \ \ \ \ \ \ \ \ \ \ \ \ \ \
$$
Hence (\ref{e2.4}) is true. So let
$i_{l_1},\cdots,i_{l_k}>0,i_{l_{k+1}}=\cdots=i_{l_n}=0,1<k\leqslant
n$. We prove that if $\overrightarrow{m}\not\prec
\overrightarrow{N}(a,b)$ then
\begin{equation}\label{e2.5}
D_i(a\langle \overrightarrow{m}\rangle D_{l_1}^{i_{l_1}}\cdots
D_{l_k}^{i_{l_k}}b)=-m_i a\langle
\overrightarrow{m}-\overrightarrow{e_i}\rangle
D_{l_1}^{i_{l_1}}\cdots D_{l_k}^{i_{l_k}}b +a\langle
\overrightarrow{m}\rangle D_iD_{l_1}^{i_{l_1}}\cdots
D_{l_k}^{i_{l_k}}b
\end{equation}
Suppose $s_1,\cdots,s_k\in Z_+$. Let
$\alpha=(-1)^{s_1+\cdots+s_k},\beta=\binom{i_{l_1}}{s_1}\cdots\binom{i_{l_k}}{s_k},
\beta'=\binom{i_{l_1}}{s_1}\cdots\binom{i_{l_t}+1}{s_t}\cdots\binom{i_{l_k}}{s_k}$,
$\gamma=[m_1]^{s_1}\cdots[m_k]^{s_k}$ where
$[m_j]^{s_j}=m_j(m_j-1)\cdots(m_j-s_j+1),1\leqslant j\leqslant k$.
Expanding the left hand side of (\ref{e2.5}), we have
\begin{eqnarray*}
&&D_i(a\langle \overrightarrow{m}\rangle D_{l_1}^{i_{l_1}}\cdots
D_{l_k}^{i_{l_k}}b\\
&=&D_i(D_{l_1}^{i_{l_1}}\cdots D_{l_k}^{i_{l_k}}(a\langle
\overrightarrow{m}\rangle b)-\sum_{(s_1,\cdots,s_k)\in
   Z_+^k\backslash 0}\beta (D_{l_1}^{s_1}\cdots
D_{l_k}^{s_k}a\langle \overrightarrow{m}\rangle
D_{l_1}^{i_{l_1}-s_1}\cdots D_{l_k}^{i_{l_k}-s_k}b))\\
&=&-D_i(\sum_{(s_1,\cdots,s_k)\in
   Z_+^k\backslash 0}\alpha\gamma\beta (a\langle \overrightarrow{m}-\Sigma_{p=1}^k s_p
   \overrightarrow{e_{l_p}}\rangle
D_{l_1}^{i_{l_1}-s_1}\cdots D_{l_k}^{i_{l_k}-s_k}b))\triangleq A.
\end{eqnarray*}
If $i=l_t \in
\{l_1,\cdots,l_k \}$, by induction on $(i_1,\cdots,i_n)$,
\begin{eqnarray}\label{e2.6}
A&=&\sum_{(s_1,\cdots,s_k)\in
   Z_+^k\backslash 0}\alpha\gamma\beta m_i(a\langle \overrightarrow{m}-
   \Sigma_{p=1}^k s_p\overrightarrow{e_{l_p}}\rangle
D_{l_1}^{i_{l_1}-s_1}\cdots D_{l_k}^{i_{l_k}-s_k}b)) \\
&&\label{e2.7} -\sum_{(s_1,\cdots,s_k)\in
   Z_+^k\backslash 0}\alpha\gamma\beta (a\langle \overrightarrow{m}-
   \Sigma_{p=1}^k s_p\overrightarrow{e_{l_p}}\rangle
D_{l_1}^{i_{l_1}-s_1}\cdots D_i^{i_{l_t}+1-s_i}\cdots
D_{l_k}^{i_{l_k}-s_k}b)
\end{eqnarray}
If $i \in \{l_{k+1},\cdots,l_n \}$, by induction on
$(i_1,\cdots,i_n)$,
\begin{eqnarray}\label{e2.8}
A&=&\sum_{(s_1,\cdots,s_k)\in
   Z_+^k\backslash 0}\alpha\gamma\beta m_i(a\langle \overrightarrow{m}-
   \Sigma_{p=1}^k s_p\overrightarrow{e_{l_p}}-e_i\rangle
D_{l_1}^{i_{l_1}-s_1}\cdots D_{l_k}^{i_{l_k}-s_k}b))\\
&&\label{e2.9} -\sum_{(s_1,\cdots,s_k)\in
   Z_+^k\backslash 0}\alpha\gamma\beta (a\langle \overrightarrow{m}-
   \Sigma_{p=1}^k s_p\overrightarrow{e_{l_p}}\rangle
D_iD_{l_1}^{i_{l_1}-s_1}\cdots
D_{l_k}^{i_{l_k}-s_k}b)
\end{eqnarray}
On the other hand,
the right hand side of (\ref{e2.5}) can be rewritten as
$$
-m_i a\langle \overrightarrow{m}-\overrightarrow{e_i}\rangle
D_{l_1}^{s_1}\cdots D_{l_k}^{s_k}b+a\langle
\overrightarrow{m}\rangle D_iD_{l_1}^{s_1}\cdots
D_{l_k}^{s_k}b\triangleq B.
$$
If $i=l_t\in \{l_1,\cdots,l_k \}$, then
\begin{eqnarray*}
B&=&-m_i a\langle \overrightarrow{m}-\overrightarrow{e_i}\rangle
D_{l_1}^{s_1}\cdots D_{l_k}^{s_k}b+D_{l_1}^{i_{l_1}}\cdots
D^{i_{l_t}+1}_i\cdots D_{l_k}^{i_{l_k}}(a\langle
\overrightarrow{m}\rangle b)\\
&& -\sum_{(s_1,\cdots,s_k)\in
   Z_+^k\backslash 0}\beta' (D_{l_1}^{s_1}\cdots
D_{l_k}^{s_k}a\langle \overrightarrow{m}\rangle
D_{l_1}^{i_{l_1-s_1}}\cdots D^{i_{l_t}+1-s_t}_i\cdots
D_{l_k}^{i_{l_k}-s_k}b))\\
&=&-m_i a\langle \overrightarrow{m}-\overrightarrow{e_i}\rangle
D_{l_1}^{s_1}\cdots D_{l_k}^{s_k}b\\
&& -\sum_{(s_1,\cdots,s_k)\in
   Z_+^k\backslash 0}\alpha\gamma\beta' (a\langle \overrightarrow{m}-
   \Sigma_{p=1}^k s_p\overrightarrow{e_{l_p}}\rangle
D_{l_1}^{i_{l_1-s_1}}\cdots D^{i_{l_t}+1-s_t}_i\cdots
D_{l_k}^{i_{l_k}-s_k}b)
\end{eqnarray*}
From (\ref{e2.6}) and (\ref{e2.7}), we see that

1) \ $m_ij a\langle \overrightarrow{m}-\overrightarrow{e_i}\rangle
D_{l_1}^{i_{l_1}}\cdots D_{l_k}^{i_{l_k}}b =-m_i a\langle
\overrightarrow{m}-\overrightarrow{e_i}\rangle
D_{l_1}^{i_{l_1}}\cdots D_{l_k}^{i_{l_k}}b+m_i(j+1) a\langle
\overrightarrow{m}-\overrightarrow{e_i}\rangle
D_{l_1}^{i_{l_1}}\cdots D_{l_k}^{i_{l_k}}b; $

2) \ for  $s_i\geqslant 1$,
\begin{eqnarray*}
&&(-1)^{s_i}[m_i]^{s_i+1}(\binom{i_{l_t}}{s_i}+\binom{i_{l_t}}{s_i+1})a\langle
\overrightarrow{m}-\Sigma_{p=1}^k
s_p\overrightarrow{e_{l_p}}-\overrightarrow{e_i}\rangle
D_{l_1}^{i_{l_1}-s_1}\cdots D_{l_k}^{i_{l_k}-s_k}b\\
&=&(-1)^{s_i}[m_i]^{s_i+1}\binom{i_{l_t}+1}{s_i+1}a\langle
\overrightarrow{m}-\Sigma_{p=1}^k
s_p\overrightarrow{e_{l_p}}-\overrightarrow{e_i}\rangle
D_{l_1}^{i_{l_1}-s_1}\cdots D_{l_k}^{i_{l_k}-s_k}b.
\end{eqnarray*}

It is easy to check that $A=B$.

If $i \in \{l_{k+1},\cdots,l_n \}$, then
\begin{eqnarray*}
B&=&-m_i a\langle \overrightarrow{m}-\overrightarrow{e_i}\rangle
D_{l_1}^{i_{l_1}}\cdots D_{l_k}^{i_{l_k}}b+D_{l_1}^{i_{l_1}}\cdots
D_{l_k}^{i_{l_k}}D_i(a\langle \overrightarrow{m}\rangle b)\\
&&-\sum_{(s_1,\cdots,s_k,s_i)\in
   Z_+^{k+1}\backslash 0}\beta \binom{1}{s_i}(D_{l_1}^{s_1}\cdots
D_{l_k}^{s_k}D_ia\langle \overrightarrow{m}\rangle
D_{l_1}^{i_{l_1}-s_1}\cdots D_{l_k}^{i_{l_k}-s_k}D_i^{1-s_i}b)\\
&=&-m_i a\langle \overrightarrow{m}-\overrightarrow{e_i}\rangle
D_{l_1}^{i_{l_1}}\cdots D_{l_k}^{i_{l_k}}b-D_ia\langle
\overrightarrow{m}\rangle D_{l_1}^{i_{l_1}}\cdots
D_{l_k}^{i_{l_k}}b\\
&&-\sum_{(s_1,\cdots,s_k)\in
   Z_+^k\backslash 0}\beta (D_{l_1}^{s_1}\cdots
D_{l_k}^{s_k}D_ia\langle \overrightarrow{m}\rangle
D_{l_1}^{i_{l_1}-s_1}\cdots D_{l_k}^{i_{l_k}-s_k}b)\\
&&-\sum_{(s_1,\cdots,s_k)\in
   Z_+^k\backslash 0}\beta (D_{l_1}^{s_1}\cdots
D_{l_k}^{s_k}a\langle \overrightarrow{m}\rangle
D_{l_1}^{i_{l_1}-s_1}\cdots D_{l_k}^{i_{l_k}-s_k}D_ib)\\
&=&\sum_{(s_1,\cdots,s_k)\in
   Z_+^k\backslash 0}\alpha\gamma\beta m_i(a\langle \overrightarrow{m}
   -\Sigma_{p=1}^k s_p\overrightarrow{e_{l_p}}-\overrightarrow{e_i}\rangle
D_{l_1}^{i_{l_1}-s_1}\cdots
D_{l_k}^{i_{l_k}-s_k}b))\\
&&-\sum_{(s_1,\cdots,s_k)\in
   Z_+^k\backslash 0}\alpha\gamma\beta (a\langle \overrightarrow{m}
   -\Sigma_{p=1}^k s_p\overrightarrow{e_{l_p}}\rangle
D_iD_{l_1}^{i_{l_1}-s_1}\cdots D_{l_k}^{i_{l_k}-s_k}b).
\end{eqnarray*}
From (\ref{e2.8}) and (\ref{e2.9}) it follows that $A=B$.

This completes the case of $|u|=|v|=1$.
\par We proceed by induction on $|u|+|v|=l$ and assume
that $l\geqslant 3$. Suppose that $|u|=1$. Then $|v|\geqslant 2$,
and so $[v]=b\langle \overrightarrow{m}'\rangle [v_1]$ for some
$v_1$. Thus we have to establish  (\ref{e2.4}), namely,
\begin{eqnarray}\label{e2.12}
D_i(a\langle \overrightarrow{m}\rangle(b\langle
\overrightarrow{m}'\rangle [v_1]))=D_ia\langle
\overrightarrow{m}\rangle(b\langle \overrightarrow{m}'\rangle [v_1])
+a\langle \overrightarrow{m}\rangle D_i(b\langle
\overrightarrow{m}'\rangle [v_1]),\ 1\leqslant i \leqslant n.
\end{eqnarray}
If $\overrightarrow{m}\prec \overrightarrow{N}(a,b)$, (\ref{e2.12})
follows immediately from the definition of derivation of normal
words. Assume that $\overrightarrow{m}\not\prec
\overrightarrow{N}(a,b)$. By the definition of multiplication, we
have
$$
a\langle \overrightarrow{m}\rangle (b\langle
\overrightarrow{m}'\rangle[v_1])=(a\langle \overrightarrow{m}\rangle
b)\langle \overrightarrow{m}'\rangle[v_1]
-\sum_{\overrightarrow{s}\in
   Z_+^n\backslash 0}(-1)^{\overrightarrow{s}}\binom{\overrightarrow{m}}{\overrightarrow{s}} a\langle
   \overrightarrow{m}-\overrightarrow{s}\rangle (b\langle
   \overrightarrow{m}'+\overrightarrow{s}\rangle [v_1]).
$$
Using induction on $\overrightarrow{m}$, the left hand side of
(\ref{e2.12}) becomes
\begin{eqnarray*}
-\sum_{\overrightarrow{s}\in
   Z_+^n\backslash 0}(-1)^{\overrightarrow{s}}\binom{\overrightarrow{m}}{\overrightarrow{s}} (D_ia\langle
   \overrightarrow{m}-\overrightarrow{s}\rangle (b\langle
   \overrightarrow{m}'+\overrightarrow{s}\rangle [v_1])-a\langle
   \overrightarrow{m}-\overrightarrow{s}\rangle D_i(b\langle
   \overrightarrow{m}'+\overrightarrow{s}\rangle [v_1])).
\end{eqnarray*}
Now, using the multiplicity algorithm and induction on $|v|$, the
above becomes
\begin{eqnarray*}
A_0&\triangleq&\sum_{\overrightarrow{s}\in
   Z_+^n\backslash 0}(-1)^{\overrightarrow{s}}\binom{\overrightarrow{m}}{\overrightarrow{s}} (m_i-s_i)a\langle
   \overrightarrow{m}-\overrightarrow{s}-\overrightarrow{e_i} \rangle (b\langle
   \overrightarrow{m}'+\overrightarrow{s}\rangle [v_1])\\
&&\ \ \ +\sum_{\overrightarrow{s}\in
   Z_+^n\backslash 0}(-1)^{\overrightarrow{s}}\binom{\overrightarrow{m}}{\overrightarrow{s}} (m'_i+s_i)a\langle
   \overrightarrow{m}-\overrightarrow{s}\rangle (b\langle
   \overrightarrow{m}'+\overrightarrow{s}-\overrightarrow{e_i}\rangle
   [v_1])\\
&&\ \ \ -\sum_{\overrightarrow{s}\in
   Z_+^n\backslash 0}(-1)^{\overrightarrow{s}}\binom{\overrightarrow{m}}{\overrightarrow{s}} a\langle
   \overrightarrow{m}-\overrightarrow{s}\rangle (b\langle
   \overrightarrow{m}'+\overrightarrow{s}\rangle D_i[v_1]).
\end{eqnarray*}
The right hand side of (\ref{e2.12}) is expended to
\begin{eqnarray*}
B_0&\triangleq&D_ia\langle \overrightarrow{m}\rangle(b\langle
\overrightarrow{m}'\rangle [v_1])+a\langle \overrightarrow{m}\rangle
D_i(b\langle \overrightarrow{m}'\rangle [v_1])\\
&=&-m_ia\langle m-e_i\rangle(b\langle m'\rangle [v_1]) -m'_ia\langle
m\rangle (b\langle m'-e_i\rangle [v_1]) +a\langle m\rangle (b\langle
m'\rangle D_i[v_1])\\
&=&-m_ia\langle
\overrightarrow{m}-\overrightarrow{e_i}\rangle(b\langle
\overrightarrow{m}'\rangle [v_1])\\
&& +\sum_{\overrightarrow{s}\in
   Z_+^n\backslash 0}(-1)^{\overrightarrow{s}}\binom{\overrightarrow{m}}{\overrightarrow{s}} m'_ia\langle
   \overrightarrow{m}-\overrightarrow{s}\rangle (b\langle
   \overrightarrow{m}'+\overrightarrow{s}-\overrightarrow{e_i}\rangle
   [v_1])\\
&&-\sum_{\overrightarrow{s}\in
   Z_+^n\backslash 0}(-1)^{\overrightarrow{s}}\binom{\overrightarrow{m}}{\overrightarrow{s}} a\langle
   \overrightarrow{m}-\overrightarrow{s}\rangle (b\langle
   \overrightarrow{m}'+\overrightarrow{s}\rangle D_i[v_1]).
\end{eqnarray*}
It is now straightforward to check that $A_0$ and $B_0$ are
coincide.

Hence (\ref{e2.4}) holds when $|u|=1$.

\par Next, assume $|u|>1$. Thus $[u]=a\langle \overrightarrow{m}'\rangle[u_1], \overrightarrow{m}'\in
Z_+^n$, $\overrightarrow{m}'\prec \overrightarrow{N}(a,b)$, $b$ is
the first letter of $[u_1]$. We shall make use of the property
$$
D_i[u]\langle \overrightarrow{m}\rangle[v]=-m_i[u]\langle
\overrightarrow{m}-\overrightarrow{e_i}\rangle[v]
$$
which will be shown in (\ref{e2.16}).
\par First of all, let's figure out all summands in (\ref{e2.4}). Suppose $s_1,\cdots,s_n\in Z_+$. Let
$\beta=\binom{m'_1}{s_1}\cdots\binom{m'_i}{s_i+1}\cdots\binom{m'_n}{s_n}$.
Using the multiplicity algorithm with induction on $|u|$ and
$|u|+|v|$, we have
$$
D_i((a\langle \overrightarrow{m}'\rangle[u_1])\langle
\overrightarrow{m}\rangle[v]) =D_i(\sum_{\overrightarrow{s}\in
   Z_+^n}(-1)^{\overrightarrow{s}}\binom{\overrightarrow{m}'}{\overrightarrow{s}} a\langle
   \overrightarrow{m}'-\overrightarrow{s}\rangle ([u_1]\langle
   \overrightarrow{m}+\overrightarrow{s}\rangle [v]))
=A_1+A_2+A_3
$$
where
\begin{eqnarray*}
A_1&=&-\sum_{\overrightarrow{s}\in
   Z_+^n}(-1)^{\overrightarrow{s}}\binom{\overrightarrow{m}'}{\overrightarrow{s}} (m'_i-s_i)a\langle
   \overrightarrow{m}'-\overrightarrow{s}-\overrightarrow{e_i}\rangle ([u_1]\langle
   \overrightarrow{m}+\overrightarrow{s}\rangle [v]),\\
A_2&=&-\sum_{\overrightarrow{s}\in
   Z_+^n}(-1)^{\overrightarrow{s}}\binom{\overrightarrow{m}'}{\overrightarrow{s}} (m_i+s_i)a\langle
   \overrightarrow{m}'-\overrightarrow{s}\rangle ([u_1]\langle
   \overrightarrow{m}+\overrightarrow{s}-\overrightarrow{e_i}\rangle
   [v]),\\
A_3&=&\sum_{\overrightarrow{s}\in
   Z_+^n}(-1)^{\overrightarrow{s}}\binom{\overrightarrow{m}'}{\overrightarrow{s}} a\langle
   \overrightarrow{m}'-\overrightarrow{s}\rangle ([u_1]\langle
   \overrightarrow{m}+\overrightarrow{s}\rangle D_i[v]).
\end{eqnarray*}
On the other hand,
\begin{eqnarray*}
D_i[u]\langle \overrightarrow{m}\rangle[v]&=&-m_i[u]\langle
\overrightarrow{m}-\overrightarrow{e_i}\rangle[v]\\
&=&-m_i(a\langle \overrightarrow{m}'\rangle[u_1])\langle
\overrightarrow{m}-\overrightarrow{e_i}\rangle[v]\\
&=&-m_i\sum_{\overrightarrow{s}\in
   Z_+^n}(-1)^{\overrightarrow{s}}\binom{\overrightarrow{m}'}{\overrightarrow{s}} a\langle
   \overrightarrow{m}'-\overrightarrow{s}\rangle ([u_1]\langle
   \overrightarrow{m}+\overrightarrow{s}-\overrightarrow{e_i}\rangle [v])\\
   &\triangleq&B_1
\end{eqnarray*}
and
\begin{eqnarray*}
[u]\langle \overrightarrow{m}\rangle D_i[v]&=&(a\langle
\overrightarrow{m}'\rangle[u_1])\langle \overrightarrow{m}\rangle
D_i[v]\\
&=&\sum_{\overrightarrow{s}\in
   Z_+^n}(-1)^{\overrightarrow{s}}\binom{\overrightarrow{m}'}{\overrightarrow{s}} a\langle
   \overrightarrow{m}'-\overrightarrow{s}\rangle ([u_1]\langle
   \overrightarrow{m}+\overrightarrow{s}\rangle D_i[v])\\
   &\triangleq&B_2
\end{eqnarray*}
Now, $A_3=B_2$. Subtracting $B_1$ from $A_2$, we get
$$ C=A_2-B_1=-\sum_{\overrightarrow{s}\in
   Z_+^n}(-1)^{\overrightarrow{s}}\binom{\overrightarrow{m}'}{\overrightarrow{s}} s_ia\langle
   \overrightarrow{m}'-\overrightarrow{s}\rangle ([u_1]\langle
   \overrightarrow{m}+\overrightarrow{s}-\overrightarrow{e_i}\rangle [v])
$$
$$
 \ \ \ \ \ \ \ \ =\sum_{\overrightarrow{s}\in
   Z_+^n}(-1)^{\overrightarrow{s}}\beta (s_i+1)a\langle
   \overrightarrow{m}'-\overrightarrow{s}-\overrightarrow{e_i}\rangle ([u_1]\langle
   \overrightarrow{m}+\overrightarrow{s}\rangle [v]).
$$
Finally, we note that $C+A_1=0$ since
$(s_i+1)\binom{m'_i}{s_i+1}=(m'_i-s_i)\binom{m'_i}{s_i}$.

(iii) The identity
\begin{eqnarray}\label{e2.16}
D_i[u]\langle \overrightarrow{m}\rangle[v]=-m_i[u]\langle
\overrightarrow{m}-\overrightarrow{e_i}\rangle[v]
\end{eqnarray}
\par To show this identity, we use induction on $|u|$. If $|u|=1$,
then $u=D_1^{i_1}\cdots D_n^{i_n}a$, and the result follows from the
definition. So let $|u|>1$, and write $[u]=a\langle
\overrightarrow{m}'\rangle[u_1]$. Suppose $s_1,\cdots,s_n\in Z_+$.
Let $\beta=\binom{m'_1}{s_1}\cdots
\binom{m'_i-1}{s_i}\cdots\binom{m'_n}{s_n}$,
$\beta'=\binom{m'_1}{s_1}\cdots
\binom{m'_i}{s_i+1}\cdots\binom{m'_n}{s_n}$.

Using the multiplicity algorithm with induction on $|u|$, the left
hand side of (\ref{e2.16}) is equal to
\begin{eqnarray*}
&&(D_ia\langle \overrightarrow{m}'\rangle[u_1])\langle
\overrightarrow{m}\rangle[v]+(a\langle \overrightarrow{m}'\rangle
D_i[u_1])\langle \overrightarrow{m}\rangle[v]\\
&=&-m'_i(a\langle
\overrightarrow{m}'-\overrightarrow{e_i}\rangle[u_1])\langle
\overrightarrow{m}\rangle[v]+(a\langle \overrightarrow{m}'\rangle
D_i[u_1])\langle \overrightarrow{m}\rangle[v]\\
&=&\sum_{\overrightarrow{s}\in
   Z_+^n}(-1)^{\overrightarrow{s}}(-\beta m'_ia\langle
   \overrightarrow{m}'-\overrightarrow{s}-\overrightarrow{e_i}\rangle ([u_1]\langle
   \overrightarrow{m}+\overrightarrow{s}\rangle [v])\\
&&+\binom{\overrightarrow{m}'}{\overrightarrow{s}} a\langle
   \overrightarrow{m}'-\overrightarrow{s}\rangle (D_i[u_1]\langle
   \overrightarrow{m}+\overrightarrow{s}\rangle [v]))\\
&\triangleq&A_1+A_2
\end{eqnarray*}
where
\begin{eqnarray*}
A_1&=&-m'_i\sum_{\overrightarrow{s}\in
   Z_+^n}(-1)^{\overrightarrow{s}}\beta a\langle
   \overrightarrow{m}'-\overrightarrow{s}-\overrightarrow{e_i}\rangle ([u_1]\langle
   \overrightarrow{m}+\overrightarrow{s}\rangle [v]),\\
A_2&=&-\sum_{\overrightarrow{s}\in
   Z_+^n}(-1)^{\overrightarrow{s}}\binom{\overrightarrow{m}'}{\overrightarrow{s}} (m_i+s_i)a\langle
   \overrightarrow{m}'-\overrightarrow{s}\rangle ([u_1]\langle
   \overrightarrow{m}+\overrightarrow{s}-\overrightarrow{e_i}\rangle [v]).
\end{eqnarray*}
The right hand side of (\ref{e2.16}) is equal to
$$
B=-m_i(a\langle \overrightarrow{m}'\rangle[u_1])\langle
\overrightarrow{m}-\overrightarrow{e_i}\rangle[v]
 =-m_i\sum_{\overrightarrow{s}\in
   Z_+^n}(-1)^{\overrightarrow{s}}\binom{\overrightarrow{m}'}{\overrightarrow{s}} a\langle
   \overrightarrow{m}'-\overrightarrow{s}\rangle ([u_1]\langle
   \overrightarrow{m}+\overrightarrow{s}-\overrightarrow{e_i}\rangle [v]).
$$
Now, $A_1+A_2=B$. Indeed, subtracting $B$ from $A_2$, we get
\begin{eqnarray*}
C=A_2-B&=&-\sum_{\overrightarrow{s}\in
   Z_+^n}(-1)^{\overrightarrow{s}}\binom{\overrightarrow{m}'}{\overrightarrow{s}} s_ia\langle
   \overrightarrow{m}'-\overrightarrow{s}\rangle ([u_1]\langle
   \overrightarrow{m}+\overrightarrow{s}-\overrightarrow{e_i}\rangle
   [v])\\
&=&\sum_{\overrightarrow{s}\in
   Z_+^n}(-1)^{\overrightarrow{s}}\beta'(s_i+1) a\langle
   \overrightarrow{m}'-\overrightarrow{s}-\overrightarrow{e_i}\rangle ([u_1]\langle
   \overrightarrow{m}+\overrightarrow{s}\rangle [v])
\end{eqnarray*}
and $A_1+C=0$ as
$(s_i+1)\binom{m'_i}{s_i+1}=m'_i\binom{m'_i-1}{s_i}$.

(iv) (Associativity condition): We have to show that
\begin{eqnarray}\label{e2.17}
   ([u]\langle \overrightarrow{m}\rangle [v])\langle \overrightarrow{m}'\rangle
   [w]=\sum_{\overrightarrow{s}\in
   Z_+^n}(-1)^{\overrightarrow{s}}\binom{\overrightarrow{m}}{\overrightarrow{s}}[u]\langle
   \overrightarrow{m}-\overrightarrow{s}\rangle ([v]\langle
   \overrightarrow{m}'+\overrightarrow{s}\rangle [w])
\end{eqnarray}

Suppose $s_1,\cdots,s_n\in Z_+$. Let $\beta=\binom{m_1}{s_1}\cdots
\binom{m_t-1}{s_t}\cdots\binom{m_n}{s_n}$,
$\beta'=\binom{m_1}{s_1}\cdots
\binom{m_t}{s_t+1}\cdots\binom{m_n}{s_n}$. We consider first some
special cases.

Case 1. $u=a$ and $|v|=1$. Then $v=D_1^{i_1}\cdots D_n^{i_n}b$. When
$i_1=\cdots=i_n=0$ and $\overrightarrow{m}\prec
\overrightarrow{N}(a,b)$, (\ref{e2.17}) follows from the definition.
When $i_1=\cdots=i_n=0$ and $\overrightarrow{m}\not\prec
\overrightarrow{N}(a,b)$, the left hand side of (\ref{e2.17}) is
equal to 0, while the right hand side of (\ref{e2.17}) contains the
summand
$$
a\langle \overrightarrow{m}\rangle (b\langle
\overrightarrow{m}'\rangle[w]) =(a\langle \overrightarrow{m}\rangle
b)\langle \overrightarrow{m}'\rangle[w] -\sum_{\overrightarrow{s}\in
   Z_+^n\backslash 0}(-1)^{\overrightarrow{s}}\binom{\overrightarrow{m}}{\overrightarrow{s}}[u]\langle
   \overrightarrow{m}-\overrightarrow{s}\rangle ([v]\langle
   \overrightarrow{m}'+\overrightarrow{s}\rangle [w]).
$$
Hence the right hand side of (\ref{e2.17}) is 0 by induction on
$\overrightarrow{m}$.

Suppose now that $i_t>0$ for some $t, \ 1\leqslant t\leqslant n$.
Then by multiplicity algorithm and induction on $(i_1,\cdots,i_n)$,
the left hand side of (\ref{e2.17}) is equal to
\begin{eqnarray*}
&&(a\langle \overrightarrow{m}\rangle D_1^{i_1}\cdots
D_n^{i_n}b)\langle \overrightarrow{m}'\rangle[w]\\
&=&D_t(a\langle \overrightarrow{m}\rangle D_1^{i_1}\cdots
D_{t}^{i_t-1}\cdots D_n^{i_n}b)\langle \overrightarrow{m}'\rangle[w]
+m_t (a\langle \overrightarrow{m}-\overrightarrow{e_t}\rangle
D_1^{i_1}\cdots D_{t}^{i_t-1}\cdots D_n^{i_n}b)\langle
\overrightarrow{m}'\rangle[w]\\
&=&-m'_t(a\langle \overrightarrow{m}\rangle D_1^{i_1}\cdots
D_{t}^{i_t-1}\cdots D_n^{i_n}b)\langle
\overrightarrow{m}'-\overrightarrow{e_t}\rangle[w]\\
&&+m_t (a\langle \overrightarrow{m}-\overrightarrow{e_t}\rangle
D_1^{i_1}\cdots D_{t}^{i_t-1}\cdots D_n^{i_n}b)\langle
\overrightarrow{m}'\rangle[w]\\
&\triangleq&A_1+A_2
\end{eqnarray*}
where
\begin{eqnarray*}
A_1&=&-m'_t\sum_{\overrightarrow{s}\in
   Z_+^n}(-1)^{\overrightarrow{s}}\binom{\overrightarrow{m}}{\overrightarrow{s}} a\langle
   \overrightarrow{m}-\overrightarrow{s}\rangle (D_1^{i_1}\cdots
D_{t}^{i_t-1}\cdots D_n^{i_n}b\langle
   \overrightarrow{m}'+\overrightarrow{s}-\overrightarrow{e_t}\rangle
   [w]),\\
A_2&=&\sum_{\overrightarrow{s}\in
   Z_+^n}(-1)^{\overrightarrow{s}}\beta m_t a\langle
   \overrightarrow{m}-\overrightarrow{s}-\overrightarrow{e_t}\rangle (D_1^{i_1}\cdots
D_{t}^{i_t-1}\cdots D_n^{i_n}b\langle
   \overrightarrow{m}'+\overrightarrow{s}\rangle [w]).
\end{eqnarray*}
The right hand side of (\ref{e2.17}) is equal to
\begin{eqnarray*}
&&\sum_{\overrightarrow{s}\in
   Z_+^n}(-1)^{\overrightarrow{s}}\binom{\overrightarrow{m}}{\overrightarrow{s}} a\langle
   \overrightarrow{m}-\overrightarrow{s}\rangle (D_1^{i_1}\cdots D_n^{i_n}b\langle
   \overrightarrow{m}'+\overrightarrow{s}\rangle [w])\\
&=&-\sum_{\overrightarrow{s}\in
   Z_+^n}(-1)^{\overrightarrow{s}}\binom{\overrightarrow{m}}{\overrightarrow{s}} (m'_t+s_t)a\langle
   \overrightarrow{m}-\overrightarrow{s}\rangle (D_1^{i_1}\cdots
D_{t}^{i_t-1}\cdots D_n^{i_n}b\langle
   \overrightarrow{m}+\overrightarrow{s}-\overrightarrow{e_t}\rangle [w])\\
&\triangleq&B
\end{eqnarray*}
Hence,
\begin{eqnarray*}
B-A_1&=&-\sum_{\overrightarrow{s}\in
   Z_+^n}(-1)^{\overrightarrow{s}}\binom{\overrightarrow{m}}{\overrightarrow{s}} s_t a\langle
   \overrightarrow{m}-\overrightarrow{s}\rangle (D_1^{i_1}\cdots
D_{t}^{i_t-1}\cdots D_n^{i_n}b\langle
   \overrightarrow{m}'+\overrightarrow{s}-\overrightarrow{e_t}\rangle
   [w])\\
&=&\sum_{\overrightarrow{s}\in
   Z_+^n}(-1)^{\overrightarrow{s}}\beta' (s_k+1) a\langle
   \overrightarrow{m}-\overrightarrow{s}-\overrightarrow{e_t}\rangle (D_1^{i_1}\cdots
D_{t}^{i_t-1}\cdots D_n^{i_n}b\langle
   \overrightarrow{m}'+\overrightarrow{s}\rangle [w])
\end{eqnarray*}
which is the same as $A_2$ since
$m_t\binom{m_t-1}{s_t}=(s_t+1)\binom{m_t}{s_t+1}$.

Case 2. $u=D_1^{i_1}\cdots D_n^{i_n}a, \ i_t>0$ for some $t, \
1\leqslant t\leqslant n$ and $|v|=1$. Using induction on
$(i_1,\cdots,i_n)$ (with Case 1 providing the induction case), we
see that the left hand side of (\ref{e2.17}) is equal to
\begin{eqnarray*}
&&(D_tu_1\langle \overrightarrow{m}\rangle [v])\langle
\overrightarrow{m}'\rangle[w]\\
&=&-m_t(u_1\langle \overrightarrow{m}-\overrightarrow{e_t}\rangle
[v])\langle \overrightarrow{m}'\rangle[w]\\
&=&-m_t\sum_{\overrightarrow{s}\in
   Z_+^n}(-1)^{\overrightarrow{s}}\beta  u_1\langle
   \overrightarrow{m}-\overrightarrow{s}-\overrightarrow{e_t}\rangle ([v]\langle
   \overrightarrow{m}'+\overrightarrow{s}\rangle [w])
\end{eqnarray*}
where $u_1=D_1^{i_1}\cdots D_{t}^{i_t-1}\cdots D_n^{i_n}a$.

The right hand side of (\ref{e2.17}), on the other hand, is equal to
$$
-\sum_{\overrightarrow{s}\in
   Z_+^n}(-1)^{\overrightarrow{s}}\binom{\overrightarrow{m}}{\overrightarrow{s}} (m_t-s_t)  u_1\langle
   \overrightarrow{m}-\overrightarrow{s}-\overrightarrow{e_t}\rangle ([v]\langle
   \overrightarrow{m}'+\overrightarrow{s}\rangle [w]).
$$
Because $m_t\binom{m_t-1}{s_t}=(m_t-s_t)\binom{m_t}{s_t}$, we see
that (\ref{e2.17}) holds, and this settles the case for $|u|=1$.

It follows from Case 1 and Case 2 that, in particular, (\ref{e2.17})
holds for $|u|=|v|=|w|=1$. We will use induction on $|u|+|v|+|w|$,
and so we assume $|u|+|v|+|w|> 3$. Also, to simplify notation, when
we are dealing with a right normed bracketing expression, the
brackets are simply dropped.

Case 3. $u=a$ and $|v|=1$. In this case, we just repeat the argument
of Case 1, and we are done.

Case 4. $u=a$ and $|v|>1$. Then $v=b\langle
\overrightarrow{k}\rangle[v_1]$ for some $v_1$. If
$\overrightarrow{m}\prec \overrightarrow{N}(a,b)$, then $[u]\langle
\overrightarrow{m}\rangle[v]=a\langle \overrightarrow{m}\rangle[v]$
is a normal word and (\ref{e2.17}) follows from the definition of
the product of normal words (see Remark). Hence, assume that
$\overrightarrow{m}\not\prec  \overrightarrow{N}(a,b)$. By
definition, we have
$$
a\langle \overrightarrow{m}\rangle(b\langle
\overrightarrow{k}\rangle[v_1]) =-\sum_{\overrightarrow{s}\in
   Z_+^n\backslash 0}(-1)^{\overrightarrow{s}}\binom{\overrightarrow{m}}{\overrightarrow{s}} a\langle
   \overrightarrow{m}-\overrightarrow{s}\rangle (b\langle
   \overrightarrow{k}+\overrightarrow{s}\rangle [v_1]).
$$
Applying induction on $\overrightarrow{m}$ (with $|a|=1$ and
$|b\langle \overrightarrow{k}+\overrightarrow{s}\rangle[v_1]|=|v|$)
and on $|u|+|v|+|w|$, we have
\begin{eqnarray*}
&&(a\langle \overrightarrow{m}\rangle[v])\langle
\overrightarrow{m}'\rangle[w]\\
&=&-\sum_{\overrightarrow{s}\in
   Z_+^n\backslash 0}(-1)^{\overrightarrow{s}}\binom{\overrightarrow{m}}{\overrightarrow{s}} (a\langle
   \overrightarrow{m}-\overrightarrow{s}\rangle (b\langle
   \overrightarrow{k}+\overrightarrow{s}\rangle [v_1]))\langle
   \overrightarrow{m}'\rangle[w]\\
&=&-\sum_{\overrightarrow{s},\overrightarrow{s}^2\in
   Z_+^n\backslash 0}(-1)^{\overrightarrow{s}}(-1)^{\overrightarrow{s}^2}
   \binom{\overrightarrow{m}}{\overrightarrow{s}}
   \binom{\overrightarrow{m}-\overrightarrow{s}}{\overrightarrow{s}^2} a\langle
   \overrightarrow{m}-\overrightarrow{s}-\overrightarrow{s}^2\rangle ((b\langle
   \overrightarrow{k}+\overrightarrow{s}\rangle [v_1])\langle
   \overrightarrow{m}'+\overrightarrow{s}^2\rangle[w])\\
&=&-\sum_{\overrightarrow{s},\overrightarrow{s}^2,\overrightarrow{s}^3\in
   Z_+^n\backslash 0}
   \xi_{_{\overrightarrow{s},\overrightarrow{s}^2,\overrightarrow{s}^3}}
     a\langle
   \overrightarrow{m}-\overrightarrow{s}-\overrightarrow{s}^2\rangle b\langle
   \overrightarrow{k}+\overrightarrow{s}-\overrightarrow{s}^3\rangle [v_1]\langle
   \overrightarrow{m}'+\overrightarrow{s}^2+\overrightarrow{s}^3\rangle[w]\\
&\triangleq& A_0
\end{eqnarray*}
where $
\xi_{_{\overrightarrow{s},\overrightarrow{s}^2,\overrightarrow{s}^3}}
=(-1)^{\overrightarrow{s}+\overrightarrow{s}^2+\overrightarrow{s}^3}
\binom{\overrightarrow{m}}{\overrightarrow{s}}
   \binom{\overrightarrow{m}-\overrightarrow{s}}{\overrightarrow{s}^2}
   \binom{\overrightarrow{k}+\overrightarrow{s}}{\overrightarrow{s}^3}$.

Now, let us take the right hand side of (\ref{e2.17}) for $u=a$,
$v=b\langle \overrightarrow{k}\rangle[v_1]$,
$\overrightarrow{m}\not\prec  \overrightarrow{N}(a,b)$. Then
\begin{eqnarray*}
&&a\langle \overrightarrow{m}\rangle(b\langle
\overrightarrow{k}\rangle[v_1])\langle
\overrightarrow{m}'\rangle[w]+\sum_{\overrightarrow{s}\in
   Z_+^n\backslash 0}(-1)^{\overrightarrow{s}}\binom{\overrightarrow{m}}{\overrightarrow{s}} a\langle
   \overrightarrow{m}-\overrightarrow{s}\rangle (b\langle
   \overrightarrow{k}\rangle [v_1])\langle
   \overrightarrow{m}'+\overrightarrow{s}\rangle[w]\\
&=&\sum_{\overrightarrow{s}^2\in Z^n_+}
  (-1){\overrightarrow{s}^2}\binom{\overrightarrow{k}}{\overrightarrow{s}^2} a\langle
   \overrightarrow{m}\rangle b\langle
   \overrightarrow{k}-\overrightarrow{s}^2\rangle [v_1]\langle
   \overrightarrow{m}'+\overrightarrow{s}^2\rangle[w]\\
&&+\sum_{\overrightarrow{s},\overrightarrow{s}^2\in
   Z_+^n\backslash 0}
  (-1)^{\overrightarrow{s}}(-1){\overrightarrow{s}^2}
  \binom{\overrightarrow{m}}{\overrightarrow{s}}
  \binom{\overrightarrow{k}}{\overrightarrow{s}^2} a\langle
   \overrightarrow{m}-\overrightarrow{s}\rangle b\langle
   \overrightarrow{k}-\overrightarrow{s}^2\rangle [v_1]\langle
   \overrightarrow{m}'+\overrightarrow{s}+\overrightarrow{s}^2\rangle[w]\\
&\triangleq&A_1+A_2
\end{eqnarray*}
where
\begin{eqnarray*}
&&A_1=-\sum_{\overrightarrow{s},\overrightarrow{s}^2\in
Z_+^n\backslash 0}
      (-1)^{\overrightarrow{s}}(-1){\overrightarrow{s}^2}
      \binom{\overrightarrow{m}}{\overrightarrow{s}}\binom{\overrightarrow{k}}{\overrightarrow{s}^2} a\langle
   \overrightarrow{m}-\overrightarrow{s}\rangle b\langle
       \overrightarrow{k}+\overrightarrow{s}-\overrightarrow{s}^2\rangle [v_1]\langle
       \overrightarrow{m}'+\overrightarrow{s}^2\rangle[w],\\
&&A_2=\sum_{\overrightarrow{s},\overrightarrow{s}^2\in
Z_+^n\backslash
0}(-1)^{\overrightarrow{s}}(-1){\overrightarrow{s}^2}
\binom{\overrightarrow{m}}{\overrightarrow{s}}\binom{\overrightarrow{k}}{\overrightarrow{s}^2}
a\langle
   \overrightarrow{m}-\overrightarrow{s}\rangle b\langle
   \overrightarrow{k}-\overrightarrow{s}^2\rangle [v_1]\langle
   \overrightarrow{m}'+\overrightarrow{s}+\overrightarrow{s}^2\rangle[w].
\end{eqnarray*}
It is not difficult to see that $A_1+A_2$ is equal to $A_0$: First,
make a transformation
$$
(\overrightarrow{i},\overrightarrow{j},\overrightarrow{l})=
(\overrightarrow{m}-\overrightarrow{s}-\overrightarrow{s}^2,
\overrightarrow{k}+\overrightarrow{s}-\overrightarrow{s}^2,
\overrightarrow{m}'+\overrightarrow{s}^2+\overrightarrow{s}^3)
$$
where $\overrightarrow{i},\overrightarrow{j},\overrightarrow{l}\in
Z_+^n$ and so
$\overrightarrow{i}+\overrightarrow{j}+\overrightarrow{l}=
\overrightarrow{m}+\overrightarrow{m}'+\overrightarrow{k}$. Then
$A_0$ becomes a sum of the expressions
$$
-\sum_{\overrightarrow{s}\in
    Z_+^n\backslash 0}(-1)^{\overrightarrow{s}+\overrightarrow{l}+\overrightarrow{m}'}
    \binom{\overrightarrow{m}}{\overrightarrow{s}}
    \binom{\overrightarrow{m}-\overrightarrow{s}}{\overrightarrow{m}-\overrightarrow{s}-\overrightarrow{i}}
    \binom{\overrightarrow{k}+\overrightarrow{s}}{\overrightarrow{k}+\overrightarrow{s}-\overrightarrow{j}}
    a\langle \overrightarrow{i}\rangle
    b\langle \overrightarrow{j}\rangle[v_1]\langle
    \overrightarrow{l}\rangle[w].
$$
Next, do a similar transformation
$$
(\overrightarrow{i},\overrightarrow{j},\overrightarrow{l})=
(\overrightarrow{m}-\overrightarrow{s},
\overrightarrow{k}+\overrightarrow{s}-\overrightarrow{s}^2,
\overrightarrow{m}'+\overrightarrow{s}^2)
$$
with $\overrightarrow{i},\overrightarrow{j},\overrightarrow{l}\in
Z_+^n$. Then $A_1$ becomes a sum of the expressions
$$
     -(-1)^{\overrightarrow{k}+\overrightarrow{j}}
     \binom{\overrightarrow{m}}{\overrightarrow{m}-\overrightarrow{i}}
     \binom{\overrightarrow{k}}{\overrightarrow{l}-\overrightarrow{m}'}
     a\langle \overrightarrow{i}\rangle
     b\langle \overrightarrow{j}\rangle[v_1]\langle
     \overrightarrow{l}\rangle[w].
$$
Do another transformation
$$
(\overrightarrow{i},\overrightarrow{j},\overrightarrow{l})
=(\overrightarrow{m}-\overrightarrow{s},\overrightarrow{k}-\overrightarrow{s}^2,
\overrightarrow{m}'+\overrightarrow{s}+\overrightarrow{s}^2)
$$
with $\overrightarrow{i},\overrightarrow{j},\overrightarrow{l}\in
Z_+^n$. Let $\gamma_2=(-1)^{\sum_{i=1}^nl_i+m'_i}$,
$\delta_2=\binom{m_1}{m_1-i_1}\cdots\binom{m_n}{m_n-i_n}$,
$\delta'_2=\binom{k_1}{k_1-j_1}\cdots\binom{k_n}{k_n-j_n}$. Then
$A_2$ becomes a sum of the expressions
$$
     (-1)^{\overrightarrow{l}+\overrightarrow{m}'}
     \binom{\overrightarrow{m}}{\overrightarrow{m}-\overrightarrow{i}}
     \binom{\overrightarrow{k}}{\overrightarrow{k}-\overrightarrow{j}}
      a\langle \overrightarrow{i}\rangle
     b\langle \overrightarrow{j}\rangle[v_1]\langle
     \overrightarrow{l}\rangle[w].
$$
Thus the associativity law for $u=a$ will follow if we can establish
the following identity

\begin{eqnarray}
&&-\sum_{\overrightarrow{s}\in
    Z_+^n\backslash
    0}(-1)^{\overrightarrow{s}+\overrightarrow{l}+\overrightarrow{m}'}
    \binom{\overrightarrow{m}}{\overrightarrow{s}}
    \binom{\overrightarrow{m}-\overrightarrow{s}}{\overrightarrow{m}-\overrightarrow{s}-\overrightarrow{i}}
    \binom{\overrightarrow{k}+\overrightarrow{s}}{\overrightarrow{k}+\overrightarrow{s}-\overrightarrow{j}}
    \nonumber\\ \label{*}
&=&-(-1)^{\overrightarrow{k}+\overrightarrow{j}}
    \binom{\overrightarrow{m}}{\overrightarrow{m}-\overrightarrow{i}}
    \binom{\overrightarrow{k}}{\overrightarrow{l}-\overrightarrow{m}'}
    +(-1)^{\overrightarrow{l}+\overrightarrow{m}'}\binom{m}{m-i}
    \binom{\overrightarrow{k}}{\overrightarrow{k}-\overrightarrow{j}}
\end{eqnarray}

Clearly, (\ref{*}) is equivalent to
$$
\sum_{\overrightarrow{s}\in
    Z_+^n}(-1)^{\overrightarrow{s}+\overrightarrow{l}+\overrightarrow{m}'}
    \binom{\overrightarrow{m}}{\overrightarrow{s}}
    \binom{\overrightarrow{m}-\overrightarrow{s}}{\overrightarrow{m}-\overrightarrow{s}-\overrightarrow{i}}
    \binom{\overrightarrow{k}+\overrightarrow{s}}{\overrightarrow{k}+\overrightarrow{s}-\overrightarrow{j}}
    =(-1)^{\overrightarrow{k}+\overrightarrow{j}}
    \binom{\overrightarrow{m}}{\overrightarrow{m}-\overrightarrow{i}}
    \binom{\overrightarrow{k}}{\overrightarrow{l}-\overrightarrow{m}'}
$$
which is equivalent to
\begin{eqnarray}\label{e2.22}
\sum_{\overrightarrow{s}\in
    Z_+^n}(-1)^{\overrightarrow{s}+\overrightarrow{l}+\overrightarrow{m}'}
    \binom{\overrightarrow{m}-\overrightarrow{i}}{\overrightarrow{s}}
    \binom{\overrightarrow{k}+\overrightarrow{s}}{\overrightarrow{j}}
    =(-1)^{\overrightarrow{k}+\overrightarrow{j}}
    \binom{\overrightarrow{k}}{\overrightarrow{l}-\overrightarrow{m}'}
\end{eqnarray}
due to the equality
$$
\binom{\overrightarrow{m}}{\overrightarrow{s}}
\binom{\overrightarrow{m}-\overrightarrow{s}}{\overrightarrow{m}-\overrightarrow{s}-\overrightarrow{i}}=
\binom{\overrightarrow{m}}{\overrightarrow{i}}\binom{\overrightarrow{m}-\overrightarrow{i}}{\overrightarrow{s}}.
$$
Using a substitution
$$
(\overrightarrow{a},\overrightarrow{b})=(\overrightarrow{m}-\overrightarrow{i},
\overrightarrow{l}-\overrightarrow{m}')
$$
with $\overrightarrow{a},\overrightarrow{b}\in Z_+^n$, and
$\overrightarrow{a}+\overrightarrow{k}=\overrightarrow{b}+\overrightarrow{j}$,
(\ref{e2.22}) becomes
$$
(-1)^{\overrightarrow{a}}\sum_{\overrightarrow{s}\in Z_+^n}
\binom{\overrightarrow{a}}{\overrightarrow{s}}
\binom{\overrightarrow{k}+\overrightarrow{s}}{\overrightarrow{j}}=
\binom{\overrightarrow{k}}{\overrightarrow{a}+\overrightarrow{k}-\overrightarrow{j}}
$$
which is equivalent to
$$
(-1)^{a_t}\sum_{s_t\in Z_+} \binom{a_t}{s_t}\binom{k_t+s_t}{j_t}=
\binom{k_t}{a_t+k_t-j_t}, \ 1\leqslant t \leqslant n .
$$
The last identity is indeed valid (see, for example, D. Knuth. The
art of computer programming).

Case 5. $u=D_1^{i_1}\cdots D_n^{i_n}a, \ i_t>0$ for some $t,\
1\leqslant t\leqslant n$. In this case, we can repeat the argument
from Case 2, and we are done.

Case 6. $|u|>1$. In this case, we write $[u]=a\langle
\overrightarrow{k}\rangle[u_1]$ where $\overrightarrow{k}\in Z_+^n$,
$\overrightarrow{k}\prec \overrightarrow{N}(a,b), \ a,b\in B$ and
$b$ is the first letter of $[u_1]$.

By definition, the left hand side of (\ref{e2.17}) is equal to
\begin{eqnarray*}
&&((a\langle \overrightarrow{k}\rangle[u_1])\langle
\overrightarrow{m}\rangle[v])\langle \overrightarrow{m}'\rangle[w]\\
&=&\sum_{\overrightarrow{s}\in
   Z_+^n}(-1)^{\overrightarrow{s}} \binom{\overrightarrow{k}}{\overrightarrow{s}} (a\langle
   \overrightarrow{k}-\overrightarrow{s}\rangle [u_1]\langle
   \overrightarrow{m}+\overrightarrow{s}\rangle [v])\langle
   \overrightarrow{m}'\rangle[w]\\
&=&\sum_{\overrightarrow{s},\overrightarrow{s}^2\in
   Z_+^n}(-1)^{\overrightarrow{s}}(-1)^{\overrightarrow{s}^2}
   \binom{\overrightarrow{k}}{\overrightarrow{s}}
   \binom{\overrightarrow{k}-\overrightarrow{s}}{\overrightarrow{s}^2} a\langle
   \overrightarrow{k}-\overrightarrow{s}-\overrightarrow{s}^2\rangle ([u_1]\langle
   \overrightarrow{m}+\overrightarrow{s}\rangle [v_1])\langle
   \overrightarrow{m}'+\overrightarrow{s}^2\rangle[w]\\
&=&\sum_{\overrightarrow{s},\overrightarrow{s}^2,\overrightarrow{s}^3\in
   Z_+^n}
     \xi_{_{\overrightarrow{s},\overrightarrow{s}^2,\overrightarrow{s}^3}}a\langle
   \overrightarrow{k}-\overrightarrow{s}-\overrightarrow{s}^2\rangle[u_1]\langle
   \overrightarrow{m}+\overrightarrow{s}-\overrightarrow{s}^3\rangle [v]\langle
   \overrightarrow{m}'+\overrightarrow{s}^2+\overrightarrow{s}^3\rangle[w]\\
&\triangleq&A_1
\end{eqnarray*}
where
$\xi_{_{\overrightarrow{s},\overrightarrow{s}^2,\overrightarrow{s}^3}}=
(-1)^{\overrightarrow{s}+\overrightarrow{s}^2+\overrightarrow{s}^3}
   \binom{\overrightarrow{k}}{\overrightarrow{s}}
   \binom{\overrightarrow{k}-\overrightarrow{s}}{\overrightarrow{s}^2}
   \binom{\overrightarrow{m}+\overrightarrow{s}}{\overrightarrow{s}^3}$.

Doing a transformation on the indices
$$
(\overrightarrow{i},\overrightarrow{j},\overrightarrow{l})=
(\overrightarrow{k}-\overrightarrow{s}-\overrightarrow{s}^2,
\overrightarrow{m}+\overrightarrow{s}-\overrightarrow{s}^3,
\overrightarrow{m}'+\overrightarrow{s}^2+\overrightarrow{s}^3)
$$
with $\overrightarrow{i},\overrightarrow{j},\overrightarrow{l}\in
Z_+^n$, $A_1$ then becomes a sum of
\begin{eqnarray}\label{e2.24}
\sum_{\overrightarrow{s}\in
    Z_+^n}(-1)^{\overrightarrow{s}+\overrightarrow{l}+\overrightarrow{m}'}
    \binom{\overrightarrow{k}}{\overrightarrow{s}}
    \binom{\overrightarrow{k}-\overrightarrow{s}}{\overrightarrow{k}-\overrightarrow{s}-\overrightarrow{i}}
    \binom{\overrightarrow{m}+\overrightarrow{s}}{\overrightarrow{m}+\overrightarrow{s}-\overrightarrow{j}}
    a\langle \overrightarrow{i}\rangle [u_1]\langle\overrightarrow{j}\rangle[v]\langle
     \overrightarrow{l}\rangle[w]
\end{eqnarray}

The right hand side of (\ref{e2.17}) is equal to
\begin{eqnarray*}
&&\sum_{\overrightarrow{s}\in
   Z_+^n}(-1)^{\overrightarrow{s}} \binom{\overrightarrow{m}}{\overrightarrow{s}}
   (a\langle \overrightarrow{k}\rangle[u_1])\langle
   \overrightarrow{m}-\overrightarrow{s}\rangle ([v]\langle
   \overrightarrow{m}'+\overrightarrow{s}\rangle [w])\\
&=&\sum_{\overrightarrow{s},\overrightarrow{s}^2\in
    Z_+^n}
    (-1)^{\overrightarrow{s}+\overrightarrow{s}^2}  \binom{\overrightarrow{m}}{\overrightarrow{s}}
    \binom{\overrightarrow{k}}{\overrightarrow{s}^2}a\langle
    \overrightarrow{k}-\overrightarrow{s}^2\rangle[u_1]\langle
    \overrightarrow{m}-\overrightarrow{s}+\overrightarrow{s}^2\rangle [v]\langle
    \overrightarrow{m}'+\overrightarrow{s}\rangle[w]\\
&\triangleq&A_2.
\end{eqnarray*}

Another transformation on indices
$$
(\overrightarrow{i},\overrightarrow{j},\overrightarrow{l})=
(\overrightarrow{k}-\overrightarrow{s}^2,
\overrightarrow{m}-\overrightarrow{s}+\overrightarrow{s}^2,
\overrightarrow{m}'+\overrightarrow{s})
$$
with $\overrightarrow{i},\overrightarrow{j},\overrightarrow{l}\in
Z_+^n$, brings $A_2$ to a sum of
\begin{eqnarray}\label{e2.26}
(-1)^{\overrightarrow{m}+\overrightarrow{j}}
\binom{\overrightarrow{m}}{\overrightarrow{l}-\overrightarrow{m}'}
\binom{\overrightarrow{k}}{\overrightarrow{k}-\overrightarrow{i}}
    a\langle \overrightarrow{i}\rangle [u_1]\langle \overrightarrow{j}\rangle[v]\langle
     \overrightarrow{l}\rangle[w]
\end{eqnarray}

Finally, the fact that
$$
\sum_{\overrightarrow{s}\in Z_+^n}
    (-1)^{\overrightarrow{s}+\overrightarrow{l}+\overrightarrow{m}'}
    \binom{\overrightarrow{k}}{\overrightarrow{s}}
    \binom{\overrightarrow{k}-\overrightarrow{s}}{\overrightarrow{k}-\overrightarrow{s}-\overrightarrow{i}}
    \binom{\overrightarrow{m}+\overrightarrow{s}}{\overrightarrow{m}+\overrightarrow{s}-\overrightarrow{j}}
 =(-1)^{\overrightarrow{m}+\overrightarrow{j}}
 \binom{\overrightarrow{m}}{\overrightarrow{l}-\overrightarrow{m}'}
 \binom{\overrightarrow{k}}{\overrightarrow{k}-\overrightarrow{i}}
$$
indicates that (\ref{e2.24}) and (\ref{e2.26}) are identical (this
is (\ref{*}) after only an exchange of the roles of
$\overrightarrow{m}$ and $\overrightarrow{k}$).

This completes our proof that
$C(B,\overrightarrow{N},D_1,\cdots,D_n)$ is  a free associative
$n$-conformal algebra. \ \ $\square$

\begin{corollary}
The normal words consist of a linear basis of the free associative
$n$-conformal algebra $C(B,\overrightarrow{N},D_1,\cdots,D_n)$.
\end{corollary}

\section{Composition-Diamond lemma}

Here and after, we only consider the situation when the locality
function $\overrightarrow{N}(a,b),a,b\in B$, is uniformly bounded by
some $\overrightarrow{N}$. Then without loss of generality, we may
assume that $\overrightarrow{N}(a,b)=\overrightarrow{N}$ for all
$a,b\in B$.
\par For $\overrightarrow{i}=(i_1,\dots,i_n)\in Z_+^n$, put $D^{\overrightarrow{i}}=D_1^{i_1}\cdots
D_n^{i_n}$.

We assume that any polynomial of $f\in
C(B,\overrightarrow{N},D_1,\cdots,D_n)$ is presented as a linear
combination of normal words
\begin{eqnarray}\label{e3.1}
[u]=a_1\langle
\overrightarrow{m}^{(1)}\rangle (a_2\langle
    \overrightarrow{m}^{(2)}\rangle(\cdots (a_k\langle
    \overrightarrow{m}^{(k)}\rangle D^{\overrightarrow{i}}a_{k+1})\cdots))
\end{eqnarray}
where $a_l\in B$, $\overrightarrow{m}^{(j)}\in Z_+^n$,
$\overrightarrow{m}^{(j)}\prec \overrightarrow{N}$, $1\leqslant l
\leqslant k+1, 1\leqslant j \leqslant k$,
ind(u)=$\overrightarrow{i}\in Z_+^n$, $|u|=k+1,k\geqslant 0$. We
shall refer to $[u]$ as D-free if ind(u)=$(0,\cdots,0)$, and we
shall say that $f$ is D-free if every normal word in $f$ is D-free.
Also, put
\begin{eqnarray}\label{e3.2}
u=a_1\langle \overrightarrow{m}^{(1)}\rangle a_2\langle
    \overrightarrow{m}^{(2)}\rangle\cdots a_k\langle
    \overrightarrow{m}^{(k)}\rangle D^{\overrightarrow{i}}a_{k+1}
\end{eqnarray}
and call it an associative normal word.

We order the forms (\ref{e3.1}) and (\ref{e3.2}) according to the
lexicographical ordering of their weights:
$$
wt(u)=(|u|,a_1,\overrightarrow{m}^{(1)},a_2,\overrightarrow{m}^{(2)},
\cdots,a_k,\overrightarrow{m}^{(k)},a_{k+1},\overrightarrow{i}).
$$
For $f\in C(B,\overrightarrow{N},D_1,\cdots,D_n)$, the leading
associative word of $f$ is denoted by $\overline{f}$. So
\begin{eqnarray*}
f&=&\alpha_{\overline{f}}[\overline{f}]+\Sigma_i\alpha_i[u_i], \\
\overline{f}&=&a_1\langle \overrightarrow{m}^{(1)}\rangle a_2\langle
\overrightarrow{m}^{(2)}\rangle\cdots a_k\langle
               \overrightarrow{m}^{(k)}\rangle
               D^{\overrightarrow{i}}a_{k+1},\\
u_i&<&\overline{f}.
\end{eqnarray*}

We denote deg$(f)=|\overline{f}|$, and we will call $f$ $monic$ if
$\alpha_{\overline{f}}=1$.
\par All associative normal words (\ref{e3.2}) form a $\Omega$-semigroup
$T$ with 0, where $\Omega=\{\langle \overrightarrow{m}\rangle
|\overrightarrow{m}\in Z_+^n\}$. The multiplication in $T$ is
defined by the following formulas:
$$
 \ \  \ \ \ \ \ \ \ \ D^{\overrightarrow{i}}a \langle \overrightarrow{m}\rangle x
 =a \langle \overrightarrow{m}-\Sigma_{p=1}^ni_p\overrightarrow{e_p}\rangle x ,
$$
$$
\ \ \ \    a \langle \overrightarrow{m}\rangle x=0, \ \
\overrightarrow{0}\not\prec  \overrightarrow{m} ,
$$
$$
 \ \ \ \ \  a \langle \overrightarrow{m}\rangle x=0, \ \ \overrightarrow{m}\not\prec  \overrightarrow{N}
$$
where $a,b\in B$ and $x\in D^\omega(B)$.
\par We note that the associative law
$$
(u\langle \overrightarrow{m}\rangle v)\langle
\overrightarrow{m}'\rangle w= u\langle \overrightarrow{m}\rangle
(v\langle \overrightarrow{m}'\rangle w)
$$
holds trivially in $T$.

For any associative normal word
$$
u=a_1\langle \overrightarrow{m}^{(1)}\rangle a_2\langle
    \overrightarrow{m}^{(2)}\rangle\cdots a_k\langle
    \overrightarrow{m}^{(k)}\rangle D^{\overrightarrow{i}}a_{k+1},
$$
we denote by $(u)$ a nonassociative word on $u$, i.e., $(u)$ is some
bracketing of $u$.

\begin{lemma}\label{l3.1}
Let $[p]$ and $[q]$ be normal words. If $[p]$ is D-free, then for
any $\overrightarrow{m}\in Z_+^n$, $\overrightarrow{m}\prec
\overrightarrow{N}$,
$$
[p]\langle \overrightarrow{m}\rangle[q]=[p\langle
\overrightarrow{m}\rangle q]+\Sigma\gamma_i[q_i]
$$
where $[q_i]$ are normal, $q_i<p\langle \overrightarrow{m}\rangle
q$, and $|q_i|=|p|+|q|$. If $q$ is D-free, so are $q_i$'s.
\end{lemma}
\noindent{\bf Proof.} If $|p|=1$, then $p=a\in B$, and $[p]\langle
\overrightarrow{m}\rangle[q]$ is normal. So let $|p|>1$ and write
$p=a_1\langle \overrightarrow{m}'\rangle [p_1]$, where
$\overrightarrow{m}'\prec \overrightarrow{N}$. Then
\begin{eqnarray}\label{e3.4}
(a_1\langle \overrightarrow{m}'\rangle [p_1])\langle
\overrightarrow{m}\rangle[q]=a_1\langle \overrightarrow{m}'\rangle
([p_1]\langle \overrightarrow{m}\rangle[q])
+\sum_{\overrightarrow{s}\in
   Z_+^n\backslash 0}\alpha_{\overrightarrow{s}} a_1\langle
   \overrightarrow{m}-\overrightarrow{s}\rangle ([p_1]\langle
   \overrightarrow{m}'+\overrightarrow{s}\rangle [q])\ \ \ \
\end{eqnarray}
Using induction on $|p|$, we have
$$
a_1\langle \overrightarrow{m}'\rangle ([p_1]\langle
\overrightarrow{m}\rangle[q])=a_1\langle \overrightarrow{m}'\rangle
[p_1\langle \overrightarrow{m}\rangle q] +\Sigma\beta_ia_1\langle
\overrightarrow{m}'\rangle[u_i]
$$
where
$$
|u_i|=|p_1|+|q|, \ \  u_i<p_1\langle \overrightarrow{m}\rangle q , \
\  a_1\langle \overrightarrow{m}'\rangle u_i<a_1\langle
\overrightarrow{m}'\rangle p_1\langle \overrightarrow{m}\rangle
q=p\langle \overrightarrow{m}\rangle q.
$$
By Lemma \ref{l2.2}, the second summand in the right hand side of
(\ref{e3.4}) is a linear combination of normal words of the
following form:
$$
a_1\langle \overrightarrow{m}-\overrightarrow{s}\rangle[w],\
|w|=|p_1|+|q| \ and \ \overrightarrow{s}\in Z_+^n\backslash 0.
$$
Hence
$$
a_1\langle \overrightarrow{m}-\overrightarrow{s}\rangle w
<a_1\langle \overrightarrow{m}'\rangle p_1\langle
\overrightarrow{m}\rangle q=p\langle \overrightarrow{m}\rangle q, \
\ |a_1\langle \overrightarrow{m}-\overrightarrow{s}\rangle
w|=|p\langle \overrightarrow{m}\rangle q|.
$$
This completes the proof of the lemma. \ \  $\square$
\begin{lemma}\label{l3.2}
Let $u=a_1\langle \overrightarrow{m}^{(1)}\rangle \cdots
    a_k\langle \overrightarrow{m}^{(k)}\rangle D^{\overrightarrow{i}}a_{k+1} $
    be an associative normal word. Then for any
    nonassociative word $(u)$, $(u)=[u]+\Sigma_l\alpha_l[u_l]$, where
    $[u_l]$ are normal words, $u_l<u$ and $|u_l|=|u|$. In
    particular, $\overline{(u)}=u$. If $u$ is D-free, then the
    $u_l$'s are D-free as well.
\end{lemma}
\noindent{\bf Proof.} Induction on $|u|$. The case of $|u|=1$ is
clear. Let $|u|>1$, and so $(u)=(v)\langle
\overrightarrow{m}\rangle(w)$ for some associative normal words $w$
and $v$ with $v$ being D-free, and $\overrightarrow{m}\in Z_+^n$,
$\overrightarrow{m}\prec \overrightarrow{N}$. By induction, we have
$$
(v)=[v]+\Sigma\alpha_l[v_l],(w)=[w]+\Sigma\beta_j[w_j]
$$
where $v_l<v$, $w_j<w$, $|v_l|=|v|$, $|w_j|=|w|$, $ v_l$'s are D-free.\\
The lemma would follow from Lemma \ref{l3.1}. \ \  $\square$

\begin{definition}
Let $S\subset C(B,\overrightarrow{N},D_1,\cdots,D_n)$ be a set of
monic polynomials. We define $S$-words
$(u)_{D^{\overrightarrow{i}}s}$ by induction. (Actually,
$(u)_{D^{\overrightarrow{i}}s}$ is a word with the occurrence of a
polynomial $D^{\overrightarrow{i}}s, s\in S$.)

1)\
$(D^{\overrightarrow{i}}s)_{D^{\overrightarrow{i}}s}=D^{\overrightarrow{i}}s
$ where $s\in S$, is an $S$-word of $S$-length 1.

2)\ If $u_{D^{\overrightarrow{i}}s}$ is an $S$-word of $S$-length
$k$, and $(v)$ is any word in $B$ of length $l$, then $
(u)_{D^{\overrightarrow{i}}s}\langle \overrightarrow{m}\rangle(v)$
and $(v)\langle
\overrightarrow{m}\rangle(u)_{D^{\overrightarrow{i}}s}$ are
$S$-words of $S$-length $k+l$.
\end{definition}

The $S$-length of an $S$-word $(u)_{D^{\overrightarrow{i}}s}$ will
be denoted by $|u|_{D^{\overrightarrow{i}}s}$.

\ \

\noindent{\bf Remark}. Any element of Id$(D^\omega(S))$ is a linear
combination of $S$-words. Yet for $S$-word
$(u)_{D^{\overrightarrow{i}}s}$, the length
$|u|_{D^{\overrightarrow{i}}s}$ may increase under some of the
transformation that follow. For this reason, we need the notion of
formal degree (fdeg) of expression in
$C(B,\overrightarrow{N},D_1,\cdots,D_n)$.
\begin{definition}
Let $\Sigma_{i=1}^k\alpha_i(u_i),\alpha_i\neq 0$, be a linear
combination of some words $u_i$ in $B$. Put
$$
fdeg(\sum_{i=1}^k\alpha_i(u_i))=\max_{1\leqslant i\leqslant k}|u_i|.
$$
\end{definition}
Of course, fdeg$(\Sigma_{i=1}^k\alpha_i(u_i))$ depends on the
presentation of $\Sigma_{i=1}^k\alpha_i(u_i)$. It is also clear that
fdeg$(\Sigma \alpha_i(u_i))\geqslant$deg$(\Sigma \alpha_i(u_i))$,
and if the $(u_i)$'s are normal, then fdeg$(\Sigma
\alpha_i(u_i))=$deg$(\Sigma \alpha_i(u_i))$.
\begin{definition}
Let $S$ be a set of monic elements in
$C(B,\overrightarrow{N},D_1,\cdots,D_n)$. An associative normal
$S$-word is an expression of the following forms:
\begin{eqnarray}\label{e3.5}
 u_s=a\langle
\overrightarrow{m}\rangle s\langle \overrightarrow{m}'\rangle b
\end{eqnarray}
where $s\in S$ is D-free, $a,b\in T$ with
being D-free, $\overrightarrow{m},\overrightarrow{m}'\in Z_+^n$,
$\overrightarrow{m},\overrightarrow{m}'\prec \overrightarrow{N}$; or
\begin{eqnarray}\label{e3.6}
u_{D^{\overrightarrow{i}}s}=a\langle \overrightarrow{m}\rangle
D^{\overrightarrow{i}}s
\end{eqnarray}
where $s\in S, \overrightarrow{i}\in Z_+^n, a\in T$ is D-free,
$\overrightarrow{m}\in Z_+^n$, $\overrightarrow{m}\prec
\overrightarrow{N}$.
\end{definition}
Here we note that $a$ and $b$ may be empty in (\ref{e3.5}) and
(\ref{e3.6}).

The words in (\ref{e3.5}) will be referred to associative normal
$S$-words of the first kind, while those words in (\ref{e3.6})
associative normal $S$-words of the second kind. The elements
$$
[u]_s=[a\langle \overrightarrow{m}\rangle s\langle
\overrightarrow{m}'\rangle b] \ \ \mbox{ and } \ \
[u]_{D^{\overrightarrow{i}}s}=[a\langle \overrightarrow{m}\rangle
D^{\overrightarrow{i}}s]
$$
will be referred to normal $S$-words of the first kind and the
second kind, respectively. A common notation for words (\ref{e3.5})
and (\ref{e3.6}) is
$$
 u_{D^{\overrightarrow{i}}s}=a\langle
\overrightarrow{m}\rangle D^{\overrightarrow{i}}s\langle
\overrightarrow{m}'\rangle b
$$
where $s$ is D-free and $\overrightarrow{i}=(0,\cdots,0)$ if $b\neq
1$, while $s$ and $\overrightarrow{i}\in Z_+^n$ are arbitrary if
$b=1$.
\par From now on, $S$ will denote a set of monic polynomials.
\begin{lemma}\label{l3.6}
Let $[u]_{D^{\overrightarrow{i}}s}=[a\langle
\overrightarrow{m}\rangle D^{\overrightarrow{i}}s\langle
\overrightarrow{m}'\rangle b]$ be a normal $S$-word, $a,b\in T$.
Then $\overline{[u]_{D^{\overrightarrow{i}}s}}=[a\langle
\overrightarrow{m}\rangle \overline{D^{\overrightarrow{i}}s}\langle
\overrightarrow{m}'\rangle b]$.
\end{lemma}
\noindent{\bf Proof.} Assume first that $b\neq 1$. Then
$\overrightarrow{i}=(0,\cdots,0)$, $s$ is D-free, and
$[u]_s=[a\langle \overrightarrow{m}\rangle s\langle
\overrightarrow{m}'\rangle b]$. Let
$s=[\overline{s}]+\Sigma\alpha_l[u_l]$ where $[\overline{s}]$ and
$[u_l]$ are normal D-free, and $u_l<\overline{s}$. Applying Lemma
\ref{l3.2} to each of the expressions
$$
[a\langle \overrightarrow{m}\rangle [\overline{s}]\langle
\overrightarrow{m}'\rangle b] \ \ \mbox{ and } \ \ [a\langle
\overrightarrow{m}\rangle [u_l]\langle \overrightarrow{m}'\rangle
b],
$$
we see that the leading associative normal word of $[a\langle
\overrightarrow{m}\rangle s\langle \overrightarrow{m}'\rangle b]$ is
$a\langle \overrightarrow{m}\rangle \overline{s}\langle
\overrightarrow{m}'\rangle b$.
\par Next if $b=1$, then
$[u]_{D^{\overrightarrow{i}}s}=[a\langle \overrightarrow{m}\rangle
D^{\overrightarrow{i}}s]$. Again, we have
$$
s=[\overline{s}]+\Sigma\alpha_l[u_l]
$$
where $[\overline{s}]$ and $[u_l]$ are normal, $u_l<\overline{s}$.
Also,
$$
[a\langle \overrightarrow{m}\rangle
D^{\overrightarrow{i}}s]=[a\langle \overrightarrow{m}\rangle
D^{\overrightarrow{i}}[\overline{s}]] +\Sigma\alpha_l[a\langle
\overrightarrow{m}\rangle D^{\overrightarrow{i}}[u_l]],
$$
of course,
$D^{\overrightarrow{i}}[\overline{s}]=[\overline{s}D^{\overrightarrow{i}}]+\Sigma\beta_j[v_j]$
with normal $[v_j]$ such that
$v_j<\overline{s}D^{\overrightarrow{i}}$ where
$\overline{s}D^{\overrightarrow{i}}$ means that we apply
$D^{\overrightarrow{i}}$ to the last letter of $\overline{s}$. To
see this, we use induction on $|\overline{s}|$. For
$|\overline{s}|=1$, it is clear. So assume $\overline{s}=a_1\langle
\overrightarrow{m}'\rangle[s']$, $s'$ is normal,
$\overrightarrow{m}'\in Z_+^n$, $\overrightarrow{m}'\prec
\overrightarrow{N}$. Then
\begin{eqnarray*}
D^{\overrightarrow{i}}[\overline{s}]&=&D^{\overrightarrow{i}}[a_1\langle
\overrightarrow{m}'\rangle[s']] =\sum_{\overrightarrow{t}\in
Z_+^n}\binom{\overrightarrow{i}}{\overrightarrow{t}}D^{\overrightarrow{t}}a_1\langle
\overrightarrow{m}'\rangle
D^{\overrightarrow{i}-\overrightarrow{t}}[s']\\
&=&a_1\langle \overrightarrow{m}'\rangle D^{\overrightarrow{i}}[s']
+\sum_{\overrightarrow{t}\in
Z_+^n\backslash0}(-1)^{\overrightarrow{t}}\prod_{q=1}^n[m'_q]^{t_q}
\binom{\overrightarrow{i}}{\overrightarrow{t}}a_1\langle
\overrightarrow{m}'-\overrightarrow{t}\rangle
D^{\overrightarrow{i}-\overrightarrow{t}}[s']
\end{eqnarray*}
where $[m'_q]^{t_q}=m'_q(m'_q-1)\cdots(m'_q-t_q+1),1\leqslant
q\leqslant n$. By induction,
\begin{eqnarray*}
 D^{\overrightarrow{i}}[s']&=&[\overline{s'}D^{\overrightarrow{i}}]+\Sigma\alpha_j[u_j],\\
D^{\overrightarrow{i}-\overrightarrow{t}}[s']&=&
[\overline{s'}D^{\overrightarrow{i}-\overrightarrow{t}}]+\Sigma\gamma_j[w_j]
\end{eqnarray*}
with normal $[u_j]$ and $[w_j]$ such that
$u_j<\overline{s'}D^{\overrightarrow{i}},w_j<\overline{s'}D^{\overrightarrow{i}-\overrightarrow{t}}$.
So
$$
[\overline{s}D^{\overrightarrow{i}}]=a_1\langle
\overrightarrow{m}'\rangle[\overline{s'}D^{\overrightarrow{i}}], \ \
a_1\langle \overrightarrow{m}'\rangle[u_j]<a_1\langle
\overrightarrow{m}'\rangle[\overline{s'}D^{\overrightarrow{i}}],
$$
$$
a_1\langle
\overrightarrow{m}'-\overrightarrow{t}\rangle[\overline{s'}D^{\overrightarrow{i}-\overrightarrow{t}}]<a_1\langle
\overrightarrow{m}'\rangle[\overline{s'}D^{\overrightarrow{i}}],  \
\ a_1\langle
\overrightarrow{m}'-\overrightarrow{t}\rangle[w_j]<a_1\langle
\overrightarrow{m}'\rangle[\overline{s'}D^{\overrightarrow{i}}]
$$
for $\overrightarrow{t}\in Z_+^n\backslash0$. We are done.
\par Similarly, $D^{\overrightarrow{i}}[u_l]=[u_l D^{\overrightarrow{i}}]+\Sigma\beta'_j[v'_j]$
with normal $[v'_j]$ such that $v'_j<u_l D^{\overrightarrow{i}}$.
Since $u_l<\overline{s}$, $u_l
D^{\overrightarrow{i}}<\overline{s}D^{\overrightarrow{i}}$. So
$\overline{D^{\overrightarrow{i}}s}=\overline{s}D^{\overrightarrow{i}}$,
$a\langle \overrightarrow{m}\rangle
u_lD^{\overrightarrow{i}}<a\langle \overrightarrow{m}\rangle
\overline{s}D^{\overrightarrow{i}}$.
\par Thus we see that the leading associative normal word of
$a\langle \overrightarrow{m}\rangle D^{\overrightarrow{i}}s$ is
$a\langle \overrightarrow{m}\rangle
\overline{s}D^{\overrightarrow{i}}$, which is equal to $a\langle
\overrightarrow{m}\rangle \overline{D^{\overrightarrow{i}}s}$. This
completes the proof of the lemma. \ \ $\square$

\begin{definition}
Let $f$ and $g$ be monic polynomials of
$C(B,\overrightarrow{N},D_1,\cdots,D_n)$ and $w\in T$. We have the
following compositions.

$\bullet$ If $g$ is D-free and $w=\overline{f}= u\langle
\overrightarrow{m}\rangle \overline{g}\langle
\overrightarrow{m}'\rangle v$ for some $u,v\in T$ with $u$ D-free
then define
$$
(f,g)_w=f-[u\langle \overrightarrow{m}\rangle g \langle
\overrightarrow{m}'\rangle v],
$$
which is a composition of including.

$\bullet$ If $w=\overline{f}D_{j_1}^{i_1}\cdots
D_{j_k}^{i_k}=u\langle
\overrightarrow{m}\rangle\overline{g}D_{j_{k+1}}^{i_{k+1}}\cdots
D_{j_n}^{i_n}$ for some $u\in T$, $u$ D-free, $i_1,\cdots,i_n \in
Z_+$, $(j_1,\cdots,j_k,j_{k+1},\cdots,j_n)\in S_n$, then define
$$
(f,g)_w=D_{j_1}^{i_1}\cdots D_{j_k}^{i_k}f-[u\langle
\overrightarrow{m}\rangle D_{j_{k+1}}^{i_{k+1}}\cdots
D_{j_n}^{i_n}g],
$$
which is a composition of right including.

$\bullet$ If $f$ is D-free and $w= \overline{f}\langle
\overrightarrow{m}\rangle v=u\langle \overrightarrow{m}'\rangle
\overline{g}$ for some $u,v\in T$, $u$ D-free, and
$|\overline{f}|+|\overline{g}|>|w|$, then define
$$
(f,g)_w=[f\langle \overrightarrow{m}\rangle v]-[u\langle
\overrightarrow{m}'\rangle g],
$$
which is a composition of intersection.

$\bullet$ If $u\in B$ and $\overrightarrow{m}\not\prec
\overrightarrow{N}$, then $a\langle \overrightarrow{m}\rangle f$ is
referred to as a composition of left multiplication.

$\bullet$ If $f$ is not D-free, $a\in B$ and $\overrightarrow{m}\in
Z_+^n$, then $f\langle \overrightarrow{m}\rangle a$ is referred to
as a composition of right multiplication.
\end{definition}

\noindent{\bf Remark}. From Lemma \ref{l3.6}, it follows that for
any composition, $\overline{(f,g)_w}<w$.
\begin{definition}
Let $S$ be a set of polynomials in
$C(B,\overrightarrow{N},D_1,\cdots,D_n)$. A composition $(f,g)_w$ is
said to be trivial modulo $S$, denoted by $(f,g)_w\equiv0\ mod\ S$
or $(f,g)_w\equiv0\ mod(S,w)$,  if
\begin{eqnarray}\label{e3.8}
(f,g)_w=\sum_{l\in I}\alpha_l[u_l\langle
\overrightarrow{m}^{(l)}\rangle s_l\langle
\overrightarrow{m}'^{(l)}\rangle v_l] +\sum_{j\in
J}\alpha_j[u_j\langle \overrightarrow{m}^{(j)}\rangle
D^{\overrightarrow{i_j}}s_j]
\end{eqnarray}
where $I\cap J=\emptyset, s_l,s_j\in S, \ \  [u_l\langle
\overrightarrow{m}^{(l)}\rangle s_l\langle
\overrightarrow{m}'^{(l)}\rangle v_l] \ and \ [u_j\langle
\overrightarrow{m}^{(j)}\rangle D^{\overrightarrow{i_j}}s_j]
$
are normal $S$-words (in particular, $s_l, l\in I$, are D-free), and
$
u_l\langle \overrightarrow{m}^{(l)}\rangle \overline{s_l}\langle
\overrightarrow{m}'^{(l)}\rangle v_l<w, \ u_j\langle
\overrightarrow{m}^{(j)}\rangle \overline{s_j}
D^{\overrightarrow{i_j}}<w
$
for all $l\in I,j\in J$. Also, if $f$ and $g$ are both D-free, then
all normal $S$-words in (\ref{e3.8}) are D-free as well.
\par Next, $h=a\langle \overrightarrow{m}\rangle f$ or $g\langle \overrightarrow{m}'\rangle
a$, where $a\in B,\overrightarrow{m}\not\prec \overrightarrow{N}$
and $\overrightarrow{m}'\in Z_+^n$, is call trivial mod $S$ ($h
\equiv 0$ mod $S$) if $h$ can be written in the form (\ref{e3.8})
with
$$
|u_l\langle \overrightarrow{m}^{(l)}\rangle \overline{s_l}\langle
\overrightarrow{m}'^{(l)}\rangle v_l|\leqslant |\overline{h}| \ , \
|u_j\langle \overrightarrow{m}^{(j)}\rangle \overline{s_j}
D^{\overrightarrow{i_j}}|\leqslant |\overline{h}|
$$
for $l\in I,j\in J$. Also, all normal $S$-words in (\ref{e3.8}) are
D-free if $h=g\langle \overrightarrow{m}'\rangle a$ or $h=a\langle
\overrightarrow{m}\rangle f$ with $f$ D-free.
\end{definition}
\par Finally, a set $S\subseteq C(B,\overrightarrow{N},D_1,\cdots,D_n)$ of monic polynomials is
called a Gr\"{o}bner-Shirshov basis if all compositions of elements
of $S$ are trivial modulo $S$.

\ \

\noindent{\bf Remark}. When calculating compositions in what
follows, we will omit the elements the form of the right hand side
of (\ref{e3.8}). Such an omission will be signified by replacing the
symbol $=$ by the symbol $\equiv$. Thus in this notation, the
triviality of a composition $h$ (mod $S$) can be stated as $h\equiv
0$ (mod $S$).
\begin{lemma}\label{l3.9}
Let $s\in S$ be D-free and $u\in T$. Then the $S$-word $[s\langle
\overrightarrow{m}\rangle u]$ is a linear combination of normal
$S$-words $s\langle \overrightarrow{m}^{(l)}\rangle u_l$, where
$\overrightarrow{m}^{(l)}\prec \overrightarrow{N}$, $u_l\in
T,|u_l|=|u|$. If $u$ is D-free, then $u_l$'s are D-free as well.
\end{lemma}
\noindent{\bf Proof.} We use induction on $\overrightarrow{m}$. If
$\overrightarrow{m}\prec \overrightarrow{N}$, then $[s\langle
\overrightarrow{m}\rangle u]$ is already a normal $S$-word. So we
let $\overrightarrow{m}\not\prec  \overrightarrow{N}$. Assume first
that $|u|=1$, so $u=D_1^{i_1}\cdots D_n^{i_n}a$ for some $a\in B$
and $i_1,\cdots i_n\in Z_+$. When $i_1=\cdots=i_n=0$, $s\langle
\overrightarrow{m}\rangle u=s\langle \overrightarrow{m}\rangle a=0$
since $s$ is D-free, $a\in B$ and $\overrightarrow{m}\not\prec
\overrightarrow{N}$. Thus, no less of generality, assume that
$i_{l_1},\cdots,i_{l_k}>0$, $i_{l_{k+1}}=\cdots=i_{l_n}=0$, $1< k
\leqslant n$, $(l_1,\cdots,l_n) \in S_n$. Then
\begin{eqnarray*}
0&=&D_{l_1}^{i_{l_1}}\cdots D_{l_k}^{i_{l_k}}(s\langle m\rangle a)\\
&=&s\langle \overrightarrow{m}\rangle D_{l_1}^{i_{l_1}}\cdots
D_{l_k}^{i_{l_k}}a +\sum_{(s_1,\cdots,s_k)\in Z_+^k\backslash
0}\alpha_{s_1,\cdots,s_k}s\langle
\overrightarrow{m}-\Sigma_{p=1}^ks_p\overrightarrow{e_{l_p}}\rangle
D_{l_1}^{i_{l_1}-s_1}\cdots D_{l_k}^{i_{l_k}-s_k}a.
\end{eqnarray*}
Thus we have
$$
s\langle \overrightarrow{m}\rangle D_{l_1}^{i_{l_1}}\cdots
D_{l_k}^{i_{l_k}}a =\sum_{(s_1,\cdots,s_k)\in Z_+^k\backslash
0}-\alpha_{s_1,\cdots,s_k}s\langle
\overrightarrow{m}-\Sigma_{p=1}^ks_p\overrightarrow{e_{l_p}}\rangle
D_{l_1}^{i_{l_1}-s_1}\cdots D_{l_k}^{i_{l_k}-s_k}a.
$$
The case that $|u|=1$ follows from the induction on
$(i_1,\cdots,i_n)$.
\par Now, let $|u|>1$, and write
$$
u=a_1\langle \overrightarrow{m}'\rangle u_1, \ \overrightarrow{m}\in
Z_+^n, \ \overrightarrow{m}\prec \overrightarrow{N}, \ a_1\in B.
$$
Then
$$
s\langle \overrightarrow{m}\rangle u=
s\langle \overrightarrow{m}\rangle (a_1\langle \overrightarrow{m}'\rangle [u_1])
=(s\langle \overrightarrow{m}\rangle a_1)\langle \overrightarrow{m}'\rangle [u_1] -\sum_{\overrightarrow{t}\in
Z_+^n\backslash 0}\beta_{\overrightarrow{t}} s\langle \overrightarrow{m}-\overrightarrow{t}\rangle (a_1\langle
\overrightarrow{m}'+\overrightarrow{t}\rangle [u_1])
$$
where $s\langle \overrightarrow{m}\rangle a_1=0$ as we have seen. By
Lemma \ref{l2.2}, all $a_1\langle
\overrightarrow{m}'+\overrightarrow{t}\rangle [u_1]$ are linear
combinations of normal words of the same length. Then
$$
s\langle \overrightarrow{m}\rangle u=
\sum\gamma_{\overrightarrow{t},j}s\langle \overrightarrow{m}-\overrightarrow{t}\rangle[v_j],\
v_j\in T, \ |v_j|=|v|.
$$
Hence, the result follows from the induction on
$\overrightarrow{m}$. \ \  $\square$

\begin{lemma}\label{l3.10}
Let $S$ be a Gr\"{o}bner-Shirshov basis, $s\in S$ not D-free, and
$u\in T$. Then the $S$-word $[s\langle \overrightarrow{m}\rangle u]$
of formal degree $L$ is a linear combination of normal $S$-words of
degree at most $L$. If $u$ is D-free, then the normal $S$-words in
consideration are D-free as well.
\end{lemma}
\noindent{\bf Proof.} Suppose $|u|=1$. Then $u=D_1^{i_1}\cdots
D_n^{i_n}a$ for some $a\in B$. If $i_1=\cdots=i_n=0$, then $u=a$ and
the statement follows from the triviality of the right
multiplication. No less of generality, assume that $i_k>0$ for some
$k,\ 1\leqslant k\leqslant n$. Then
$$
s\langle \overrightarrow{m}\rangle D_1^{i_1}\cdots
D_n^{i_n}a=D_k(s\langle \overrightarrow{m}\rangle D_1^{i_1}\cdots
D_k^{i_k-1}\cdots D_n^{i_n}a) +  m_k s\langle
\overrightarrow{m}-\overrightarrow{e_k}\rangle D_1^{i_1}\cdots
D_k^{i_k-1}\cdots D_n^{i_n}a
$$
and the results follows from the fact that the derivation of any
normal $S$-word is a linear combination of normal $S$-words, and the
induction on $(i_1,\cdots,i_n)$.
\par Now, assume that $|u|>1$. Hence $u=a_1\langle
\overrightarrow{m}'\rangle u_1$ for some $a_1\in B$ and $u_1\in T$. Then
$$
[s\langle \overrightarrow{m}\rangle u]=
s\langle \overrightarrow{m}\rangle (a_1\langle \overrightarrow{m}'\rangle [u_1])
=\sum_{\overrightarrow{t}\in
   Z_+^n}\binom{\overrightarrow{m}}{\overrightarrow{t}}(s\langle
   \overrightarrow{m}-\overrightarrow{t}\rangle a_1)\langle
   \overrightarrow{m}'+\overrightarrow{t}\rangle [u_1].
$$
\par Again, using the triviality of the right multiplication, we
have that every $s\langle \overrightarrow{m}-\overrightarrow{t}\rangle a_1$ is linear combination of
D-free normal $S$-words of degree at most $L-|u_1|$ of the form
\begin{eqnarray}\label{e3.9}
[v\langle \overrightarrow{l}\rangle s_1\langle
\overrightarrow{t}\rangle w]
\end{eqnarray}
where $s_1\in S$ and $v,w\in T$ (which may be empty). Multiplying
(\ref{e3.9}) by some $\langle \overrightarrow{k}\rangle[u_1]$ from
the right, and applying the associative law, we obtain a linear
combination of words of the form
\[
\begin{cases}
[v'\langle \overrightarrow{l}'\rangle s_1\langle \overrightarrow{t}'\rangle
w'\langle \overrightarrow{k}'\rangle u_1],&\text {if $w\neq 1$}\\
[v'\langle \overrightarrow{l}'\rangle s_1\langle \overrightarrow{k}'\rangle u_1],&\text {if $w=1$}\\
\end{cases}
\]
where $\overrightarrow{l}'\preceq \overrightarrow{l}\prec
\overrightarrow{N}$,
 $\overrightarrow{t}'\preceq \overrightarrow{t}\prec \overrightarrow{N}$, $v',w'\in T$, and
$k\preceq k'$. We can assume, using Lemma \ref{l2.2}, that in the
words of the first form, all $w'\langle \overrightarrow{k}'\rangle
u_1$, are normal, and also, since $s_1$ is D-free now, Lemma
\ref{l3.9} can be applied to the words of the second form. In all
cases, we have linear combinations of normal $S$-words of degrees at
most $L$. \ \ $\square$

\ \

The proof of the next lemma is  similar to that of Lemma \ref{l2.1}.
The only difference is that we add to $B$ a new letter
$D^{\overrightarrow{i}}s$, and we are dealing with words with one
occurrence of $D^{\overrightarrow{i}}s$.
\begin{lemma}\label{l3.11}
Any $S$-word $(u)_{D^{\overrightarrow{i}}s}$ can be presented as a
linear combination of right normed $S$-words
\begin{eqnarray}\label{e3.10}
[v]_{D^{\overrightarrow{i}}s}=x_1\langle
\overrightarrow{m}^{(1)}\rangle(x_2\langle
\overrightarrow{m}^{(2)}\rangle\cdots(x_k\langle
\overrightarrow{m}^{(k)}\rangle x_{k+1})\cdots)
\end{eqnarray}
where $x_{j_0}=D^{\overrightarrow{i}}s$ for some $j_0\in \{1,\cdots
k+1 \}$, and $x_j\in D^\omega(B)$ for $j\neq j_0$, each of the same
formal degree.
\end{lemma}

\begin{lemma}\label{l3.12}
Let $S$ be a Gr\"{o}bner-Shirshov basis in
$C(B,\overrightarrow{N},D_1,\cdots,D_n)$. Then any $S$-word
$(v)_{D^{\overrightarrow{i}}s}$ of formal degree $L$ is a linear
combination of normal $S$-words of degree less or equal to $L$. If
$i_1=\cdots=i_n=0$, then all the normal $S$-words are D-free as
well.
\end{lemma}
\noindent{\bf Proof.} By Lemma \ref{l3.11}, we may assume that
$(v)_{D^{\overrightarrow{i}}s}=[v]_{D^{\overrightarrow{i}}s}$ of the
form (\ref{e3.10}). If $L=1$, then
$(v)_{D^{\overrightarrow{i}}s}=D^{\overrightarrow{i}}s$, and it is a
normal $S$-word.
\par Now, let $L>1$. Consider first that $|v|_{D^{\overrightarrow{i}}s}=2$. There are essentially two possibilities:
$$
[v]_s=s\langle \overrightarrow{m}\rangle D^{\overrightarrow{j}}a \ \
or \ \ [v]_{D^{\overrightarrow{i}}s}=a\langle
\overrightarrow{m}\rangle D^{\overrightarrow{i}}s
$$
where $s\in S,\ a\in B$. In the first case, the result follows from
Lemma \ref{l3.9} or Lemma \ref{l3.10}.
\par Suppose $i_1=\cdots=i_n=0$ in $a\langle
\overrightarrow{m}\rangle D^{\overrightarrow{i}}s$, so
$[v_s]=a\langle \overrightarrow{m}\rangle s$. If
$\overrightarrow{m}\prec \overrightarrow{N}$, it is a normal
$S$-word. If $\overrightarrow{m}\not\prec  \overrightarrow{N}$, we
are done by the triviality of left multiplication.

Assume that $i_k>0$ for some $k, \ 1\leqslant k\leqslant n$, in
$a\langle \overrightarrow{m}\rangle D^{\overrightarrow{i}}s$, and
that all
$$
a\langle \overrightarrow{m}\rangle D_1^{i_1}\cdots
D_k^{i_k-1}\cdots D_n^{i_n}s
$$
are linear combinations of normal $S$-words of degrees at most $L$.
Now the result follows from the fact that
$$
a\langle \overrightarrow{m}\rangle D_1^{i_1}\cdots D_n^{i_n}s=D_k(a\langle
\overrightarrow{m}\rangle D_1^{i_1}\cdots D_k^{i_k-1}\cdots D_n^{i_n}s)
 + m_k a\langle m-e_k\rangle
D_1^{i_1}\cdots D_k^{i_k-1}\cdots D_n^{i_n}s,
$$
the derivation of any normal $S$-word is a linear combination of
normal $S$-words and the induction on $(i_1,\cdots,i_n)$.

Now, let $|v|_{D^{\overrightarrow{i}}s}\geqslant 3$. If
$(v)_{D^{\overrightarrow{i}}s}$ begins with
$D^{\overrightarrow{i}}s$, then we may assume that
$i_1=\cdots=i_n=0$, $[v]_s=s\langle \overrightarrow{m}\rangle[u]$,
$s\in S$, $\overrightarrow{m}\not\prec  \overrightarrow{N}$, and
$[u]$ is a normal word. Then, by Lemma \ref{l3.9} or Lemma
\ref{l3.10} we are done.
\par Suppose that $[v]_{D^{\overrightarrow{i}}s}=
a\langle \overrightarrow{m}\rangle[u]_{D^{\overrightarrow{i}}s}$ for
some $a\in B$ and $\overrightarrow{m}\in Z_+^n$. We may assume that
$[u]_{D^{\overrightarrow{i}}s}$ is a normal $S$-word of degree less
than  or equal to $L-1$. If $\overrightarrow{m}\prec
\overrightarrow{N}$, then $a\langle
\overrightarrow{m}\rangle[u]_{D^{\overrightarrow{i}}s}$ is already a
normal $S$-word. So let $\overrightarrow{m}\not\prec
\overrightarrow{N}$. Applying the associativity law to $a\langle
\overrightarrow{m}\rangle[u]_{D^{\overrightarrow{i}}s}$, we get a
linear combination of four kinds of words (of formal degree less
than or equal to $L$):
\begin{eqnarray}\label{e3.12}
&&\nonumber a\langle
\overrightarrow{m}-\overrightarrow{t}\rangle[u]_{D^{\overrightarrow{i}}s}
\ \ (\overrightarrow{t}\in Z_+^n\backslash0), \ \ (a\langle
\overrightarrow{m}\rangle b)\langle
\overrightarrow{m}'\rangle[u]_{D^{\overrightarrow{i}}s},\\
&& a\langle \overrightarrow{m}\rangle D^{\overrightarrow{i}}s, \ \
(a\langle \overrightarrow{m}\rangle s)\langle
\overrightarrow{m}'\rangle[u_1]
\end{eqnarray}
where $a,b\in B$. Induction on $\overrightarrow{m}$ is applicable to
the words of the first kind. Since $a\langle
\overrightarrow{m}\rangle b=0$, the words of the second kind are
zero. The words of the third type in (\ref{e3.12}) have been treated
above. Note that, by triviality of left multiplications, the words
of the fourth kind are linear combination of words
\begin{eqnarray}\label{e3.13}
[v_l\langle \overrightarrow{m}^{(l)}\rangle s_l\langle
\overrightarrow{t}^{(l)}\rangle w_l]\langle
\overrightarrow{m}'\rangle[u_1], \ \ [v_j\langle
\overrightarrow{m}^{(j)}\rangle
D^{{\overrightarrow{i}}_j}s_j]\langle
\overrightarrow{m}'\rangle[u_1]
\end{eqnarray}
where $[v_l\langle \overrightarrow{m}^{(l)}\rangle s_l\langle
\overrightarrow{t}^{(l)}\rangle w_l]$ are normal $S$-words of the
first kind, $[v_j\langle \overrightarrow{m}^{(j)}\rangle
D^{\overrightarrow{i_j}}s_j]$ are normal $S$-words of the second
kind, and all words in (\ref{e3.13}) have formal degree at most $L$.
\par Applying the associative law to (\ref{e3.13}), we get a linear
combination of right normed $S$-words, of formal degree less than or
equal to $L$ of the form
$$
[v'_l\langle \overrightarrow{m}'^{(l)}\rangle s_l\langle
\overrightarrow{t}'^{(l)}\rangle w'_l\langle
\overrightarrow{m}''\rangle u_1], \ \ [v'_j\langle
\overrightarrow{m}'^{(j)}\rangle s_j\langle
\overrightarrow{m}'''\rangle u_1]
$$
where $s_l,s_j\in S$, $v'_l,w'_l$ and $v'_j$ are normal D-free
associative words ($w'_l$ may be empty),
$\overrightarrow{m}'^{(l)}\preceq \overrightarrow{m}^{(l)}\prec
\overrightarrow{N}$, $\overrightarrow{t}'^{(l)}\preceq
\overrightarrow{t}^{(l)}\prec \overrightarrow{N}$, and
$\overrightarrow{m}''\in Z_+^n$, $\overrightarrow{m}'^{(j)}\preceq
\overrightarrow{m}^{(j)}\prec \overrightarrow{N}$ and
$\overrightarrow{m}'''\in Z_+^n$. Now the result follows from Lemmas
\ref{l2.2}, \ref{l3.9} and \ref{l3.10} the same way as before.
\par In the above argument, we have used the compositions of left
and right multiplications. If $i_1=\cdots=i_n=0$ and $s$ is D-free,
we only need compositions of left multiplications. In this case, we
get only D-free normal $S$-words. \ \  $\square$

\ \

Here and after, we say that a polynomial $f$ is smaller than another
$g$ if $\overline{f}<\overline{g}$.

\begin{lemma}\label{l3.13}
Let $S$ be a Gr\"{o}bner-Shirshov basis. Let
\begin{eqnarray}\label{e3.15}
(r)_s=[p\langle \overrightarrow{i}\rangle
 s\langle \overrightarrow{j}\rangle q]\langle \overrightarrow{k}\rangle[t]
\end{eqnarray}
be an $S$-word where $[p\langle \overrightarrow{i}\rangle s\langle
\overrightarrow{j}\rangle q]$ is a normal $S$-word of the first
kind, $q$ is D-free, $t\in T$, and $\overrightarrow{k}\in Z_+^n$,
$\overrightarrow{k}\prec \overrightarrow{N}$. Then $(r)_s$ is equal
to $[p\langle \overrightarrow{i}\rangle s\langle
\overrightarrow{j}\rangle q \langle \overrightarrow{k}\rangle t]$
modulo smaller normal $S$-words of the first kind. If $t$ is D-free,
then all the words are D-free as well.
\end{lemma}
\noindent{\bf Proof.} Consider first that $|p|=0$. If $|q|=0$, then
$s\langle \overrightarrow{k}\rangle t=[s\langle
\overrightarrow{k}\rangle t]$ is already a normal $S$-word. So we
assume that $|q|>0$. Then
$$
(s\langle \overrightarrow{j}\rangle [q])\langle \overrightarrow{k}\rangle[t]=s\langle \overrightarrow{j}\rangle
([q]\langle \overrightarrow{k}\rangle[t]) +\sum_{\overrightarrow{l}\in
   Z_+^n\backslash0}\alpha_{\overrightarrow{l}} s\langle
   \overrightarrow{j}-\overrightarrow{l}\rangle ([q]\langle
   \overrightarrow{k}+\overrightarrow{l}\rangle [t]).
$$
By Lemma \ref{l3.1}, $[q]\langle \overrightarrow{k}\rangle[t]$ is
equal to $[q\langle \overrightarrow{k}\rangle t]$ modulo smaller
normal words. So $s\langle \overrightarrow{j}\rangle ([q]\langle
\overrightarrow{k}\rangle[t])$ is equal to $s\langle
\overrightarrow{j}\rangle [q\langle \overrightarrow{k}\rangle
t]=[s\langle \overrightarrow{j}\rangle q\langle
\overrightarrow{k}\rangle t]$ modulo smaller normal $S$-words. By
Lemma \ref{l2.2}, $[q]\langle
   \overrightarrow{k}+\overrightarrow{l}\rangle [t]$ is a linear combination of
   normal words $[u_i]$ of the same length. It follows that $s\langle
   \overrightarrow{j}-\overrightarrow{l}\rangle
   [u_i]$, $\overrightarrow{l}\in
   Z_+^n\backslash0$, are smaller normal $S$-words.
\par Now, let $|p|\geqslant 1$, and $p=a_1\langle \overrightarrow{m}\rangle
p_1$ where $a_1\in B$, $\overrightarrow{m}\prec \overrightarrow{N}$
and $|p_1|\geqslant 0$. Denote for simplicity that
$[h]_s=[p_1\langle \overrightarrow{i}\rangle s\langle
\overrightarrow{j}\rangle q]$. Then (\ref{e3.15}) has the form
$$
[p\langle \overrightarrow{i}\rangle s\langle \overrightarrow{j}\rangle q]\langle
\overrightarrow{k}\rangle[t]=(a_1\langle \overrightarrow{m}\rangle[h]_s)\langle \overrightarrow{k}\rangle[t]
\ \ \ \ \ \ \ \ \  \ \ \ \ \ \ \ \ \  \ \ \ \ \ \ \ \ \ \ \ \ \
\ \ \ \ \ \ \  \ \ \ \ \ \ \ \ \  \ \ \ \ \
$$
$$
\ \ \ \ \ \ \ \ \  \ \ \ \ \ \ \ \ \  \ \ \ \ \ \ \ \ \ \ \
=a_1\langle \overrightarrow{m}\rangle([h]_s\langle \overrightarrow{k}\rangle[t]) +\sum_{\overrightarrow{l}\in
   Z_+^n\backslash0}\beta_{\overrightarrow{l}} a_1\langle
   \overrightarrow{m}-\overrightarrow{l}\rangle ([h]_s\langle
   \overrightarrow{k}+\overrightarrow{l}\rangle [t]).
$$
Since $|p_1|<|p|$, the induction on $|p|$ applies and we have that
$a_1\langle \overrightarrow{m}\rangle([h]_s\langle
\overrightarrow{k} \rangle[t])$ is equal to $a_1\langle
\overrightarrow{m}\rangle[h_s\langle \overrightarrow{k}\rangle t]=
[a_1\langle \overrightarrow{m}\rangle h_s\langle
\overrightarrow{k}\rangle t]$ modulo smaller normal $S$-words. On
the other hand, let $L=$fdeg$(r)_s$. Then by Lemma \ref{l3.12},
$[h]_s\langle \overrightarrow{k}+\overrightarrow{l}\rangle [t]$ is a
linear combination of normal $S$-words of the first kind whose
degree is at most $L-1$. Hence $a_1\langle
\overrightarrow{m}-\overrightarrow{l}\rangle ([h]_s\langle
\overrightarrow{k}+\overrightarrow{l}\rangle [t])$ is a linear
combination of smaller normal $S$-words of the first kind. This
completes the proof of the lemma. \ \  $\square$

\begin{lemma}\label{l3.14}
Let $S$ be a Gr\"{o}bner-Shirshov basis. Let
$$
(r)_s=[t]\langle \overrightarrow{k}\rangle[p\langle \overrightarrow{i}\rangle s\langle \overrightarrow{j}\rangle q]
$$
be an $S$-word where $[p\langle \overrightarrow{i}\rangle s\langle
\overrightarrow{j}\rangle q]$ is a normal $S$-word of the first
kind, $t\in T$ is D-free, and $\overrightarrow{k}\in Z_+^n$,
$\overrightarrow{k}\prec \overrightarrow{N}$. Then $(r)_s$ is equal
to $[t\langle \overrightarrow{k}\rangle p\langle
\overrightarrow{i}\rangle s\langle \overrightarrow{j}\rangle q ]$
modulo smaller normal $S$-words of the first kind. If $q$ is D-free,
then all the words are D-free as well.
\end{lemma}
\noindent{\bf Proof.} We shall use induction on $|t|$. When $|t|=1,
t=a\in B$ and $(r)_s$ is a normal $S$-word. So, let $|t|>1$ and
write $t=a_1\langle \overrightarrow{m}\rangle t_1$. For simplicity,
we put $[h]_s=[p\langle \overrightarrow{i}\rangle s\langle
\overrightarrow{j}\rangle q]$. Then
$$
(r)_s=(a_1\langle \overrightarrow{m}\rangle [t_1])\langle \overrightarrow{k}\rangle[h]_s =a_1\langle
\overrightarrow{m}\rangle ([t_1]\langle \overrightarrow{k}\rangle[h]_s) +\sum_{\overrightarrow{l}\in
   Z_+^n\backslash0}\alpha_{\overrightarrow{l}} a_1\langle
   \overrightarrow{m}-\overrightarrow{l}\rangle ([t_1]\langle
   \overrightarrow{k}+\overrightarrow{l}\rangle [h]_s).
$$
We can apply induction on $t_1\langle
\overrightarrow{k}\rangle[h]_s$ since $|t_1|<|t|$. As a result,
$a_1\langle \overrightarrow{m}\rangle (t_1\langle
\overrightarrow{k}\rangle[h]_s)$ is equal to $a_1\langle
\overrightarrow{m}\rangle [t_1\langle \overrightarrow{k}\rangle
p\langle \overrightarrow{i}\rangle s\langle
\overrightarrow{j}\rangle q]= [t\langle \overrightarrow{k}\rangle
p\langle \overrightarrow{i}\rangle s\langle
\overrightarrow{j}\rangle q]$ modulo smaller normal $S$-words of the
first kind. By Lemma \ref{l3.12}, $[t_1]\langle
   \overrightarrow{k}+\overrightarrow{l}\rangle [h]_s$ is a linear combination of
   normal $S$-words of the first kind whose degree is at most
   fdeg$(r)_s-1$. Thus, each $a_1\langle
   \overrightarrow{m}-\overrightarrow{l}\rangle ([t_1]\langle
   \overrightarrow{k}+\overrightarrow{l}\rangle [h]_s)$ is a linear combination of
   normal $S$-words of the first kind which are smaller than $[t\langle
\overrightarrow{k}\rangle p\langle \overrightarrow{i}\rangle s\langle \overrightarrow{j}\rangle q ]$. \ \ $\square$

\begin{proposition}\label{p3.15}
Let $S$ be a Gr\"{o}bner-Shirshov basis in
$C(B,\overrightarrow{N},D_1,\cdots,D_n)$, and let $u_s=v\langle
\overrightarrow{m}\rangle s\langle \overrightarrow{m}'\rangle w$ be
an associative normal $S$-word of the first kind. Let $(u)_s$ be
some bracketing of $u_s$. Then
$$
(u)_s=[v\langle \overrightarrow{m}\rangle s\langle
\overrightarrow{m}'\rangle w]+\Sigma\alpha_i[u_i]_{s_i}
$$
where $[u_i]_{s_i}$ are normal $S$-words of the first kind, and
$\overline{u_i}<\overline{u}=v\langle \overrightarrow{m}\rangle \overline{s}\langle
\overrightarrow{m}'\rangle w$. If $w$ is D-free (or empty), then all the $u_i$'s are
D-free as well.
\end{proposition}
\noindent{\bf Proof.} Let $L=$fdeg$(u)_s\geqslant 1$. If $L=1$, then
$|u|_s=1$ and $(u)_s=s$ is a normal $S$-word. Let $L>1$, $|u|_s>1$,
and
$$
(u)_s=(u_1)_s\langle \overrightarrow{k}\rangle(u_2) \ \ or \ \ (u)_s=(u_1)\langle
\overrightarrow{k}\rangle(u_2)_s
$$
for some $\overrightarrow{k}\in Z_+^n$, $\overrightarrow{k}\prec
\overrightarrow{N}$.
\par Suppose $(u)_s=(u_1)_s\langle \overrightarrow{k}\rangle(u_2)$. Then
$u_1=v\langle \overrightarrow{m}\rangle s\langle
\overrightarrow{m}'\rangle p$ and $w=p\langle
\overrightarrow{k}\rangle u_2$ where $p$ is D-free and may be empty.
By induction on $|u|$, we have
$$
(u_1)_s=[v\langle \overrightarrow{m}\rangle s\langle
\overrightarrow{m}'\rangle p]+\Sigma\alpha_i[v_i]_{s_i}
$$
where $[v_i]_{s_i}$ are D-free normal $S$-words of the first kind,
and $\overline{v_i}<v\langle \overrightarrow{m}\rangle
\overline{s}\langle \overrightarrow{m}'\rangle p$. By Lemma
\ref{l3.1},
$$
(u_2)=[u_2]+\Sigma\beta_j[w_j]
$$
where $[w_j]$ are normal words with $w_j<u_2$. Thus
$$
(u)_s=[v\langle \overrightarrow{m}\rangle s\langle \overrightarrow{m}'\rangle p]\langle
\overrightarrow{k}\rangle[u_2]+\Sigma\alpha_i[v_i]_{s_i}\langle \overrightarrow{k}\rangle[u_2]
+\Sigma\beta_j[v\langle \overrightarrow{m}\rangle s\langle \overrightarrow{m}'\rangle p]\langle
\overrightarrow{k}\rangle[w_j]+\Sigma\alpha_i\beta_j[v_i]_{s_i}\langle \overrightarrow{k}\rangle[w_j].
$$
Now, the claim follows from Lemma \ref{l3.13}.
\par Similarly, the case when $(u)_s=(u_1)\langle
\overrightarrow{k}\rangle(u_2)_s$ follows from Lemma \ref{l3.14}.
This completes the proof. \ \  $\square$

\begin{lemma}\label{l3.16}
Let $S$ be a Gr\"{o}bner-Shirshov basis in $C(B,\overrightarrow{N},D_1,\cdots,D_n)$. Let
$$
(r)_{D^{\overrightarrow{l}}s}=
[t]\langle \overrightarrow{k}\rangle[p\langle \overrightarrow{i}\rangle D^{\overrightarrow{l}}s]
$$
be an $S$-word where $[p\langle \overrightarrow{i}\rangle D^ls]$ is
a normal $S$-word of the second kind, $t\in T$ is D-free, and
$\overrightarrow{k}\in Z_+^n$, $\overrightarrow{k}\prec
\overrightarrow{N}$. Then $(r)_{D^{\overrightarrow{l}}s}$ is equal
to the normal $S$-word $[t\langle \overrightarrow{k}\rangle p\langle
\overrightarrow{i}\rangle D^{\overrightarrow{l}}s]$ modulo smaller
normal $S$-words.
\end{lemma}
\noindent{\bf Proof.} We shall use induction on $|t|$. When $|t|=1,
t=a\in B$ and $(r)_{D^{\overrightarrow{l}}s}$ is a normal $S$-word.
Let $|t|>1$ and write $t=a_1\langle \overrightarrow{m}\rangle t_1$.
For simplicity, we put $[h]_{D^{\overrightarrow{l}}s}=[p\langle
\overrightarrow{i}\rangle D^{\overrightarrow{l}}s]$. Then
$$
(r)_{D^{\overrightarrow{l}}s}=
(a_1\langle \overrightarrow{m}\rangle [t_1])\langle \overrightarrow{k}\rangle[h]_{D^{\overrightarrow{l}}s}
\ \ \ \ \ \ \ \ \ \ \ \ \ \ \ \ \ \ \  \ \ \ \ \ \ \ \ \  \ \ \ \ \ \ \ \ \ \ \ \  \ \ \ \ \ \ \ \ \
$$
$$
 \ \ \ \ \ \ \ \ \ \ \  \ \ \ \ \ \ \ \ \ =a_1\langle \overrightarrow{m}\rangle
([t_1]\langle \overrightarrow{k}\rangle[h]_{D^{\overrightarrow{l}}s})
+\sum_{\overrightarrow{j}\in
   Z_+^n\backslash0}\alpha_{\overrightarrow{j}} a_1\langle
   \overrightarrow{m}-\overrightarrow{j}\rangle ([t_1]\langle
   \overrightarrow{k}+\overrightarrow{j}\rangle [h]_{D^{\overrightarrow{l}}s}).
$$
We are done by induction on $|t|$ and Lemma \ref{l3.12}. \ \
$\square$

\begin{proposition}\label{p3.17}
Let $S$ be a Gr\"{o}bner-Shirshov basis in
$C(B,\overrightarrow{N},D_1,\cdots,D_n)$, and let
$u_{D^{\overrightarrow{l}}s}=v\langle \overrightarrow{m}\rangle
D^{\overrightarrow{l}}s$ be an associative normal $S$-word of the
second kind. Let $(u)_{D^{\overrightarrow{l}}s}$ be some bracketing
of $u_{D^{\overrightarrow{l}}s}$. Then
$$
(u)_{D^{\overrightarrow{l}}s}=[v\langle \overrightarrow{m}\rangle
D^{\overrightarrow{l}}s]+\Sigma\alpha_j[u_j]_{D^{\overrightarrow{l}_k}s_k}
$$
where $[u_j]_{D^{\overrightarrow{l}_k}s_k}$ are normal $S$-words with
$\overline{u_j}<\overline{u}=v\langle \overrightarrow{m}\rangle \overline{s}D^{\overrightarrow{l}}$.
\end{proposition}
\noindent{\bf Proof.} Let $L=$fdeg$(u)_{D^{\overrightarrow{l}}s}\geqslant 1$. If $L=1$,
then $(u)_{D^{\overrightarrow{l}}s}=D^{\overrightarrow{l}}s$ is a normal $S$-word. Let $L>1$,
$|u|_{D^{\overrightarrow{l}}s}>1$, and
$$
(u)_{D^{\overrightarrow{l}}s}=(u_1)\langle \overrightarrow{k}\rangle(u_2)_{D^{\overrightarrow{l}}s}
$$
for some $\overrightarrow{k}\in Z_+^n$, $\overrightarrow{k}\prec
\overrightarrow{N}$.
\par Suppose $u_2=p\langle \overrightarrow{m}\rangle D^{\overrightarrow{l}}s$
and $v=u_1\langle \overrightarrow{k}\rangle p$, where $p$ is D-free and may be empty.
By induction on $|u|_{D^{\overrightarrow{l}}s}$, we have
$$
(u_2)_{D^{\overrightarrow{l}}s}=[p\langle \overrightarrow{m}\rangle
D^{\overrightarrow{l}}s]+\Sigma\alpha_j[v_j]_{D^{\overrightarrow{l}_k}s_k}
$$
where $[v_j]_{D^{\overrightarrow{l}_k}s_k}$ are normal $S$-words ,
and $\overline{v_j}<p\langle \overrightarrow{m}\rangle
\overline{s}D^{\overrightarrow{l}}$. By Lemma \ref{l3.1},
$$
(u_1)=[u_1]+\Sigma\beta_j[w_j]
$$
where $[w_j]$ are normal words with $w_j<u_1$. Thus
\begin{eqnarray*}
(u)_{D^{\overrightarrow{l}}s}&=&[u_1]\langle
\overrightarrow{k}\rangle[p\langle \overrightarrow{m}\rangle
D^{\overrightarrow{l}}s]+\Sigma\alpha_j[u_1]\langle
\overrightarrow{k}\rangle[v_j]_{D^{\overrightarrow{l}_k}s_k}\\
&& +\Sigma\beta_j[w_j]\langle \overrightarrow{k}\rangle[p\langle
\overrightarrow{m}\rangle
D^{\overrightarrow{l}}s]+\Sigma\alpha_j\beta_j[w_j]\langle
\overrightarrow{k}\rangle[v_j]_{D^{\overrightarrow{l}_k}s_k}.
\end{eqnarray*}
Now, the claim follows from Lemmas \ref{l3.16} and
\ref{l3.14}. This completes the proof. \ \  $\square$

\begin{lemma}\label{l3.18}
Let $S$ be a set of monic polynomial, $s_1,s_2\in S$. Let
$(s_1,s_2)_w$ be one of the following composition:
\par $\bullet$ $(s_1,s_2)_w=s_1-[u\langle \overrightarrow{m}\rangle
s_2]$, where $w=\overline{s_1}=u\langle \overrightarrow{m}\rangle \overline{s_2}$;
\par $\bullet$ $(s_1,s_2)_w=s_1-[u\langle \overrightarrow{m}\rangle
s_2\langle \overrightarrow{m}'\rangle v]$, where $w=\overline{s_1}=u\langle \overrightarrow{m}\rangle
\overline{s_2}\langle \overrightarrow{m}'\rangle v$, $v\neq 1$, and $s_2$ is D-free;
\par $\bullet$ $(s_1,s_2)_w=D_{l_1}^{i_1}\cdots D_{l_k}^{i_k}s_1-
[u\langle \overrightarrow{m}\rangle D_{l_{k+1}}^{i_{k+1}}\cdots D_{l_n}^{i_n}s_2]$,
where $w=\overline{s_1}D_{l_1}^{i_1}\cdots D_{l_k}^{i_k}=u\langle
\overrightarrow{m}\rangle\overline{s_2}D_{l_{k+1}}^{i_{k+1}}$ $\cdots D_{l_n}^{i_n}$,
$(l_1,\cdots,l_k,l_{k+1},\cdots,l_n)\in S_n$;
\par $\bullet$ $(s_1,s_2)_w=[s_1\langle \overrightarrow{m}\rangle v]-[u\langle \overrightarrow{m}'\rangle s_2]$,
where $w=\overline{s_1}\langle \overrightarrow{m}\rangle v=u\langle \overrightarrow{m}'\rangle
\overline{s_2}$, $v\neq 1$, $s_1$ is D-free, and
$|\overline{s_1}|+|\overline{s_2}|>|w|$.

If $(s_1,s_2)_w\equiv 0\ mod(S,w)$, then we have, respectively,
\par $\bullet$ $D^{\overrightarrow{j}}s_1-[u\langle \overrightarrow{m}\rangle D^{\overrightarrow{j}}s_2]\equiv 0$
mod $(S,wD^{\overrightarrow{j}})$;
\par $\bullet$ $D^{\overrightarrow{j}}s_1-[u\langle \overrightarrow{m}\rangle s_2\langle
\overrightarrow{m}'\rangle vD^{\overrightarrow{j}}]\equiv 0$ mod $(S,wD^{\overrightarrow{j}})$;
\par $\bullet$ $D^{\overrightarrow{j}}(D_{l_1}^{i_1}\cdots D_{l_k}^{i_k}s_1)-
[u\langle \overrightarrow{m}\rangle D^{\overrightarrow{j}}
(D_{l_{k+1}}^{i_{k+1}}\cdots D_{l_n}^{i_n}s_2)]\equiv 0$ mod
$(S,wD^{\overrightarrow{j}})$;
\par $\bullet$ $[s_1\langle \overrightarrow{m}\rangle vD^{\overrightarrow{j}}]-
[u\langle \overrightarrow{m}'\rangle D^{\overrightarrow{j}}s_2]\equiv 0$ mod $(S,wD^{\overrightarrow{j}})$.
\end{lemma}
\noindent{\bf Proof.} The fact that $(s_1,s_2)_w\equiv 0\ mod(S,w)$
means that
$$
(s_1,s_2)_w=\sum_{i\in I}\alpha_i[u_i\langle
\overrightarrow{m}^{(i)}\rangle s_i\langle
\overrightarrow{m}'^{(i)}\rangle v_i] +\sum_{k\in
J}\alpha_k[u_k\langle \overrightarrow{m}^{(k)}\rangle
D^{\overrightarrow{l_k}}s_k]
$$
where all the $S$-words on the right hand side are normal and less
than $w$. We shall show that
\begin{eqnarray}\label{e3.16}
D^{\overrightarrow{j}}[a\langle \overrightarrow{m}\rangle s_1]
\equiv[a\langle \overrightarrow{m}\rangle D^{\overrightarrow{j}}s_1]
\ mod \ (S,wD^{\overrightarrow{j}})
\end{eqnarray}
and
\begin{eqnarray}\label{e3.17}
D^{\overrightarrow{j}}[a\langle \overrightarrow{m}\rangle s_2\langle
\overrightarrow{m}'\rangle b]\equiv[a\langle
\overrightarrow{m}\rangle s_2\langle \overrightarrow{m}'\rangle
bD^{\overrightarrow{j}}] \ mod \
(S,wD^{\overrightarrow{j}})
\end{eqnarray}
 for all
$\overrightarrow{j}\in Z_+^n$, $a,b\in T$, where $s_2$ is D-free and
$b\neq 1$.

 Now we use an induction on $\overrightarrow{j}$ to show
that
$$
D^{\overrightarrow{j}}[a\langle \overrightarrow{m}\rangle s_1]=
[a\langle \overrightarrow{m}\rangle D^{\overrightarrow{j}}s_1]+\varepsilon_1
$$
where $\varepsilon_1$ is a linearly combination of normal $S$-words
that less than $a\langle \overrightarrow{m}\rangle
\overline{s_1}D^{\overrightarrow{j}}$. If $j_1=\cdots=j_n=0$, it is
clear. Suppose that $j_t>0$ for some $t, \ 1\leqslant t \leqslant
n$. Consider first that $|a|=1$. Then $a\in B$ and
\begin{eqnarray*}
D^{\overrightarrow{j}}[a\langle \overrightarrow{m}\rangle s_1] &=&
D^{\overrightarrow{j}-\overrightarrow{e_t}}[D_ta\langle
\overrightarrow{m}\rangle s_1+a\langle \overrightarrow{m}\rangle
D_ts_1]\\
& =&D^{\overrightarrow{j}-\overrightarrow{e_t}} [\alpha a\langle
\overrightarrow{m}-\overrightarrow{e_t}\rangle s_1+a\langle
\overrightarrow{m}\rangle D_ts_1]\\
&=&\alpha D^{\overrightarrow{j}-\overrightarrow{e_t}}[a\langle
\overrightarrow{m}-\overrightarrow{e_t}\rangle s_1]
+D^{\overrightarrow{j}-\overrightarrow{e_t}}[a\langle
\overrightarrow{m}\rangle D_ts_1].
\end{eqnarray*}
By induction on $\overrightarrow{j}$, we have
$$
D^{\overrightarrow{j}-\overrightarrow{e_t}}[a\langle
\overrightarrow{m}-\overrightarrow{e_t}\rangle s_1]= [a\langle
\overrightarrow{m}-\overrightarrow{e_t}\rangle
D^{\overrightarrow{j}-\overrightarrow{e_t}}s_1]+\varepsilon'_1,
$$
$$
D^{\overrightarrow{j}-\overrightarrow{e_t}}[a\langle
\overrightarrow{m}\rangle D_ts_1]= [a\langle
\overrightarrow{m}\rangle D^{\overrightarrow{j}}s_1]+\varepsilon''_1
$$
where $\varepsilon'_1$ and $\varepsilon''_1$ are linear combination
of normal $S$-words that less than $a\langle \overrightarrow{m}-\overrightarrow{e_t}\rangle
\overline{s_1}D^{\overrightarrow{j}-\overrightarrow{e_t}}$ and
$a\langle \overrightarrow{m}\rangle \overline{s_1}D^{\overrightarrow{j}}$. We
are done since $a\langle \overrightarrow{m}-\overrightarrow{e_t}\rangle
\overline{s_1}D^{\overrightarrow{j}-\overrightarrow{e_t}}<a\langle
\overrightarrow{m}\rangle \overline{s_1}D^{\overrightarrow{j}}$.
\par Assume that $|a|>1$, and write $a=a_1\langle \overrightarrow{m}'\rangle[u_1]$. Then
$$
D^{\overrightarrow{j}}[a\langle \overrightarrow{m}\rangle s_1]=
D^{\overrightarrow{j}}[a_1\langle \overrightarrow{m}'\rangle [u_1\langle
\overrightarrow{m}\rangle s_1]] \ \ \  \ \ \ \ \ \ \ \ \ \ \ \ \ \ \ \ \ \ \ \ \ \ \ \ \
\ \ \ \ \ \ \ \ \ \ \ \ \ \ \ \ \ \ \ \ \ \ \ \ \ \ \  \ \ \
$$
$$
\ \ \ \ \ \ \ \ \  \ \ \ =a_1\langle \overrightarrow{m}'\rangle
D^{\overrightarrow{j}}[u_1\langle \overrightarrow{m}\rangle s_1]
+\sum_{\overrightarrow{l}\in Z_+^n\backslash 0}\alpha_{\overrightarrow{l}}
a_1\langle \overrightarrow{m}'-\overrightarrow{l}\rangle
D^{\overrightarrow{j}-\overrightarrow{l}}[u_1\langle \overrightarrow{m}\rangle s_1].
$$
We are done by induction on $|a|$ and $\overrightarrow{j}$, and
(\ref{e3.16}) follows.
\par Similarly, using induction on $\overrightarrow{j}$ again, one
has
$$
D^{\overrightarrow{j}}[a\langle \overrightarrow{m}\rangle s_2\langle
\overrightarrow{m}'\rangle b]\equiv[a\langle
\overrightarrow{m}\rangle s_2\langle \overrightarrow{m}'\rangle
bD^{\overrightarrow{j}}]+\varepsilon_2
$$
where $\varepsilon_2$ is a linear combination of normal $S$-words
that less than $a\langle \overrightarrow{m}\rangle
\overline{s_2}\langle \overrightarrow{m}'\rangle $
$bD^{\overrightarrow{j}}$, and (\ref{e3.17}) follows.
\par Thus,
\begin{eqnarray*}
D^{\overrightarrow{j}}(s_1,s_2)_w&=&\sum_{i\in
I}\alpha_iD^{\overrightarrow{j}}[u_i\langle
\overrightarrow{m}^{(i)}\rangle s_i\langle
\overrightarrow{m}'^{(i)}\rangle v_i]+\sum_{k\in
J}\alpha_kD^{\overrightarrow{j}}[u_k\langle
\overrightarrow{m}^{(k)}\rangle
D^{\overrightarrow{l_k}}s_k]\\
&\equiv&\sum_{i\in I}\alpha_i[u_i\langle
\overrightarrow{m}^{(i)}\rangle s_i\langle
\overrightarrow{m}'^{(i)}\rangle v_iD^{\overrightarrow{j}}]
+\sum_{k\in J}\alpha_k[u_k\langle \overrightarrow{m}^{(k)}\rangle
D^{\overrightarrow{l_k}+\overrightarrow{j}}s_k] \ \ mod
(S,wD^{\overrightarrow{j}}).
\end{eqnarray*}
Since $u_i\langle \overrightarrow{m}^{(i)}\rangle s_i\langle
\overrightarrow{m}'^{(i)}\rangle v_i<w$ \  and  \ $u_k\langle
\overrightarrow{m}^{(k)}\rangle
\overline{s_k}D^{\overrightarrow{l_k}}<w, $
$
u_i\langle \overrightarrow{m}^{(i)}\rangle s_i\langle \overrightarrow{m}'^{(i)}\rangle
 v_iD^{\overrightarrow{j}}<wD^{\overrightarrow{j}}$ \
and \ $u_k\langle \overrightarrow{m}^{(k)}\rangle
\overline{s_k}D^{\overrightarrow{l_k}+\overrightarrow{j}}<wD^{\overrightarrow{j}}.
$ Therefore, we have
$$
D^{\overrightarrow{j}}(s_1,s_2)_w\equiv 0 \
mod(S,wD^{\overrightarrow{j}}) \ for \ all \ \overrightarrow{j}\in
Z_+^n.
$$
This completes the proof. \ \  $\square$
\begin{theorem}\label{t3.19}
Let $S$ be a Gr\"{o}bner-Shirshov basis in
$C(B,\overrightarrow{N},D_1,\cdots,D_n)$ with fixed
$\overrightarrow{N}$. If $f\in Id(D^{\omega}(S))$, then leading
associative word $\overline{f}$ of $f$ is either
$$
\overline{f}=u\langle \overrightarrow{m}\rangle \overline{s}\langle \overrightarrow{m}'\rangle v \ or
\ \overline{f}=u\langle \overrightarrow{m}\rangle \overline{s}D^{\overrightarrow{j}}
$$
for some $s\in S$ which is D-free in the first case.
\end{theorem}
\noindent{\bf Proof.} Let $f\in Id(D^{\omega}(S))$. By Lemma
\ref{l3.12},
\begin{eqnarray}\label{e3.18}
f=\sum_{i=1}^k\alpha_i[u_i\langle \overrightarrow{m}^{(i)}\rangle
D^{\overrightarrow{j_i}}s_i\langle \overrightarrow{m}'^{(i)}\rangle
v_i]
\end{eqnarray}
 where $[u_i\langle \overrightarrow{m}^{(i)}\rangle
D^{\overrightarrow{j_i}}s_i\langle \overrightarrow{m}'^{(i)}\rangle
v_i]$ are normal $S$-words ($j_i=0$ and $s_i$ is D-free if $v_i\neq
1$ ). We will assume that the right hand side of (\ref{e3.18}) has
no similar summands.

Let us take leading associative words. By Lemma \ref{l3.6},
$$
w_i=u_i\langle \overrightarrow{m}^{(i)}\rangle\overline{s_i} D^{\overrightarrow{j_i}}\langle
\overrightarrow{m}'^{(i)}\rangle v_i, \ \ 1\leqslant i \leqslant k,
$$
and arrange them in decreasing order, i.e.,
$$
w_1=w_2=\cdots=w_l>w_{l+1}\geqslant w_{l+2}\geqslant \cdots
\geqslant w_k.
$$
We will use induction on $(w_1,l)$.

If $l=1$, then $\overline{f}=u_1\langle
\overrightarrow{m}^{(1)}\rangle\overline{s_1}
D^{\overrightarrow{j_1}}\langle \overrightarrow{m}'^{(1)}\rangle
v_1$ and we are done. So let $|l|>1$. Then
$$
w_1=u_1\langle \overrightarrow{m}^{(1)}\rangle\overline{s_1}
D^{\overrightarrow{j_1}}\langle \overrightarrow{m}'^{(1)}\rangle v_1
=u_2\langle \overrightarrow{m}^{(2)}\rangle\overline{s_2}
D^{\overrightarrow{j_2}}\langle \overrightarrow{m}'^{(2)}\rangle
v_2.
$$
Assume first that $v_1\neq 1$ and $v_2\neq 1$. Thus $j_1=0$, $j_2=0$
and $s_1, s_2$ are D-free. Let us rewrite the first two summands of
(\ref{e3.18}) in the form
\begin{eqnarray}\label{e3.20}
&&\nonumber\alpha_1[u_1\langle \overrightarrow{m}^{(1)}\rangle s_1
\langle \overrightarrow{m}'^{(1)}\rangle v_1]+\alpha_2[u_2\langle
\overrightarrow{m}^{(2)}\rangle s_2 \langle
\overrightarrow{m}'^{(2)}\rangle v_2] \\
&=&(\alpha_1+\alpha_2)[u_1\langle \overrightarrow{m}^{(1)}\rangle
s_1 \langle \overrightarrow{m}'^{(1)}\rangle v_1]\\
&&\nonumber +\alpha_2([u_2\langle \overrightarrow{m}^{(2)}\rangle
s_2 \langle \overrightarrow{m}'^{(2)}\rangle v_2]-[u_1\langle
\overrightarrow{m}^{(1)}\rangle s_1 \langle
\overrightarrow{m}'^{(1)}\rangle v_1])
\end{eqnarray}
There are three cases to be discussed.
\par Case 1. $\overline{s_1}$ and $\overline{s_2}$ are mutually
disjoint. We assume that $\overline{s_1}$ is at the left of $s_2$,
i.e.,
$$
u_2=u_1\langle \overrightarrow{m}^{(1)}\rangle \overline{s_1} \langle \overrightarrow{m}'^{(1)}\rangle
a \ and \ v_1=a\langle \overrightarrow{m}^{(2)}\rangle \overline{s_2} \langle
\overrightarrow{m}'^{(2)}\rangle v_2,
$$
here $a\in T$ may be empty. The last item of (\ref{e3.20}) can be
rewritten:
\begin{eqnarray*}
&&[u_2\langle \overrightarrow{m}^{(2)}\rangle s_2 \langle
\overrightarrow{m}'^{(2)}\rangle v_2]-[u_1\langle
\overrightarrow{m}^{(1)}\rangle s_1 \langle
\overrightarrow{m}'^{(1)}\rangle v_1]\\
&=&[u_1\langle \overrightarrow{m}^{(1)}\rangle \overline{s_1}
\langle \overrightarrow{m}'^{(1)}\rangle a \langle
\overrightarrow{m}^{(2)}\rangle s_2 \langle
\overrightarrow{m}'^{(2)}\rangle v_2] -[u_1\langle
\overrightarrow{m}^{(1)}\rangle s_1 \langle
\overrightarrow{m}'^{(1)}\rangle a\langle
\overrightarrow{m}^{(2)}\rangle \overline{s_2} \langle
\overrightarrow{m}'^{(2)}\rangle v_2]\\
&\triangleq&A+B
\end{eqnarray*}
where
\begin{eqnarray*}
&&A=[u_1\langle \overrightarrow{m}^{(1)}\rangle \overline{s_1}
\langle \overrightarrow{m}'^{(1)}\rangle a \langle
\overrightarrow{m}^{(2)}\rangle s_2 \langle
\overrightarrow{m}'^{(2)}\rangle v_2] -[u_1\langle
\overrightarrow{m}^{(1)}\rangle s_1 \langle
\overrightarrow{m}'^{(1)}\rangle a \langle
\overrightarrow{m}^{(2)}\rangle s_2 \langle
\overrightarrow{m}'^{(2)}\rangle v_2],\\
&&B=[u_1\langle \overrightarrow{m}^{(1)}\rangle s_1 \langle
\overrightarrow{m}'^{(1)}\rangle a \langle
\overrightarrow{m}^{(2)}\rangle s_2 \langle
\overrightarrow{m}'^{(2)}\rangle v_2] -[u_1\langle
\overrightarrow{m}^{(1)}\rangle s_1 \langle
\overrightarrow{m}'^{(1)}\rangle a\langle
\overrightarrow{m}^{(2)}\rangle \overline{s_2} \langle
\overrightarrow{m}'^{(2)}\rangle v_2].
\end{eqnarray*}
By Proposition \ref{p3.15}, $A$ and $B$ are linear combination of
normal $S$-words with leading associative words less than $w_1$.
Thus, we can rewrite $f$ with a smaller $(w_1,l)$.
\par Case 2. One of $s_1$ and $s_2$ is a subword of the other, say,
$w=s_1=a\langle \overrightarrow{m}\rangle\overline{s_2}\langle
\overrightarrow{m}'\rangle b$. Then
$\overrightarrow{m}=\overrightarrow{m}^{(2)}$,
$\overrightarrow{m}'=\overrightarrow{m}'^{(2)}$, $u_2=u_1\langle
\overrightarrow{m}^{(1)}\rangle a$ and $v_2=b\langle
\overrightarrow{m}'^{(1)}\rangle v_1$. Let
$(s_1,s_2)_w=s_1-[a\langle
\overrightarrow{m}\rangle\overline{s_2}\langle
\overrightarrow{m}'\rangle b]$, and rewrite the last item of
(\ref{e3.20}) into
\begin{eqnarray*}
&&[u_2\langle \overrightarrow{m}^{(2)}\rangle s_2 \langle
\overrightarrow{m}'^{(2)}\rangle v_2]-[u_1\langle
\overrightarrow{m}^{(1)}\rangle s_1 \langle
\overrightarrow{m}'^{(1)}\rangle v_1]\\
&=&[u_1\langle \overrightarrow{m}^{(1)}\rangle a \langle
\overrightarrow{m}^{(2)}\rangle s_2 \langle
\overrightarrow{m}'^{(2)}\rangle b\langle
\overrightarrow{m}'^{(1)}\rangle v_1] -[u_1\langle
\overrightarrow{m}^{(1)}\rangle s_1 \langle
\overrightarrow{m}'^{(1)}\rangle v_1]\\
&\triangleq&A-B
\end{eqnarray*}
where
\begin{eqnarray*}
A&=&[u_1\langle \overrightarrow{m}^{(1)}\rangle a \langle
\overrightarrow{m}^{(2)}\rangle s_2 \langle
\overrightarrow{m}'^{(2)}\rangle b\langle
\overrightarrow{m}'^{(1)}\rangle v_1] -[u_1\langle
\overrightarrow{m}^{(1)}\rangle [a \langle
\overrightarrow{m}^{(2)}\rangle s_2 \langle
\overrightarrow{m}'^{(2)}\rangle b]\langle
\overrightarrow{m}'^{(1)}\rangle v_1],
\\
B&=&[u_1\langle \overrightarrow{m}^{(1)}\rangle s_1 \langle
\overrightarrow{m}'^{(1)}\rangle v_1]-[u_1\langle
\overrightarrow{m}^{(1)}\rangle [a \langle
\overrightarrow{m}^{(2)}\rangle s_2 \langle
\overrightarrow{m}'^{(2)}\rangle b]\langle
\overrightarrow{m}'^{(1)}\rangle v_1]\\
&=&[u_1\langle \overrightarrow{m}^{(1)}\rangle (s_1,s_2)_w \langle
\overrightarrow{m}'^{(1)}\rangle v_1].
\end{eqnarray*}
We apply Proposition \ref{p3.15} to $A$ and obtain a linear
combination of normal $S$-words with leading associative monomials
less than $w_1$. On the other hand, since $s_1$ and $s_2$ are
D-free, $(s_1,s_2)_w$ is a linear combination of D-free normal
$S$-words which are less than $w$. As a result, we can also apply
Proposition \ref{p3.15} to $B$ and obtain a linear combination of
normal $S$-words which are less than $w_1$. We again rewrite $f$
with a smaller $(w_1,l)$.

Case 3. $s_1$ and $s_2$ have a nonempty intersection as subword of
$w$, and $s_1$ and $s_2$ are not subwords of each other. Assume that
$\overline{s_1}$ is at the left of $\overline{s_2}$, i.e.,
$u_2=u_1\langle \overrightarrow{m}^{(1)}\rangle b, \ v_1=a\langle
\overrightarrow{m}'^{(2)}\rangle v_2$, and $ \overline{s_1}\langle
\overrightarrow{m}^{(1)}\rangle a=b\langle
\overrightarrow{m}^{(2)}\rangle \overline{s_2}=w$, where $a,b\in T$
are D-free, and $|\overline{s_1}|+|\overline{s_2}|>\overline{w}$.
Let $(s_1,s_2)_w=[s_1\langle \overrightarrow{m}^{(1)}\rangle
a]-[b\langle \overrightarrow{m}^{(2)}\rangle s_2]$, and rewrite the
last item of (\ref{e3.18}) into
\begin{eqnarray*}
&&[u_2\langle \overrightarrow{m}^{(2)}\rangle s_2 \langle
\overrightarrow{m}'^{(2)}\rangle v_2]-[u_1\langle
\overrightarrow{m}^{(1)}\rangle s_1 \langle
\overrightarrow{m}'^{(1)}\rangle v_1]\\
&=&[u_1\langle \overrightarrow{m}^{(1)}\rangle b\langle
\overrightarrow{m}^{(2)}\rangle s_2 \langle
\overrightarrow{m}'^{(2)}\rangle v_2] -[u_1\langle
\overrightarrow{m}^{(1)}\rangle s_1 \langle
\overrightarrow{m}'^{(1)}\rangle a\langle
\overrightarrow{m}'^{(2)}\rangle v_2]=A-B-C
\end{eqnarray*}
where
\begin{eqnarray*}
&&A=[u_1\langle \overrightarrow{m}^{(1)}\rangle b\langle
\overrightarrow{m}^{(2)}\rangle s_2 \langle
\overrightarrow{m}'^{(2)}\rangle v_2] -[u_1\langle
\overrightarrow{m}^{(1)}\rangle [b\langle
\overrightarrow{m}^{(2)}\rangle s_2] \langle
\overrightarrow{m}'^{(2)}\rangle v_2],\\
&& B=[u_1\langle \overrightarrow{m}^{(1)}\rangle s_1 \langle
\overrightarrow{m}'^{(1)}\rangle a\langle
\overrightarrow{m}'^{(2)}\rangle v_2] -[u_1\langle
\overrightarrow{m}^{(1)}\rangle [s_1 \langle
\overrightarrow{m}'^{(1)}\rangle a]\langle
\overrightarrow{m}'^{(2)}\rangle v_2],\\
&& C=[u_1\langle \overrightarrow{m}^{(1)}\rangle (s_1,s_2)_w\langle
\overrightarrow{m}'^{(2)}\rangle v_2].
\end{eqnarray*}
By Proposition \ref{p3.15}, $A$ and $B$ are linear combination of
normal $S$-words which are less than $w_1$. Applying Proposition
\ref{p3.15} and using the fact that $S$ is Gr\"{o}bner-Shirshov
basis, we conclude that $C$ is also a linear combination of normal
$S$-words which are less than $w_1$, and we have decrease $(w_1,l)$
in this case as well.
\par Now, we consider the cases of $v_1=v_2=1$ and $v_1\neq
1,v_2=1$.

When $v_1=v_2=1$, we have
$$
w_1=u_1\langle \overrightarrow{m}^{(1)}\rangle\overline{s_1} D^{\overrightarrow{j_1}}=u_2\langle
\overrightarrow{m}^{(2)}\rangle\overline{s_2} D^{\overrightarrow{j_2}}.
$$
Rewrite the first two summands of (\ref{e3.18}) to get
\begin{eqnarray}\label{e3.21}
&&\nonumber \alpha_1[u_1\langle \overrightarrow{m}^{(1)}\rangle
D^{\overrightarrow{j_1}}s_1 ]+\alpha_2[u_2\langle
\overrightarrow{m}^{(2)}\rangle D^{\overrightarrow{j_2}}s_2]\\
&=&(\alpha_1+\alpha_2)[u_1\langle \overrightarrow{m}^{(1)}\rangle
D^{\overrightarrow{j_1}}s_1 ] +\alpha_2([u_2\langle
\overrightarrow{m}^{(2)}\rangle
D^{\overrightarrow{j_2}}s_2]-[u_1\langle
\overrightarrow{m}^{(1)}\rangle D^{\overrightarrow{j_1}}s_1])
\end{eqnarray}
We may assume that
$j_{l_12}\geqslant j_{l_11},\cdots,j_{l_k2}\geqslant j_{l_k1}$,
$j_{l_{k+1}2}\leqslant j_{l_{k+1}1},\cdots,j_{l_n2}\leqslant
j_{l_n1}$, $(l_1,\cdots,l_n)\in S_n$, $0\leqslant k \leqslant n$,
and put $i_1=j_{l_12}-j_{l_11}\geqslant 0,\cdots,i_k=j_{l_k2}-
j_{l_k1}\geqslant 0$, $i_{k+1}=j_{l_{k+1}1}- j_{l_{k+1}2}\geqslant
0,\cdots,i_n=j_{l_n1}- j_{l_n2}\geqslant 0$. Write
$w_1=w'D_{l_1}^{j_{l_11}}\cdots
D_{l_k}^{j_{l_k1}}D_{l_{k+1}}^{j_{l_{k+1}2}}\cdots
D_{l_n}^{j_{l_n2}}$ where
$$
w'=u_1\langle
\overrightarrow{m}^{(1)}\rangle\overline{s_1}D_{l_{k+1}}^{i_{k+1}}\cdots
D_{l_n}^{i_n}=u_2\langle \overrightarrow{m}^{(2)}\rangle
\overline{s_2}D_{l_1}^{i_1}\cdots D_{l_k}^{i_k}
$$
and consider two possibilities:
$$
\overline{s_1}D_{l_{k+1}}^{i_{k+1}}\cdots D_{l_n}^{i_n}=u\langle
\overrightarrow{m}^{(2)}\rangle \overline{s_2}D_{l_1}^{i_1}\cdots D_{l_k}^{i_k}=w , \
\  \overline{s_2}D_{l_1}^{i_1}\cdots D_{l_k}^{i_k}=u\langle
\overrightarrow{m}^{(1)}\rangle\overline{s_1}D_{l_{k+1}}^{i_{k+1}}\cdots
D_{l_n}^{i_n}=w.
$$
In the first case, we have the composition of right including
$(s_1,s_2)_w=D_{l_{k+1}}^{i_{k+1}}\cdots D_{l_n}^{i_n}s_1-[u\langle
\overrightarrow{m}^{(2)}\rangle D_{l_1}^{i_1}\cdots
D_{l_k}^{i_k}s_2]$ which is zero modulo $(S,w)$, and the last item
in (\ref{e3.21}) becomes
\begin{eqnarray*}
&&[u_2\langle \overrightarrow{m}^{(2)}\rangle
D^{\overrightarrow{j_2}}s_2]-[u_1\langle
\overrightarrow{m}^{(1)}\rangle D^{\overrightarrow{j_1}}s_1]\\
&=&[u_1\langle \overrightarrow{m}^{(1)}\rangle u \langle
\overrightarrow{m}^{(2)}\rangle D_{l_1}^{j_{l_11}+i_1}\cdots
D_{l_k}^{j_{l_k1}+i_k}D_{l_{k+1}}^{j_{l_{k+1}2}}\cdots
D_{l_n}^{j_{l_n2}}s_2]\\
&&-[u_1\langle \overrightarrow{m}^{(1)}\rangle
D_{l_1}^{j_{l_11}}\cdots
D_{l_k}^{j_{l_k1}}D_{l_{k+1}}^{j_{l_{k+1}2}+i_{k+1}}\cdots
D_{l_n}^{j_{l_n2}+i_n}s_1]\\
&=&-[u_1\langle \overrightarrow{m}^{(1)}\rangle
(D_{l_1}^{j_{l_11}}\cdots
D_{l_k}^{j_{l_k1}}D_{l_{k+1}}^{j_{l_{k+1}2}}\cdots
D_{l_n}^{j_{l_n2}}(D_{l_{k+1}}^{i_{k+1}}\cdots D_{l_n}^{i_n}s_1)]\\
&&-[u \langle \overrightarrow{m}^{(2)}\rangle
D_{l_1}^{j_{l_11}}\cdots
D_{l_k}^{j_{l_k1}}D_{l_{k+1}}^{j_{l_{k+1}2}}\cdots
D_{l_n}^{j_{l_n2}}(D_{l_1}^{i_1}\cdots D_{l_k}^{i_k}s_2))]\triangleq
A.
\end{eqnarray*}
By Lemma \ref{l3.18}, we have
\begin{eqnarray*}
&&D_{l_1}^{j_{l_11}}\cdots
D_{l_k}^{j_{l_k1}}D_{l_{k+1}}^{j_{l_{k+1}2}}\cdots
D_{l_n}^{j_{l_n2}}(D_{l_{k+1}}^{i_{k+1}}\cdots D_{l_n}^{i_n}s_1)\\
&&\ \ \ \ -[u \langle \overrightarrow{m}^{(2)}\rangle
D_{l_1}^{j_{l_11}}\cdots
D_{l_k}^{j_{l_k1}}D_{l_{k+1}}^{j_{l_{k+1}2}}\cdots
D_{l_n}^{j_{l_n2}}(D_{l_1}^{i_1}\cdots D_{l_k}^{i_k}s_2)]\\
&\equiv &0 \ \ mod(S,wD_{l_1}^{j_{l_11}}\cdots
D_{l_k}^{j_{l_k1}}D_{l_{k+1}}^{j_{l_{k+1}2}}\cdots
D_{l_n}^{j_{l_n2}}).
\end{eqnarray*}
As a result, $A$ is also trivial modulo $(S,w_1)$. Similarly, we are
done in the second possibility. Hence we have $(w_1,l)$ decrease.
\par Finally, suppose that $v_1\neq 1$ and $v_2=1$. We have
$$
w_1=u_1\langle \overrightarrow{m}^{(1)}\rangle\overline{s_1} \langle
\overrightarrow{m}'^{(1)}\rangle v_1 =u_2\langle
\overrightarrow{m}^{(2)}\rangle\overline{s_2}
D^{\overrightarrow{j_2}}
$$
for some $s_1,s_2\in S$, $s_1$ is D-free, and again, the first two
summands of (\ref{e3.18}) becomes
\begin{eqnarray}\label{e3.24}
&&\nonumber\alpha_1[u_1\langle \overrightarrow{m}^{(1)}\rangle s_1
\langle \overrightarrow{m}'^{(1)}\rangle v_1]+\alpha_2[u_2\langle
\overrightarrow{m}^{(2)}\rangle D^{\overrightarrow{j_2}}\\
&=&(\alpha_1+\alpha_2)[u_1\langle \overrightarrow{m}^{(1)}\rangle
s_1 \langle
\overrightarrow{m}'^{(1)}\rangle v_1]\\
&& \nonumber+\alpha_2([u_2\langle \overrightarrow{m}^{(2)}\rangle
D^{\overrightarrow{j_2}}s_2] -[u_1\langle
\overrightarrow{m}^{(1)}\rangle s_1 \langle
\overrightarrow{m}'^{(1)}\rangle v_1])
\end{eqnarray}
Again, we have three cases to consider.
\par Case 1. $\overline{s_1}$ and $\overline{s_2}$ are mutually
disjoint. In this case, write $u_2=u_1\langle
\overrightarrow{m}^{(1)}\rangle \overline{s_1} \langle
\overrightarrow{m}'^{(1)}\rangle a$ and $v_1=a\langle
\overrightarrow{m}^{(2)}\rangle\overline{s_2}
D^{\overrightarrow{j_2}}$, and the last item of (\ref{e3.24})
becomes
\begin{eqnarray*}
&&[u_2\langle \overrightarrow{m}^{(2)}\rangle
D^{\overrightarrow{j_2}}s_2]-[u_1\langle
\overrightarrow{m}^{(1)}\rangle s_1 \langle
\overrightarrow{m}'^{(1)}\rangle v_1]\\
&=&[u_1\langle \overrightarrow{m}^{(1)}\rangle \overline{s_1}
\langle \overrightarrow{m}'^{(1)}\rangle a\langle
\overrightarrow{m}^{(2)}\rangle D^{\overrightarrow{j_2}}s_2]
-[u_1\langle \overrightarrow{m}^{(1)}\rangle s_1 \langle
\overrightarrow{m}'^{(1)}\rangle a\langle
\overrightarrow{m}^{(2)}\rangle\overline{s_2}
D^{\overrightarrow{j_2}}]\triangleq A+B
\end{eqnarray*}
where
\begin{eqnarray*}
&&A=[u_1\langle \overrightarrow{m}^{(1)}\rangle \overline{s_1}
\langle \overrightarrow{m}'^{(1)}\rangle a\langle
\overrightarrow{m}^{(2)}\rangle D^{\overrightarrow{j_2}}s_2]
-[u_1\langle \overrightarrow{m}^{(1)}\rangle s_1 \langle
\overrightarrow{m}'^{(1)}\rangle a\langle
\overrightarrow{m}^{(2)}\rangle D^{\overrightarrow{j_2}}s_2],\\
&&B=[u_1\langle \overrightarrow{m}^{(1)}\rangle s_1 \langle
\overrightarrow{m}'^{(1)}\rangle a\langle
\overrightarrow{m}^{(2)}\rangle D^{\overrightarrow{j_2}}s_2]
-[u_1\langle \overrightarrow{m}^{(1)}\rangle s_1 \langle
\overrightarrow{m}'^{(1)}\rangle a\langle
\overrightarrow{m}^{(2)}\rangle\overline{s_2}
D^{\overrightarrow{j_2}}].
\end{eqnarray*}
By Propositions \ref{p3.15} and \ref{p3.17}, both $A$ and $B$ are
linear combination of normal $S$-words less than $w_1$, hence we
have $(w_1,l)$ decrease.
\par Case 2. $\overline{s_1}$ is subword of $\overline{s_2} D^{\overrightarrow{j_2}}$.
Let $w=\overline{s_2}=a\langle \overrightarrow{m}\rangle
\overline{s_1}\langle \overrightarrow{m}'\rangle b$, where
$\overrightarrow{m}=\overrightarrow{m}^{(1)}$, $
\overrightarrow{m}'=\overrightarrow{m}'^{(1)}$, $b\neq 1$,
$u_1=u_2\langle \overrightarrow{m}^{(2)}\rangle a$, and
$v_1=bD^{\overrightarrow{j_2}}$. A rewriting of the last item of
(\ref{e3.24}) yields
\begin{eqnarray}\label{e3.25}
&&\nonumber [u_2\langle \overrightarrow{m}^{(2)}\rangle
D^{\overrightarrow{j_2}}s_2]-[u_1\langle
\overrightarrow{m}^{(1)}\rangle s_1 \langle
\overrightarrow{m}'^{(1)}\rangle v_1]\\
&=&[u_2\langle
\overrightarrow{m}^{(2)}\rangle(D^{\overrightarrow{j_2}}s_2-
[a\langle \overrightarrow{m}\rangle s_1\langle
\overrightarrow{m}'\rangle bD^{\overrightarrow{j_2}}])]
\end{eqnarray}
Since $(s_1,s_2)_w=s_2-a\langle \overrightarrow{m}\rangle s_1\langle
\overrightarrow{m}'\rangle b$ is a composition of elements of $S$,
it is trivial $mod(S,w)$. By Lemma \ref{l3.18}, the right hand side
of (\ref{e3.25}) is a linear combination of normal $S$-words whose
leading monomials are less than $w_1$, and $(w_1,l)$ is again
decrease.
\par Case 3. $\overline{s_1}$ and $\overline{s_2}$ have a nonempty
intersection and $\overline{s_1}$ is not a subword of
$\overline{s_2} D^{\overrightarrow{j_2}}$. Then $u_2=u_1\langle \overrightarrow{m}^{(1)}\rangle b$,
$\overline{s_1}\langle \overrightarrow{m}'^{(1)}\rangle v_1=b\langle \overrightarrow{m}^{(2)}\rangle
\overline{s_2} D^{\overrightarrow{j_2}}$ and $v_1=aD^{\overrightarrow{j_2}}$, where $a,b\in T$, $a\neq
1$, and so
$$
\overline{s_1}\langle \overrightarrow{m}'^{(1)}\rangle a=b\langle \overrightarrow{m}^{(2)}\rangle
\overline{s_2}=w, \ |\overline{s_1}|+|\overline{s_2}|>|w|.
$$
The composition $(s_1,s_2)_w=[s_1\langle
\overrightarrow{m}'^{(1)}\rangle a]-[b\langle
\overrightarrow{m}^{(2)}\rangle s_2]$ is trivial $mod(S,w)$. By
Lemma \ref{l3.18}, $[s_1\langle \overrightarrow{m}'^{(1)}\rangle
aD^{\overrightarrow{j_2}}]-[b\langle \overrightarrow{m}^{(2)}\rangle
D^{\overrightarrow{j_2}}s_2]$ is also trivial
$mod(S,wD^{\overrightarrow{j_2}})$, i.e., it is a linear combination
of normal $S$-words which are less than $wD^{\overrightarrow{j_2}}$.
Therefore, we rewrite the last item of (\ref{e3.24}) and obtain
\begin{eqnarray*}
&&[u_2\langle \overrightarrow{m}^{(2)}\rangle
D^{\overrightarrow{j_2}}s_2]-[u_1\langle
\overrightarrow{m}^{(1)}\rangle s_1 \langle
\overrightarrow{m}'^{(1)}\rangle v_1]\\
&=&[u_1\langle \overrightarrow{m}^{(1)}\rangle b\langle
\overrightarrow{m}^{(2)}\rangle D^{\overrightarrow{j_2}}s_2]
-[u_1\langle \overrightarrow{m}^{(1)}\rangle s_1 \langle
\overrightarrow{m}'^{(1)}\rangle aD^{\overrightarrow{j_2}}]\\
&=&-[u_1\langle \overrightarrow{m}^{(1)}\rangle([s_1 \langle
\overrightarrow{m}'^{(1)}\rangle aD^{\overrightarrow{j_2}}]-
[b\langle \overrightarrow{m}^{(2)}\rangle
D^{\overrightarrow{j_2}}s_2])].
\end{eqnarray*}
The last expression is a linear combination of normal $S$-words that
less than $w_1$. Again, we have $(w_1,l)$ decrease, and the proof is
complete.  \ \  $\square$

\begin{definition}
Let $S$ be a Gr\"{o}bner-Shirshov basis in
$C(B,\overrightarrow{N},D_1,\cdots,D_n)$. A normal word $[u]$ is
said to be $S$-irreducible if $u$ is not the form $a\langle
\overrightarrow{m}\rangle \overline{s}\langle
\overrightarrow{m}'\rangle b$, $s\in S$, $s$ D-free, $a,b\in T$, nor
of the form $a\langle \overrightarrow{m}\rangle
\overline{s}D^{\overrightarrow{i}}$, $s\in S,\ a\in T$.
\end{definition}
\begin{theorem} (Composition-Diamond lemma)
Let $S$ be a Gr\"{o}bner-Shirshov basis in
$C(B,\overrightarrow{N},D_1,\cdots,D_n)$ with fixed
$\overrightarrow{N}$. Then the normal $S$-irreducible words form a
linear basis of the $n$-conformal algebra
$C(B,\overrightarrow{N},D_1,\cdots,D_n|S)$ with defining relations
$S$.
\par If $S$ is D-free, then the converse is true as well.
\end{theorem}
\noindent{\bf Proof.} Let us define the algorithm of elimination of
leading words (ELW) of $S$ in normal words. Let $[u]$ be a normal
word, and let
$$
u=a\langle \overrightarrow{m}\rangle \overline{s}\langle
\overrightarrow{m}'\rangle b,\ s \ \mbox{ is }\ D-free \  \mbox{ or
} \ u=a\langle \overrightarrow{m}\rangle
\overline{s}D^{\overrightarrow{i}}
$$
where $s\in S$, $a,b\in T$ with $a$ being D-free. The
transformations
$$
[u]\mapsto[u]-[a\langle \overrightarrow{m}\rangle s\langle
\overrightarrow{m}'\rangle b] \  \mbox{ and }\
[u]\mapsto[u]-[a\langle \overrightarrow{m}\rangle
D^{\overrightarrow{i}}s]
$$
will be called the results of ELW of $s$ in $[u]$. From Lemma
\ref{l3.6}, we have
$$
\overline{[a\langle \overrightarrow{m}\rangle s\langle
\overrightarrow{m}'\rangle b]}=u  \  \mbox{ and }\
\overline{[a\langle \overrightarrow{m}\rangle
D^{\overrightarrow{i}}s]}=u.
$$
Hence, any normal word is a linear combination of $S$-irreducible
normal words modulo $S$. It follows from Theorem \ref{t3.19} that
$S$-irreducible normal words are linear independent. This completes
the proof of the first part of the Composition-Diamond lemma.

To prove the second part, we assume that $S$ is D-free, and that the
$S$-irreducible normal words form a linear basis of
$C(B,\overrightarrow{N},D_1,\cdots,D_n|S)$. Suppose that $h$ is
nontrivial composition of elements of $S$. So, $h=(f,g)_w$ or
$h=a\langle \overrightarrow{m}\rangle f$ for some $f,g\in S$. We can
apply to $h$ the process of ELW's of $S$, and in finite number of
steps, we will have the following presentation of $h$ in
$C(B,\overrightarrow{N},D_1,\cdots,D_n)$:
\begin{eqnarray}\label{e3.26}
h=\sum_{i=1}^k\alpha_i[u_i\langle \overrightarrow{m}^{(i)}\rangle
s_i\langle \overrightarrow{m}'^{(i)}\rangle
v_i]+h_1
\end{eqnarray}
where $h_1$ is a linear combination of $S$-irreducible normal D-free
words, $[u_i\langle \overrightarrow{m}^{(i)}\rangle s_i\langle
\overrightarrow{m}'^{(i)}\rangle v_i]$, $1\leqslant i \leqslant k$,
are D-free normal $S$-words. Moreover, we have
$$
u_i\langle \overrightarrow{m}^{(i)}\rangle \overline{s_i} \langle
\overrightarrow{m}'^{(i)}\rangle v_i\leqslant \overline{h} \ \mbox{
for all } \ i.
$$
Since $h$ is nontrivial mod $S$, $h_1\neq 0$ in (\ref{e3.26}). On
the other hand, $h_1\in Id(D^{\omega}(S))$, and by Theorem
\ref{t3.19}, $h_1$ must contain a subword $\overline{s}$ for some
(D-free) $s\in S$, which contradicts the fact that $h_1$ is
$S$-irreducible. This completes the proof of the Composition-Diamond
lemma. \ \  $\square$

\ \

\noindent{\bf Remark.} The condition of triviality of compositions
modulo $S$ is much weaker and much easier to apply than checking
that $(f,g)_w$ $(a\langle \overrightarrow{m}\rangle f,\ g\langle
\overrightarrow{m}\rangle a)$ goes to zero using the ELW's of $S$.
Yet, the result for D-free set $S$ will be the same: If any
composition of elements of $S$ is trivial in the sense of the
previous definition, then it is trivial in the sense of ELW's of
$S$. However if the relations are not D-free, the result is not the
same.

\par \noindent{\bf Remark.} Let $S\subseteq C(B,\overrightarrow{N},D_1,\cdots,D_n)$ be a Gr\"{o}bner-Shirshov
basis. We call $S$ an irreducible Gr\"{o}bner-Shirshov basis if for
any $s\in S$, $s$ is a linear combination of $(S\backslash
\{s\})$-irreducible normal words. This is the same as to say that
for any normal word $[w]$ of $S$,
\begin{eqnarray*}
&&w\neq u\langle \overrightarrow{m}\rangle \overline{s_i}\langle
\overrightarrow{m}'\rangle v, \ s_i \ \mbox{ is } \ D\mbox{-free},\\
&&w\neq u\langle \overrightarrow{m}\rangle
\overline{s_i}D^{\overrightarrow{j} \ \ \ \mbox{ for any } s_i\in
S\backslash \{s\}}.
\end{eqnarray*}

It follows from the Composition-Diamond lemma that
\begin{theorem}
Any $n$-conformal ideal of $C(B,\overrightarrow{N},D_1,\cdots,D_n)$
generated by a D-free set has a unique irreducible
Gr\"{o}bner-Shirshov basis.
\end{theorem}

Also the algorithm for applying defining relations, i.e., the ELW's
algorithm, gives rise the following
\begin{theorem} The word problem for any $n$-conformal algebra with
finite or recursive Gr\"{o}bner-Shirshov basis is algorithmically
solvable.
\end{theorem}

\section{Application}
\subsection{Universal enveloping  $n$-conformal algebras }

Let $L$ be a Lie $n$-conformal algebra. By this we mean that $L$ is
a linear space over $k$ equipped with multiplications $\lfloor
\overrightarrow{m} \rfloor, \overrightarrow{m} \in Z_+^n$, and
linear maps $D_i, 1\leqslant i \leqslant n$ such that $(L, \lfloor
\overrightarrow{m} \rfloor, \overrightarrow{m} \in Z_+^n,
D_1,\cdots,D_n)$ is a $n$-conformal algebra with two additional
axioms:
\par $\bullet$ (Anti-commutativity) $a\lfloor \overrightarrow{m} \rfloor b=
\{b\lfloor \overrightarrow{m} \rfloor a \}$, where
$$
\{b\lfloor \overrightarrow{m} \rfloor a
\}=\sum_{\overrightarrow{s}\in Z_+^n}
(-1)^{\overrightarrow{m}+\overrightarrow{s}}D^{(\overrightarrow{s})}
(b\lfloor \overrightarrow{m}+\overrightarrow{s} \rfloor a), \
D^{(\overrightarrow{s})}
=D^{\overrightarrow{s}}/\overrightarrow{s}!, \ \overrightarrow{s}!=s_1!\cdots s_n!.
$$
\par $\bullet$ (Jacobi identity)
$$
(a\lfloor \overrightarrow{m} \rfloor b) \lfloor \overrightarrow{m}' \rfloor c=
\sum_{\overrightarrow{s}\in Z_+^n}(-1)^{\overrightarrow{s}}
\binom {\overrightarrow{m}}{\overrightarrow{s}}
(a\lfloor \overrightarrow{m}-\overrightarrow{s} \rfloor
(b\lfloor \overrightarrow{m}'+\overrightarrow{s} \rfloor c)-
b\lfloor \overrightarrow{m}+\overrightarrow{s} \rfloor
(a\lfloor \overrightarrow{m}'-\overrightarrow{s} \rfloor c)).
$$
\par Any associative $n$-conformal algebra $A=(A,\langle \overrightarrow{m}\rangle,
\overrightarrow{m} \in Z_+^n,D_1,\cdots,D_n)$ can be made into a Lie
$n$-conformal algebra $A^{(-)}=(A,\lfloor \overrightarrow{m}
\rfloor, \overrightarrow{m} \in Z_+^n,D_1,\cdots,D_n)$ by defining
new multiplications using $n$-conformal commutators:
$$
a\lfloor \overrightarrow{m} \rfloor b=a\langle
\overrightarrow{m}\rangle b- \{b\langle \overrightarrow{m}\rangle a
\}, \ \overrightarrow{m} \in Z_+^n, \ a,b\in A.
$$
The locality function for $A^{(-)}$ is essentially the same as for $A$. Namely, it is
given by
$$
\overrightarrow{N}_{A^{(-)}}(a,b)=max\{\overrightarrow{N}_{A}(a,b),\overrightarrow{N}_{A}(b,a)\}.
$$
\par Let $L$ be a Lie $n$-conformal algebra which is a free
$k[D_1,\cdots,D_n]$-module with a basis $B=\{a_i|i\in I\}$ and a
bounded locality $\overrightarrow{N}(a_i,a_j)\preceq
\overrightarrow{N}$ for all $i,j\in I$. Let the multiplication table
of $L$ in the basis $B$ be
$$
a_i\lfloor \overrightarrow{m} \rfloor a_j=\Sigma\alpha_{ij}^k a_k, \
\alpha_{ij}^k\in k[D_1,\cdots,D_n], i\geqslant j, \ i,j\in I, \
\overrightarrow{m}\prec  \overrightarrow{N}.
$$
Then by $\mathcal{U}_{\overrightarrow{N}}(L)$, a universal
enveloping $n$-conformal algebra of $L$ with respective to $B$ and
$\overrightarrow{N}$, one means the following associative
$n$-conformal algebra:
$$
\mathcal{U}_{\overrightarrow{N}}(L)=C(B,\overrightarrow{N}| \
a_i\langle \overrightarrow{m}\rangle a_j -\{a_j\langle
\overrightarrow{m}\rangle a_i \}-a_i\lfloor \overrightarrow{m}
\rfloor a_j=0, \ i\geqslant j, \ i,j\in I, \ \overrightarrow{m}\prec
\overrightarrow{N}).
$$
\begin{theorem}(1/2-PBW Theorem for Lie $n$-conformal algebras) Let
$$
L=(\{a_i|i\in I\},\overrightarrow{N}| \ a_i\lfloor
\overrightarrow{m} \rfloor a_j=\Sigma\alpha_{ij}^k a_k, \ i\geqslant
j, \ i,j\in I, \ \overrightarrow{m}\prec \overrightarrow{N})
$$
be a Lie $n$-conformal algebra with the basis $\{a_i|i\in I\}$ over
$k[D_1,\cdots,D_n]$ and locality $\overrightarrow{N}$. Let
\begin{eqnarray*}
\mathcal{U}_{\overrightarrow{N}}(L)&=&C(\{a_i|i\in
I\},\overrightarrow{N}|\ s_{ij}^{\overrightarrow{m}}=a_i\langle
\overrightarrow{m}\rangle a_j -\{a_j\langle
\overrightarrow{m}\rangle a_i \}-a_i\lfloor \overrightarrow{m}
\rfloor a_j=0,\\
&& \ \ \ \ \ \ \ i\geqslant j, \ i,j\in I, \ \overrightarrow{m}\prec
\overrightarrow{N})
\end{eqnarray*}
be a universal enveloping $n$-conformal algebra of $L$. Let
$S=\{s_{ij}^{\overrightarrow{m}}| i,j\in I, \
\overrightarrow{0}\preceq\overrightarrow{m}\prec
\overrightarrow{N}\}$. Then any polynomial
$$
s_{ij}^{\overrightarrow{m}}\langle \overrightarrow{m}'\rangle a_k-
a_i\langle \overrightarrow{m}\rangle s_{jk}^{\overrightarrow{m}'}
$$
where $i>j>k, \overrightarrow{m} \prec \overrightarrow{N},
\overrightarrow{m}' \prec \overrightarrow{N}$, is trivial $mod(S,w)$
for $w=a_i\langle\overrightarrow{m}\rangle
a_j\langle\overrightarrow{m}'\rangle a_k$.
\end{theorem}

\subsection{An example about the Gr\"{o}bner-Shirshov basis of an associative $n$-conformal algebra}

\par In what follows, a word $a_1\langle \overrightarrow{m}^{(1)}\rangle a_2\langle
\overrightarrow{m}^{(2)}\rangle\cdots a_k\langle
\overrightarrow{m}^{(k)}\rangle D^{\overrightarrow{i}}a_{k+1}$ will
stand for the right normed word $[a_1\langle
\overrightarrow{m}^{(1)}\rangle a_2\langle
\overrightarrow{m}^{(2)}\rangle\cdots a_k\langle
\overrightarrow{m}^{(k)}\rangle D^{\overrightarrow{i}}a_{k+1}]$,
i.e., we will omit the right normed brackets.
\par Consider the associative $n$-conformal algebra
$$
G=C({a},(2,2),D_1,D_2|a\langle 0,0\rangle a=a).
$$
Let $f=a\langle 0,0\rangle a-a$. By the Gr\"{o}bner-Shirshov basis
of $G$, we understand the Gr\"{o}bner-Shirshov basis of
$I=Id(D^{\omega}(f)\subseteq C({a},(2,2),D_1,D_2))$.

\begin{theorem}\label{t4.2}
The Gr\"{o}bner-Shirshov basis of $G$ is as follows:
$$
f=a\langle 0,0\rangle a-a, \ \  s=a\langle 1,1\rangle a\langle
1,1\rangle a,
$$
$$
g=a\langle 1,0\rangle a\langle 1,0\rangle a, \ \  p=a\langle
1,1\rangle a\langle 1,0\rangle a,
$$
$$
h=a\langle 0,1\rangle a\langle 0,1\rangle a, \ \ q=a\langle
1,1\rangle a\langle 0,1\rangle a.
$$
\end{theorem}

The following lemma will be used in proving the above theorem.

\begin{lemma}\label{l4.3}
Let $C(B,\overrightarrow{N},D_1,\cdots,D_n)$ be a free associative
$n$-conformal algebra with locality
$\overrightarrow{N}(b,c)=\overrightarrow{N}$ for all $b,c\in B$. Let
$u$ be a D-free word and $a\in B$.

(i) If $\overrightarrow{N}\not\prec  \overrightarrow{m}$, then
$(u)\langle \overrightarrow{m}\rangle a=0$. (Here $(u)$ is any
bracketing of $u$.)

(ii) $[u]\langle \overrightarrow{N}-\overrightarrow{1}\rangle a
=[u\langle \overrightarrow{N}-\overrightarrow{1}\rangle a].$
\end{lemma}
\noindent{\bf Proof.} (i) We will carry out the proof by induction
on the length of $u$. If $|u|=1$, then $u\in B$, and the result is
true. Assume $|u|>1$. Hence $(u)=(v)\langle
\overrightarrow{k}\rangle(w)$ for some words $v$ and $w$ with
$|u|>v$ and $|u|>w$. Thus, by induction,
$$
(u)\langle \overrightarrow{m}\rangle a=
(v)\langle \overrightarrow{k}\rangle(w)\langle \overrightarrow{m}\rangle a=
\sum_{\overrightarrow{s}\in
Z_+^n}(-1)^{\overrightarrow{s}}\binom{\overrightarrow{k}}{\overrightarrow{s}}(v)\langle
 \overrightarrow{k}-\overrightarrow{s}\rangle ((w)\langle \overrightarrow{m}+\overrightarrow{s} \rangle
 a)=0.
$$

(ii) Again, we prove by induction on the length of $u$. The result
holds trivially for $|u|=1$. Assume that $|u|>1$. Then
$[u]=a_1\langle \overrightarrow{m}\rangle [v]$ where $a_1\in B,
\overrightarrow{m}\in Z_+^n, \ |v|<|u|$ and $v$ is D-free. By (i),
$a_1\langle\overrightarrow{m}-\overrightarrow{s}\rangle
[v]\langle\overrightarrow{m}-\overrightarrow{1}+\overrightarrow{s}\rangle
a=0$ for all $\overrightarrow{s}\in Z_+^n\backslash0$, and so by
induction we have
\begin{eqnarray*}
&&[u]\langle \overrightarrow{N}-\overrightarrow{1}\rangle
a=(a_1\langle \overrightarrow{m}\rangle [v]) \langle
\overrightarrow{N}-\overrightarrow{1}\rangle a \\
&=&a_1\langle \overrightarrow{m}\rangle [v]\langle
\overrightarrow{N}-\overrightarrow{1}\rangle a
+\sum_{\overrightarrow{s}\in
Z_+^n\backslash0}(-1)^{\overrightarrow{s}}\binom{\overrightarrow{m}}{\overrightarrow{s}}
a_1\langle \overrightarrow{m}-\overrightarrow{s}\rangle ([v]\langle
\overrightarrow{N}-\overrightarrow{1}+\overrightarrow{s} \rangle a\\
&=&a_1\langle \overrightarrow{m}\rangle [v \langle
\overrightarrow{N}-\overrightarrow{1}\rangle a] =a_1\langle
\overrightarrow{m}\rangle v \langle
\overrightarrow{N}-\overrightarrow{1}\rangle a =[u\langle
\overrightarrow{N}-\overrightarrow{1}\rangle a].
\end{eqnarray*}
This completes the proof.   \ \  $\square$

 \noindent{\bf Proof of the Theorem \ref{t4.2}.}
We shall check all possible compositions and prove that all of them
are trivial. Explicitly, the proof will be carried out in the
following order:
\begin{enumerate}
\item[1)]\
Show that $g,h\in I$. Consequently, $a\langle 2,0\rangle f\equiv 0$
mod($g$) and $a\langle 0,2\rangle f\equiv 0$ mod($h$).

\item[2)]\ Show that $p,q\in I$. Consequently, $a\langle 2,1\rangle f\equiv
0$ mod($p$) and $a\langle 1,2\rangle f\equiv 0$ mod($q$).

\item[3)]\ Show that $s\in I$. Consequently, $a\langle 2,2\rangle f\equiv 0$
mod($s$).

\item[4)]\ Show that $a\langle m_1,m_2\rangle f=0$ for $m_1\geqslant 3$ or
$m_2\geqslant 3$.

\item[5)]\ Show that $a\langle m_1,m_2\rangle g \equiv 0$ mod$(q,s)$,
$a\langle m_1,m_2\rangle h\equiv0$ mod$(p,s)$, $a\langle
m_1,m_2\rangle p\equiv0$ mod$(s)$, $a\langle m_1,m_2\rangle
q\equiv0$ mod$(s)$, $a\langle m_1,m_2\rangle s=0$ for
$(m_1,m_2)\not\prec  (2,2)$.

\item[6)]\ Show that $(f,f)_{w}=0$, $(g,f)_w \equiv 0$ mod$(f,g;w)$,
$(h,f)_w \equiv 0$ mod$(f,h;w)$, $(p,f)_w \equiv 0$ mod$(f,p;w)$,
$(q,f)_w \equiv 0$ mod$(f,q;w)$, $(s,f)_w \equiv 0$ mod$(f,s;w)$.

\item[7)]\ Show that $(f,g)_{w}\equiv 0$ mod$(g;w)$, $(g,g)_w=0$,
$(h,g)_w \equiv 0$ mod$(q;w)$, $(p,g)_w=0$, $(q,g)_w=0$, $(s,g)_w
=0$.

\item[8)]\ Show that $(f,h)_{w}\equiv 0$ mod$(h;w)$, $(g,h)_w \equiv 0$ mod$(q;w)$,
$(h,h)_w =0$, $(p,h)_w \equiv 0$ mod$(q;w)$, $(q,h)_w \equiv 0$
mod$(q;w)$, $(s,h)_w \equiv 0$ mod$(q;w)$.

\item[9)]\ Show that $(f,p)_{w}\equiv 0$ mod$(p;w)$, $(g,p)_w =0$,
$(h,p)_w \equiv 0$ mod$(s;w)$, $(p,p)_w=0$, $(q,p)_w=0$, $(s,p)_w
=0$.

\item[10)]\ Show that $(f,q)_{w}\equiv 0$ mod$(q;w)$, $(g,q)_w \equiv 0$ mod$(s;w)$,
$(h,q)_w=0$, $(p,q)_w \equiv 0$ mod$(s;w)$, $(q,q)_w \equiv 0$
mod$(s;w)$, $(s,q)_w \equiv 0$ mod$(s;w)$.

\item[11)]\ Show that $(f,s)_{w}\equiv 0$ mod$(s;w)$, $(g,s)_w=0$,
$(h,s)_w=0$, $(p,s)_w=0$, $(q,s)_w=0$, $(s,s)_w=0$.
\end{enumerate}

\ \

Now we check 1)-11) step by step.

1) $g,h\in I$. To see this, we consider
\begin{eqnarray*}
a\langle 2,0\rangle f&=&a\langle 2,0\rangle a \langle 0,0\rangle
a-a\langle 2,0\rangle a\\
& =&(a\langle 2,0\rangle a)\langle 0,0\rangle a+ 2a\langle
1,0\rangle a \langle 1,0\rangle a-a\langle 0,0\rangle a \langle
2,0\rangle a\\
& =&2a\langle 1,0\rangle a \langle 1,0\rangle a.
\end{eqnarray*}
Hence $a\langle 2,0\rangle f=2g$, and so $g_0\in I$. Moreover, this
shows that $a\langle 2,0\rangle f\equiv 0$  mod$(g)$.

Analogously, one can show that $h\in I$, $a\langle 0,2\rangle
f\equiv 0$ mod$(h)$.

2) $p\in I$. To see this, we consider
\begin{eqnarray*}
a\langle 2,1\rangle f&=&a\langle 2,1\rangle a \langle 0,0\rangle
a-a\langle 2,1\rangle a\\
&=&(a\langle 2,1\rangle a)\langle 0,0\rangle a+ 2a\langle 1,1\rangle
a \langle 1,0\rangle a-a\langle 0,1\rangle a \langle 2,0\rangle a \\
&&+\ a\langle 2,0\rangle a \langle 0,1\rangle a- 2a\langle
1,0\rangle
a \langle 1,1\rangle a +a\langle 0,0\rangle a \langle 2,1\rangle a\\
&=&2a\langle 1,1\rangle a \langle 1,0\rangle a+a\langle 2,0\rangle a
\langle 0,1\rangle a- 2a\langle 1,0\rangle a \langle 1,1\rangle a.
\end{eqnarray*}
Now, it is easy to check that $ a\langle 2,0\rangle a \langle
0,1\rangle a=2a\langle 1,0\rangle a \langle 1,1\rangle a. $ Hence
$a\langle 2,1\rangle f=2a\langle 1,1\rangle a \langle 1,0\rangle
a=2p$, and so $p\in I$. Moreover, this shows that $a\langle
2,1\rangle f$ is trivial modulo $p$.
\par Analogously, one can show that $q\in I$, $a\langle 1,2\rangle f\equiv 0$  mod$(q)$.

3) $s\in I$. To see this, we consider
\begin{eqnarray*}
a\langle 2,2\rangle f&=&a\langle 2,2\rangle a \langle 0,0\rangle
a-a\langle 2,2\rangle a\\
&=&(a\langle 2,2\rangle a)\langle 0,0\rangle a+ 2a\langle 1,2\rangle
a \langle 1,0\rangle a-a\langle 0,2\rangle a \langle 2,0\rangle a \\
&&+2a\langle 2,1\rangle a \langle 0,1\rangle a- 4a\langle 1,1\rangle
a \langle 1,1\rangle a +2a\langle 0,1\rangle a \langle 2,1\rangle
a-a\langle 0,0\rangle a \langle 2,2\rangle a\\
&=&2a\langle 1,2\rangle a \langle 1,0\rangle a+2a\langle 2,1\rangle
a \langle 0,1\rangle a- 4a\langle 1,1\rangle a \langle 1,1\rangle a.
\end{eqnarray*}
Now, it is easy to check that
$$
a\langle 1,2\rangle a \langle 1,0\rangle a=2a\langle 1,1\rangle a
\langle 1,1\rangle a, \ \  a\langle 2,1\rangle a \langle
0,1\rangle a2a\langle 1,1\rangle a \langle 1,1\rangle a.
$$
Hence $a\langle 2,2\rangle f=4a\langle 1,1\rangle a \langle
1,1\rangle a=4s$, and so $s\in I$. Moreover, this shows that
$a\langle 2,2\rangle f$ is trivial modulo $s$.

4) $a\langle m_1,m_2\rangle f=0$ for $m_1\geqslant 3$ or
$m_2\geqslant 3$. The proof will follow from the following lemma.
\begin{lemma}\label{l4.4}
$a\langle m_1,m_2\rangle a \langle n_1,n_2\rangle a=0$ for all
$m_1+n_1 \geqslant 3$ or $m_2+n_2 \geqslant 3$.
\end{lemma}
\noindent{\bf Proof.} Assume that $m_1+n_1 \geqslant 3$. For
$(m_1,m_2)\prec  (2,2)$, the result is clear for $n_1\geqslant 2$.
For $(m_1,m_2)\not\prec  2$, we expand $a\langle m_1,m_2\rangle a
\langle n_1,n_2\rangle a$ and note that it is a linear combination
of
$$
v_i=a\langle k_{i1},k_{i2}\rangle a \langle l_{i1},l_{i2}\rangle a
$$
where $k_{i1}+l_{i1}=m_1+n_1 \geqslant 3$, $k_{i2}+l_{i2}=m_2+n_2$,
$k_{i1},k_{i2}\prec  (2,2)$, $i\in I$. So $v_i=0$ by the above
proof. This shows that $a\langle m_1,m_2\rangle a \langle
n_1,n_2\rangle a=0$. Analogously, one can show the case of $m_2+n_2
\geqslant 3$ and the proof is done. \ \ $\square$
\par The proof of the next lemma is similar to that of Lemma
\ref{l4.4}
\begin{lemma}\label{l4.5}
$a\langle m_1,m_2\rangle a \langle n_1,n_2\rangle a \langle
t_1,t_2\rangle a =0$ for all $m_1+n_1+t_1 \geqslant 4$ or
$m_2+n_2+t_2 \geqslant 4$.
\end{lemma}

5) $a\langle m_1,m_2\rangle g\equiv 0 $ mod$(p,s)$ for
$(m_1,m_2)\not\prec  (2,2)$. First of all, we show that $a\langle
0,2\rangle g\equiv 0$ mod$(p,s)$. Expanding $a\langle 0,2\rangle g$,
we get
\begin{eqnarray*}
a\langle 0,2\rangle g&=&a\langle 0,2\rangle a\langle 1,0\rangle
a\langle 1,0\rangle a\\
&=&(a\langle 0,2\rangle a)\langle 1,0\rangle a\langle 1,0\rangle
a+2a\langle 0,1\rangle a\langle 1,1\rangle a\langle 1,0\rangle
a-a\langle 0,0\rangle a\langle 1,2\rangle a\langle 1,0\rangle a\\
&=&2a\langle 0,1\rangle a\langle 1,1\rangle a\langle 1,0\rangle
a-2a\langle 0,0\rangle a\langle 1,1\rangle a\langle 1,1\rangle
a\\
&\equiv &0 \ \ mod(p,s).
\end{eqnarray*}
Similarly, we can prove that
$$
a\langle 1,2\rangle g \equiv 0 \ mod(p,s), \ a\langle 0,3\rangle g
\equiv 0 \ mod(s), \ a\langle 1,3\rangle g \equiv 0 \ mod(s).
$$
Now, $a\langle m_1,m_2\rangle g=0$ for $m_1\geqslant 2$ or
$m_2\geqslant 4$ by Lemma \ref{l4.5}. Hence $a\langle m_1,m_2\rangle
g\equiv 0$ mod$(p,s)$ for $(m_1,m_2)\not\prec  (2,2)$.
\par Analogously, one can show that if $(m_1,m_2)\not\prec
(2,2)$ then
\begin{eqnarray*}
&&a\langle m_1,m_2\rangle h \equiv0 \ mod(q,s),\ \ \ a\langle
m_1,m_2\rangle p \equiv 0 \  mod(s),\\
&&a\langle m_1,m_2\rangle q \equiv0 \ mod(s),\ \ \ a\langle
m_1,m_2\rangle s=0.
\end{eqnarray*}

6) $(f,f)_w=0$. Take
$$
w=a\langle 0,0\rangle a\langle 0,0\rangle a.
$$
Expanding $(f,f)_w$, we get
$$
(f,f)_w=f\langle 0,0\rangle a-a\langle 0,0\rangle f=(a\langle
0,0\rangle a)\langle 0,0\rangle a-a\langle 0,0\rangle a-[w]+a\langle
0,0\rangle a=0.
$$
Next, we show that $(g,f)_w\equiv0$. Take
$$
w=a\langle 1,0\rangle a\langle 1,0\rangle a\langle 0,0\rangle a.
$$
Expanding $(g,f)_w$, we get
$$
(g,f)_w=g\langle 0,0\rangle a-a\langle 1,0\rangle a\langle
1,0\rangle f=(a\langle 1,0\rangle a\langle 1,0\rangle a)\langle
0,0\rangle a-[w]+g.
$$
Expanding the first monomial on the right hand side, we have
\begin{eqnarray*}
&&(a\langle 1,0\rangle a\langle 1,0\rangle a)\langle 0,0\rangle a \\
&=&a\langle 1,0\rangle (a\langle 1,0\rangle a)\langle 0,0\rangle
a-a\langle 0,0\rangle (a\langle 1,0\rangle a)\langle 1,0\rangle a \\
&=&[w]-a\langle 1,0\rangle a\langle 0,0\rangle a\langle 1,0\rangle
-a\langle 0,0\rangle a\langle 1,0\rangle a\langle 1,0\rangle
a+a\langle 0,0\rangle a\langle 0,0\rangle a\langle 2,0\rangle a.
\end{eqnarray*}
Now $a\langle 0,0\rangle a\langle 0,0\rangle a\langle 2,0\rangle
a=0$. Therefore,
$$
(g,f)_w=-a\langle 1,0\rangle a\langle 0,0\rangle a\langle 1,0\rangle
-a\langle 0,0\rangle a\langle 1,0\rangle a\langle 1,0\rangle
a+g\equiv 0 \ mod \ (f,g;w).
$$
\par Analogously, one can show that
$$
(h,f)_w \equiv 0 \ mod \ (f,h;w), \ \  (p,f)_w \equiv 0  \ mod
(f,p;w),
$$
$$
(q,f)_w \equiv 0 \ mod \ (f,q;w), \ \ (s,f)_w \equiv 0 \  mod  \
(f,s;w).
$$

7) $(f,g)_w\equiv0$. Take
$$
w=a\langle 0,0\rangle a\langle 1,0\rangle a\langle 1,0\rangle a.
$$
Expanding $(f,g)_w$, we get
$$
(f,g)_w=f\langle 1,0\rangle a\langle 1,0\rangle a-a\langle
0,0\rangle g=(a\langle 0,0\rangle a)\langle 1,0\rangle a\langle
1,0\rangle a-g-[w]=-g\equiv 0 \ mod \ (g).
$$
Next, we show that $(g,g)_w=0$. There are two compositions
$$
w_1=a\langle 1,0\rangle a\langle 1,0\rangle a\langle 1,0\rangle
a\langle 1,0\rangle a, \ \ w_2=a\langle 1,0\rangle a\langle
1,0\rangle a\langle 1,0\rangle a.
$$
We consider $w_1$ at first. Expanding $(g,g)_{w_1}$, we get
$$
(g,g)_{w_1}=g\langle 1,0\rangle a\langle 1,0\rangle a-a\langle
1,0\rangle a\langle 1,0\rangle g =(a\langle 1,0\rangle a\langle
1,0\rangle a)\langle 1,0\rangle a\langle 1,0\rangle a-[w_1].
$$
Expanding the first monomial on the right hand side, we have
\begin{eqnarray*}
&&(a\langle 1,0\rangle a\langle 1,0\rangle a)\langle 1,0\rangle
a\langle 1,0\rangle a \\
&=&a\langle 1,0\rangle (a\langle 1,0\rangle a)\langle 1,0\rangle
a\langle 1,0\rangle a-a\langle 0,0\rangle (a\langle 1,0\rangle
a)\langle 2,0\rangle a\langle 1,0\rangle a \\
&=&[w]-a\langle 1,0\rangle a\langle 0,0\rangle a\langle 2,0\rangle
a\langle 1,0\rangle a \\
&&
-a\langle 0,0\rangle a\langle 1,0\rangle a\langle 2,0\rangle
a\langle 1,0\rangle a+a\langle 0,0\rangle
a\langle 0,0\rangle a\langle 3,0\rangle a\langle 1,0\rangle a \\
&=&[w]
\end{eqnarray*}
by Lemma \ref{l4.5}. Therefore, $(g,g)_{w_1}=0$. Similarly,
$(g,g)_{w_2}=0$.
\par Analogously, one can show that
$$
(h,g)_w \equiv 0 \ mod \ (q;w), \ \  (p,g)_w =0, \ \ (q,g)_w =0, \ \
(s,g)_w=0.
$$

8) $(f,h)_w\equiv0$. Take
$$
w=a\langle 0,0\rangle a\langle 0,1\rangle a\langle 0,1\rangle a.
$$
Expanding $(f,h)_w$, we get
$$
(f,g)_w=f\langle 0,1\rangle a\langle 0,1\rangle a-a\langle
0,0\rangle h=(a\langle 0,0\rangle a)\langle 0,1\rangle a\langle
0,1\rangle a-h-[w]=-h\equiv 0 \ mod \ (h).
$$
Next, we show that $(g,h)_w\equiv0$. Take
$$
w=a\langle 1,0\rangle a\langle 1,0\rangle a\langle 0,1\rangle
a\langle 0,1\rangle a.
$$
Expanding $(g,h)_{w}$, we get
$$
(g,h)_{w}=g\langle 0,1\rangle a\langle 0,1\rangle a-a\langle
1,0\rangle a\langle 1,0\rangle h =(a\langle 1,0\rangle a\langle
1,0\rangle a)\langle 0,1\rangle a\langle 0,1\rangle a-[w].
$$
Expanding the first monomial on the right hand side, we have
\begin{eqnarray*}
&&(a\langle 1,0\rangle a\langle 1,0\rangle a)\langle 0,1\rangle
a\langle 0,1\rangle a \\
&=&a\langle 1,0\rangle (a\langle 1,0\rangle a)\langle 0,1\rangle
a\langle 0,1\rangle a-a\langle 0,0\rangle (a\langle 1,0\rangle
a)\langle 1,1\rangle a\langle 0,1\rangle a \\
&=&[w]-a\langle 1,0\rangle a\langle 0,0\rangle a\langle 1,1\rangle
a\langle 0,1\rangle a \\
&&-a\langle 0,0\rangle a\langle 1,0\rangle a\langle 1,1\rangle
a\langle 0,1\rangle a+a\langle 0,0\rangle
a\langle 0,0\rangle a\langle 2,1\rangle a\langle 0,1\rangle a \\
&=&[w]-a\langle 1,0\rangle a\langle 0,0\rangle a\langle 1,1\rangle
a\langle 0,1\rangle a -a\langle 0,0\rangle a\langle 1,0\rangle
a\langle 1,1\rangle a\langle 0,1\rangle a
\end{eqnarray*}
by Lemma \ref{l4.5}. Therefore, $(g,h)_{w}\equiv 0$ mod$(q;w)$.
\par Analogously, one can show that
$$
(h,h)_w=0, \ \  (p,h)_w \equiv0 \ mod \ (q;w),
$$
$$
(q,h)_w\equiv0 \ mod \ (q;w) , \ \ (s,h)_w \equiv0 \ mod \ (q;w). \
\ \ \ \ \ \ \ \ \ \ \ \  \
$$

9) $(f,p)_w\equiv0$. Take
$$
w=a\langle 0,0\rangle a\langle 1,1\rangle a\langle 1,0\rangle a.
$$
Expanding $(f,p)_w$, we get
$$
(f,p)_w=f\langle 1,1\rangle a\langle 1,0\rangle a-a\langle
0,0\rangle p=(a\langle 0,0\rangle a)\langle 1,1\rangle a\langle
1,0\rangle a-p-[w]=-p\equiv 0 \ mod \ (p).
$$
Next, we show that $(g,p)_w=0$. Take
$$
w=a\langle 1,0\rangle a\langle 1,0\rangle a\langle 1,1\rangle
a\langle 1,0\rangle a.
$$
Expanding $(g,p)_{w}$, we get
$$
(g,p)_{w}=g\langle 1,1\rangle a\langle 1,0\rangle a-a\langle
1,0\rangle a\langle 1,0\rangle p =(a\langle 1,0\rangle a\langle
1,0\rangle a)\langle 1,1\rangle a\langle 1,0\rangle a-[w].
$$
Expanding the first monomial on the right hand side, we have
\begin{eqnarray*}
&&(a\langle 1,0\rangle a\langle 1,0\rangle a)\langle 1,1\rangle
a\langle 1,0\rangle a \\
&=&a\langle 1,0\rangle (a\langle 1,0\rangle a)\langle 1,1\rangle
a\langle 1,0\rangle a-a\langle 0,0\rangle (a\langle 1,0\rangle
a)\langle 2,1\rangle a\langle 1,0\rangle a \\
&=&[w]-a\langle 1,0\rangle a\langle 0,0\rangle a\langle 2,1\rangle
a\langle 1,0\rangle a \\
&&-a\langle 0,0\rangle a\langle 1,0\rangle a\langle 2,1\rangle
a\langle 1,0\rangle a+a\langle 0,0\rangle
a\langle 0,0\rangle a\langle 3,1\rangle a\langle 1,0\rangle \\
&=&[w]
\end{eqnarray*}
by Lemma \ref{l4.4}. Therefore, $(g,p)_{w}=0$.
\par Analogously, one can show that
$$
(h,p)_w\equiv0 \ mod \ (s;w), \ \  (p,p)_w =0, \ \ (q,p)_w=0 , \ \
(s,p)_w=0.
$$

10) $(f,q)_w\equiv0$. Take
$$
w=a\langle 0,0\rangle a\langle 1,1\rangle a\langle 0,1\rangle a.
$$
Expanding $(f,q)_w$, we get
$$
(f,q)_w=f\langle 1,1\rangle a\langle 0,1\rangle a-a\langle
0,0\rangle q=(a\langle 0,0\rangle a)\langle 1,1\rangle a\langle
0,1\rangle a-q-[w]=-q\equiv 0 \ mod \ (q).
$$
Next, we show that $(g,q)_w=0$. Take
$$
w=a\langle 1,0\rangle a\langle 1,0\rangle a\langle 1,1\rangle
a\langle 0,1\rangle a.
$$
Expanding $(g,q)_{w}$, we get
$$
(g,q)_{w}=g\langle 1,1\rangle a\langle 0,1\rangle a-a\langle
1,0\rangle a\langle 1,0\rangle q =(a\langle 1,0\rangle a\langle
1,0\rangle a)\langle 1,1\rangle a\langle 0,1\rangle a-[w].
$$
Expanding the first monomial on the right hand side, we have
\begin{eqnarray*}
&&(a\langle 1,0\rangle a\langle 1,0\rangle a)\langle 1,1\rangle
a\langle 0,1\rangle a \\
&=&a\langle 1,0\rangle (a\langle 1,0\rangle a)\langle 1,1\rangle
a\langle 0,1\rangle a-a\langle 0,0\rangle (a\langle 1,0\rangle
a)\langle 2,1\rangle a\langle 0,1\rangle a \\
&=&[w]-a\langle 1,0\rangle a\langle 0,0\rangle a\langle 2,1\rangle
a\langle 0,1\rangle a \\
&&-a\langle 0,0\rangle a\langle 1,0\rangle a\langle 2,1\rangle
a\langle 0,1\rangle a+a\langle 0,0\rangle a\langle 0,0\rangle
a\langle 3,1\rangle a\langle 0,1\rangle \\
&=&[w]-2a\langle 1,0\rangle a\langle 0,0\rangle a\langle 1,1\rangle
a\langle 1,1\rangle a-2a\langle 0,0\rangle a\langle 1,0\rangle
a\langle 1,1\rangle a\langle 1,1\rangle a
\end{eqnarray*}
by Lemma \ref{l4.4}. Therefore, $(g,q)_{w}\equiv0$ mod $(s;w)$.
\par Analogously, one can show that
$$
(h,q)_w=0, \ \  (p,q)_w \equiv0 \ mod \ (s;w),
$$
$$
(q,q)_w\equiv0 \ mod \ (s;w), \ \ (s,q)_w\equiv0 \ mod \ (s;w). \ \
\ \ \ \ \ \ \ \  \ \ \ \ \
$$

11) $(f,s)_w\equiv0$. Take
$$
w=a\langle 0,0\rangle a\langle 1,1\rangle a\langle 1,1\rangle a.
$$
Expanding $(f,s)_w$, we get
$$
(f,s)_w=f\langle 1,1\rangle a\langle 1,1\rangle a-a\langle
0,0\rangle s=(a\langle 0,0\rangle a)\langle 1,1\rangle a\langle
1,1\rangle a-s-[w]=-s\equiv 0 \ mod \ (s).
$$
Next, we show that $(g,s)_w=0$. Take
$$
w=a\langle 1,0\rangle a\langle 1,0\rangle a\langle 1,1\rangle
a\langle 1,1\rangle a.
$$
Expanding $(g,s)_{w}$, we get
$$
(g,s)_{w}=g\langle 1,1\rangle a\langle 1,1\rangle a-a\langle
1,0\rangle a\langle 1,0\rangle s =(a\langle 1,0\rangle a\langle
1,0\rangle a)\langle 1,1\rangle a\langle 1,1\rangle a-[w]=0
$$
by Lemma \ref{l4.3}.
\par Analogously, one can show that
$$
(h,s)_w=0, \ \  (p,s)_w =0, \ \ (q,s)_w=0 , \ \ (s,s)_w=0.
$$
This completes the proof of the theorem. \ \  $\square$

\ \

Let $S=\{f,g,h,p,q,s\}$. Let $[u]$ be an $S$-irreducible normal
word. Then [u] has the following form:
\begin{eqnarray*}
&&[u]=[a\langle 1,0\rangle a], \ \ [u]=[a\langle 0,1\rangle a], \ \
[u]=[a\langle 1,1\rangle a], \\
&&[u]=[a\langle 1,0\rangle a\langle 1,1\rangle a],\ \ [u]=[a\langle
0,1\rangle a\langle 1,1\rangle a], \\
&&[u]=[(a\langle 1,0\rangle a\langle 0,1\rangle)^{k} a],\ \
[u]=[(a\langle 0,1\rangle a\langle 1,0\rangle)^{l} a], \ k\geqslant
1, \ l\geqslant 1.
\end{eqnarray*}
By the Composition-Diamond lemma, these words consist of a linear
basis of $G$.

\subsection{Loop Lie  $n$-conformal algebra}
\par Let $\mathfrak{g}$ be a Lie algebra with a linear basis $\{a_i\}_{i\in
I}$. Then loop Lie $n$-conformal algebra for $\mathfrak{g}$ is given
by
$$
L(\mathfrak{g})=(\{a_i\}_{i\in I},\overrightarrow{N}=(1,1)| \
a_i\lfloor 0,0 \rfloor a_j=[a_ia_j], \ i>j, \ i,j \in I ).
$$
A universal enveloping associative $n$-conformal loop algebra of
$L(\mathfrak{g})$ is then given by
$$
\mathcal{U}_{\overrightarrow{N}}(L(\mathfrak{g}))=C(\{a_i\}_{i\in
I}, \overrightarrow{N}=(1,1)| \ a_i\langle 0,0 \rangle
a_j-\{a_i\langle 0,0 \rangle a_j\}=[a_ia_j], \ i>j, \ i,j \in I )
$$
which is an associative $n$-conformal loop algebra.
\begin{lemma}
The set $S=\{a_i\langle 0,0 \rangle a_j-\{a_i\langle 0,0 \rangle
a_j\}-[a_ia_j]| \ i>j \}\subseteq C(\{a_i\}_{i\in I},(1,1),D_1,D_2)$
is a Gr\"{o}bner-Shirshov basis.
\end{lemma}
\noindent{\bf Proof.} Any element of $S$ has the form
$$
a_i\langle 0,0 \rangle a_j-a_j\langle 0,0 \rangle a_i-[a_ia_j], \ \
i>j, \ i,j \in I
$$
and $\langle 0,0 \rangle$ is an associative multiplication. As the
same as for usual universal enveloping algebra of a Lie algebra,
these elements have trivial compositions of intersection. The
compositions of left multiplication are also trivial. As a result,
$S$ is a Gr\"{o}bner-Shirshov basis. \ \  $\square$

\ \


\begin{thebibliography}{4}

\bibitem{Aymon03} M. Aymon, P.-P. Grivel, Un Theoreme de Poincare-Birkhoff-Witt pour
les algebres de Leibniz, {\it Comm. Algebra}, \textbf{31}(2003),
527-544.

\bibitem{Bah87}Y.A Bahturin, Identical Relations in Lie Algebras,
VNU Science Press.x, 1987.


\bibitem{Bah03}Y.A. Bahturin, Groups, Rings, Lie and Hopf Algebras,
Kluwer Academis Publishers.ix, 2003.


\bibitem{BMPZ92}Y.A. Bahturin, A.A. Mikhalev, V.M. Petrogradsky, M.V. Zaicev,
Infinite-dimensional Lie Superalgebras,  de Gruyter Expositions in
Mathematics, 7. Walter de Gruyter \& Co., Berlin, 1992.


\bibitem{BAK01}B. Bakalov, A. D'Andrea, V.G. Kac, Theory of finite
pseudoalgebras, {\it Adv. Math.}, \textbf{162}(1)(2001), 1-140.

\bibitem{BD04}A.A Beilinson, V.G. Drinfeld, Chiral algebras,
{\it Amer. Math. Soc. Colloq. Publ.}, \textbf{51}(2004).

\bibitem{BPZZ84}A.A. Belavin, A.M. Polyakov, A.M. and A.B.
Zamolodchikov, Infinite conformal symmetry in two-dimensional
quantum field theory, {\it Nuclear Phys. B }, \textbf{214}(1984),
333-380.

\bibitem{Be78}G.M. Bergman, The Diamond Lemma for ring theory, {\it Adv. in
Math.}, \textbf{29}(1978), 178-218.

\bibitem{Bo07}L.A. Bokut, Gr\"obner--Shirshov bases for braid groups in Artin--Garside
generators, {\it J. Symbolic Compu.}, \textbf{43}(6-7)(2008),
397-405.

\bibitem{Bo09}L.A. Bokut, Gr\"obner--Shirshov bases for the braid groups in the
Birman-Ko-Lee generators, {\it Journal of Algebra},
\textbf{321}(2009), 361-376.

\bibitem{Bo72}L.A. Bokut,  Unsolvability of the word problem for Lie algebras,
{\it Sov. Math., Dokl.}, \textbf{13}(1972), 1388-1391.


\bibitem{Bo76}L.A. Bokut, Imbeddings into simple associative
algebras, {\it Algebra i Logika}, \textbf{15}(1976), 117-142.



\bibitem{BChainikovShum07}L.A. Bokut, V.V. Chainikov, K.P. Shum, Markov and Artin normal form
theorem for braid groups, {\it Comm. Algebra}, \textbf{35}(2007),
2105-2115.

\bibitem{BokutChen08}L.A. Bokut, Yuqun Chen, Gr\"{o}bner-Shirshov
 bases: Some new results, Proceedings of the Second International Congress
 in Algebra and Combinatorics, World Scientific, 2008, 35-56.

\bibitem{BCL08}L.A. Bokut, Yuqun Chen, Cihua Liu, Gr\"{o}bner-Shirshov
bases for dialgebras, arXiv.org/abs/0804.0638.

\bibitem{BCQ08}L.A. Bokut, Yuqun Chen, Jianjun Qiu,
Gr\"{o}bner-Shirshov bases for associative algebras with multiple
operators and free Rota-Baxter algebras, arXiv.org/abs/0805.0640.


\bibitem{BChibrikov05}L.A. Bokut, E. Chibrikov, Lyndon-Shirshov words, Gr\"obner-Shirshov
bases, and free Lie algebras, non-associative algebras and its
applications, { Taylor \& Francis Group, Boca Raton, FL}, 2005,
17-34.

\bibitem{BFK04}L.A. Bokut, Y. Fong, W.-F. Ke, Composition Diamond Lemma for
 associative conformal algebras, {\it Journal of Algebra},  \textbf{272}(2004), 739-774.


\bibitem{BoFKK00}L. A. Bokut, Y. Fong, W.-F. Ke, P. Kolesnikov, Gr\"obner and
Gr\"obner-Shirshov bases in Algebra and Conformal algebras, {\it
Fundamental and Applied Mathematics}, \textbf{6}(2000), N3, 669-706
(in Russian).

\bibitem{BFKS08}L.A. Bokut, Y. Fong, W.-F. Ke, L.-S. Shiao, Gr\"{o}bner-Shirshov
 bases for the braid semigroup, Proceedings of the ICM satellite conference in
 algebra and related topics, Hong Kong, China, August 14--17, 2002.
 River Edge, NJ: World Scientific. 60-72 (2003).

\bibitem{BoKLM99}L.A. Bokut, S.-J. Kang, K.-H. Lee, P. Malcolmson, Gr\"obner--Shirshov bases
for Lie supalgebras and their universival enveloping algebras, {\it
Journal of Algebra}, \textbf{217}(2)(1999), 461-495.

\bibitem{BK00}L.A. Bokut, P.S. Kolesnikov, Gr\"obner---Shirshov bases: From Incipient
to Nowdays, {\it Proceedings of the POMI}, \textbf{272}(2000),
26-67.

\bibitem{BK05a}L.A. Bokut, P.S. Kolesnikov, Gr\"obner-Shirshov bases:
conformal algebras and pseudoalgebras, {\it Journal of Mathematicfal
Sciences}, \textbf{131}(2005), No.5, 5962-6003.


\bibitem{BokutShiao01}L.A. Bokut, L.-S. Shiao, Gr\"obner-Shirshov bases for coxeter groups,
{\it Comm.  Algebra}, \textbf{29}(2001), No.9, 4305-4319.

\bibitem{Bor86}R.E. Borcherds, Vertex algebras, Kac-Moody algebras,
and the monster, {\it  Proc. Nat. Acda. Sci. U.S.A.}, \textbf
{83}(1986), 3068-3071.

\bibitem{Bu65}B. Buchberger, An algorithm for finding a basis for the
residue class ring of a zero-dimensional polynomial ideal, Ph.D.
thesis, University of Innsbruck, Austria, (1965). (in German)

\bibitem{Bu70}B. Buchberger, An algorithmical criteria for the
solvability of algebraic systems of equations, {\it Aequationes
Math.}, \textbf{4}(1970), 374-383. (in German)

\bibitem{CFL58}K.-T. Chen, R.H. Fox, R.C. Lyndon, Free differential caculus. IV.
The quotient groups of the lower central series, {\it Ann. Math.},
\textbf{68}(1958), 81-95.

\bibitem{Cohn65}P.M. Cohn, Universal Algebra, Harper \& Row, Publishers, New
York-London, 1965.

\bibitem{DH08}
V. Drensky, R. Holtkamp,  Planar trees, free nonassociative
 algebras, invariants,  and elliptic integrals, {\it Algebra and Discrete
 Mathmatics}, {\bf 2}(2008), 1-41.

\bibitem{EGuo08}K. Ebrahimi-Fard, L. Guo, Free Rota--Baxter algebras and rooted
trees, {\it  J. Algebra Applications}, {\bf7}(2008), 167-194.

\bibitem{Hi64}H. Hironaka, Resolution of singulatities of an algebtaic variety
over a field if characteristic zero, I, II, {\it Ann. Math.},
\textbf{79}(1964), 109-203, 205-326.

\bibitem{Kac96}V. Kac, Vertex algebras for beginners, University
Lecture Series, Vol. 10., AMS, Providence, RI, (1996).

\bibitem{KL00a}S.-J. Kang, K.-H. Lee, Gr\"obner-Shirshov bases for irreducible
$sl_{n+1}$-modules, {\it Journal of Algebra}, \textbf{232}(2000),
1-20.

\bibitem{KL00b}S.-J. Kang, K.-H. Lee, Gr\"obner-Shirshov bases for
representation theory, {\it J. Korean Math. Soc.},
\textbf{37}(2000), 55-72.

\bibitem{KLLO02}S.-J. Kang, I.-S. Lee, K.-H. Lee, H. Oh, Hecke algebras, specht
modules and Gr\"obner-Shirshov bases, {\it Journal of Algebra} ,
\textbf{252}(2002), 258-292.

\bibitem{KLLP07}S.-J. Kang, D.-I. Lee, K.-H. Lee, H. Park, Linear algebraic approach to
  Gr\"obner-Shirshov basis theory, {\it Journal of Algebra},
  \textbf{313}(2007), 988-1004.

\bibitem{Ko91}P.S. Kolesnikov, On irreducible algebras of conformal
endomorphisms over a linear algebraic group, {\it Mathematics
Subject Classification}, 1991.(arXiv.org/abs/0712-4127)

\bibitem{Ko00a}P.S. Kolesnikov, On commutative conformal algebras,
master's thesis, Novosibirsk State University, 2000.

\bibitem{Ko00b}P.S. Kolesnikov, A base for a free associative
commutative conformal algebra, {\it The 4th Siberian Congress in
Applied and Industrial Mathematics}, Abstracts, Part IV, 2000.


\bibitem{Ko02}P.S. Kolesnikov, Universal representations of some Lie conformal supalgebras,
  vestnik, {\it Quart. J. of Novosibirsk State Univ., Series: math., mech. and
  informatics}, \textbf{2}(3)(2002),  30-45.

\bibitem{Ko03}P.S. Kolesnikov, Associative enevloping
pseudoalgebras of finite Lie pseudoalgebras, {\it Comm. Algebra},
\textbf{31}(2003), 2909-2925.

\bibitem{Ko07a}P.S. Kolesnikov, Associative algebras related to conformal
algebras, {\it Appl. Categ. Struct.}, \textbf{16}(2008), No. 1-2,
167-181. MSC2000.

\bibitem{Ko07b}P.S. Kolesnikov, Conformal representations of Leibniz
algebras, {\it Sib. Mat. Zh.}, \textbf{49}(3)(2008), 540-547;
translation in {\it Sib. Math. J.}, \textbf{49}(3)(2008),  429-435.

\bibitem{Ko08a}P.S. Kolesnikov, Variethies of dialgebras and conformal
algebras, {\it Siberian Mathematical Journal}, Springer New York,
\textbf{49}(2008), 257-272.


\bibitem{Ko08b}P.S. Kolesnikov, Universally defined representations of
Lie conformal superalgebras, {\it Journal of Symbolic Computation},
\textbf{43}(2008), 406-421.

\bibitem{Ku60}A.G. Kurosh, Free sums of multiple operator algebras,
{\it Siberian. Math. J.}, {\bf 1}(1960), 62-70.

\bibitem{Lyndon54}R.C. Lyndon, On Burnside Problem, {\it Trans. Amer. Math.
Soc.}, \textbf{77}(1954), 202-215.

\bibitem{Mik96}
A.A. Mikhalev, Shirshov's composition techniques in Lie
superalgebras (Non-Commutative Gr\"obner Bases), {\it Trudy Sem.
Petrovsk}, \textbf{18}(1995), 277-289 (in Russian). English
translation: {\it J.~Math. Sci.}, \textbf{80}(1996), 2153-2160.


\bibitem{Reu93}C. Reutenauer, Free Lie algebras, Oxford Science Publications, 1993.

\bibitem{Ro99}M. Roitman, On free conformal and vertex algebras, {\it
Journal of Algebra}, \textbf{217}(2)(1999), 496-527.

\bibitem{Ro00}M. Roitman, Universal enveloping conformal algebras,
{\it Selecta Math. (N.S)}, \textbf{6}(3)(2000), 319-345.

\bibitem{Ro05}M. Roitman, A criterion for embedding of Lie conformal algebras into
associative conformal algebras, {\it Journal of Lie theory},
\textbf{15}(2)(2005), 575-588.

\bibitem{Sh53a}A.I. Shirshov, Subalgebras of free Lie algebras,
 {\it Mat. Sbornik N. S.}, \textbf{33}(2)(75)(1953),  441-452.

\bibitem{Sh53b}A.I. Shirshov, On the representation of Lie rings in associative
rings, {\it Uspekhi Mat. Nauk}, \textbf{8}(1953), 173--175.

\bibitem{Sh54}A.I. Shirshov, Subalgebras of free commutative and free
anti-commutative algebras,  {\it Mat. Sbornik.}, {\bf34}(76)(1954),
81-88.

\bibitem{Sh58}A.I. Shirshov, On free Lie rings, {\it Mat. Sb.}, \textbf{45}(87)(1958),
113-122. (in Russian)

\bibitem{Sh62a}A.I. Shirshov, On the bases of a free Lie algebras,
 {\it Algebra Logika}, \textbf{1}(1)(1962),  14-19. (in Russian)

\bibitem{Sh62b}A.I. Shirshov, Some algorithmic problem for $\varepsilon$-algebras,
 {\it Sibirsk. Mat. Z.}, \textbf{3}(1962), 132-137. (in Russian)

\bibitem{Sh62c}A.I. Shirshov, Some algorithmic problem for Lie algebras,
 {\it Sibirsk. Mat. Z.}, \textbf{3}(1962), 292-296 (in Russian).
 English translation: {\it SIGSAM Bull.}, \textbf{33}(2)(1999), 3-6.


\bibitem{Ufn95}V.A. Ufnarovskij, Combinatorial and asymptotic methods in
algebra. VI. 1 196. (Encyclopacdia Math. Sci., 57). Springer,
Berlin, 1995.

\bibitem{Vie78}G. Viennot, Alg$\grave{e}$bres de Lie Libres et Mono$\ddot{i}$des Libres,
Lecture Notes in Mathematics, Vol. 691, Springer, Berlin, 1978.

\bibitem{Zhu50}A.I. Zhukov, Complete systems of defining relations in non-associative algebras,
 {\it Mat. Sbornik}, \textbf{69}(27)(1950), 267-280.






\end{thebibliography}
\end{document}